\documentclass{siamltex}

\usepackage{amsmath}
\usepackage{amsfonts}
\usepackage{amssymb}
\usepackage{mathrsfs}
\usepackage{graphicx}
\usepackage{latexsym, amscd}

\usepackage[refpage,noprefix]{nomencl}

\makeglossary

\newtheorem{sideeg}[theorem]{Example} 
\newenvironment{example}{\begin{sideeg}\rm}{\end{sideeg}}
\newtheorem{sideremark}[theorem]{Remark}
\newenvironment{remark}{\begin{sideremark}\rm}{\end{sideremark}}
\newtheorem{sideruleofthumb}[theorem]{Rule of thumb}

\title{A discrete invitation to quantum filtering and feedback control\thanks{
R.v.H.\ and L.B.\ are partially supported by the Army Research Office 
under Grants DAAD19-03-1-0073 and W911NF-06-1-0378, and by the National 
Science Foundation under Grants CCF-0323542 and PHY-0456720.  M.R.J.\ is 
supported by the Australian Research Council.}}

\author{Luc Bouten\footnotemark[2] \and Ramon van Handel\footnotemark[2]
\and Matthew R. James\footnotemark[3]}

\newcommand{\BB}[1]{\mathbb{#1}}
\newcommand{\1}{\hbox{$1$\hspace{-.7pt}%
\vrule width .5pt height 2.5mm depth 0mm\hspace{-0.8mm}%
\vrule width 1.2mm height 0.08mm depth 0mm\hspace{-1.2mm}%
\vrule width 0.7mm height 2.5mm depth -2.42mm}\hspace{.9mm}}

\newcommand{\ten}{\otimes}

\begin{document}
\maketitle

\renewcommand{\thefootnote}{\fnsymbol{footnote}}
\footnotetext[2]{
        Physical Measurement and Control 266-33, California Institute of
        Technology, Pasadena, CA 91125, USA (bouten@its.caltech.edu,
        ramon@its.caltech.edu).
}
\footnotetext[3]{
        Department of Engineering, Australian National University,
        Canberra, ACT 0200, Australia (matthew.james@anu.edu.au).
}
\renewcommand{\thefootnote}{\arabic{footnote}}

\begin{abstract} 
The engineering and control of devices at the quantum-mechanical
level---such as those consisting of small numbers of atoms and
photons---is a delicate business.  The fundamental uncertainty that is
inherently present at this scale manifests itself in the unavoidable
presence of noise, making this a novel field of application for stochastic
estimation and control theory.  In this expository paper we demonstrate
estimation and feedback control of quantum mechanical systems in what is
essentially a noncommutative version of the binomial model that is
popular in mathematical finance.  The model is extremely rich and allows a
full development of the theory, while remaining completely within the
setting of finite-dimensional Hilbert spaces (thus avoiding the technical 
complications of the continuous theory).  We introduce discretized models 
of an atom in interaction with the electromagnetic field, obtain filtering 
equations for photon counting and homodyne detection, and solve a 
stochastic control problem using dynamic programming and Lyapunov function
methods.
\end{abstract}

\begin{keywords}
discrete quantum filtering, quantum feedback control, 
quantum probability, conditional expectation, dynamic programming, 
stochastic Lyapunov functions.
\end{keywords}

\begin{AMS}
        93E11,  93E15,  93E20,  81P15,  81S25,  34F05.
\end{AMS}

\pagestyle{myheadings}
\thispagestyle{plain}
\markboth{BOUTEN, VAN HANDEL, AND JAMES}{A DISCRETE INVITATION TO QUANTUM 
FILTERING AND FEEDBACK CONTROL}


\vskip.2cm
\begin{center}
{\it Dedicated to Slava Belavkin in the year of his 60th birthday.}
\end{center}
\vskip.4cm

\section{Introduction}\label{sec intro}

Control theory, and in particular feedback control, is an important aspect
of modern engineering.  It provides a set of tools for the design of
technology with reliable performance and has been applied with
overwhelming success in the design of many of the devices that we use on a
daily basis.  In this article we will explore the following question:  
rather than controlling, say, a jet engine, can we use feedback to control
an object as small as a single atom?

Though it is not directly visible in our everyday lives, the technology to
manipulate matter at the level of single atoms and photons is well in
place and is progressing at a fast rate.  Nobel prize winning technologies
such as laser cooling and trapping of atoms, which were state-of-the-art
only a decade ago, are now routine procedures which are being applied in
physics laboratories around the world.  The generation of highly coherent
light and the detection of single photons has been refined since the
development of the laser.  Given the availability of this unprecedented
level of control over the microscopic world, the question of control
design seems a timely one---one could wonder whether we can ``close the
loop'' and make the atoms do useful work for us.  This is the domain of
quantum feedback control.

``Quantum,'' of course, refers to the theory of quantum mechanics, which
comes in necessarily at this level of description.  The feedback which we
will use to control a quantum system is a function of observations
obtained from that system, and observations in quantum mechanics are
inherently random.  This makes the theory both challenging and interesting
from the point of view of fundamental physics and control theory.  As we
will demonstrate, the theory of quantum feedback control resembles
closely, and follows directly in many parts, the classical theory of
stochastic control (here and below, {\it classical} means non-quantum
mechanical).

\begin{figure}
\centering
\includegraphics[width=\textwidth]{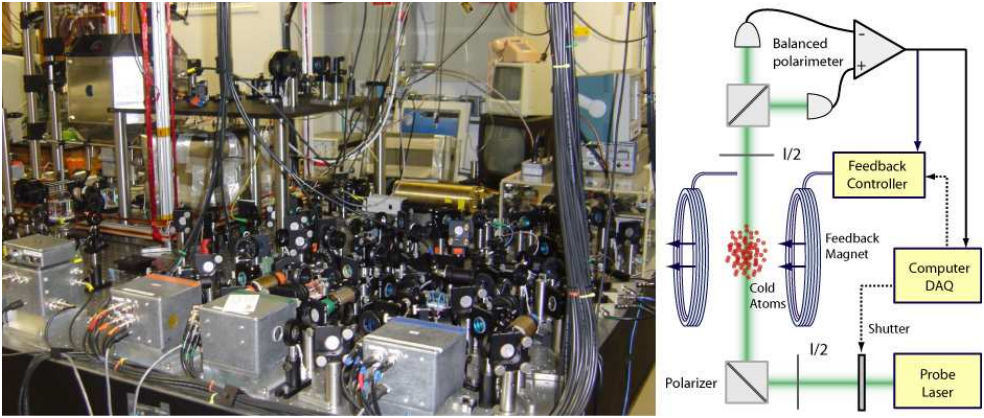}
\caption{\label{fig:photo} A laboratory experiment in quantum feedback 
control, implementing the setup reported in {\rm \cite{GSM05}}.  The metal 
boxes in the foreground are diode lasers, whose output light is 
manipulated using a large number of optical elements.  In the background 
on the left is a magnetically shielded chamber containing magnetic coils 
and a cell with laser cooled Cesium atoms.  A schematic of the experiment 
is shown on the right. (Photograph courtesy of John Stockton.)}
\end{figure}

Several promising applications of quantum feedback control have been
proposed and are now starting to be investigated in a laboratory setting.  
One class of such applications falls under the general heading of
precision metrology: can we utilize the sensitivity of small numbers of
atoms to external perturbations to design ultrafast, precise sensors,
when the desired accuracy is on the order of the intrinsic quantum
uncertainty of the sensor?  Concrete applications include e.g.\ precision
magnetometry \cite{GSDM03,SGDM04}, which was recently investigated in a
laboratory experiment \cite{GSM05} (see Figure \ref{fig:photo}), and
atomic clocks \cite{ASL04}.  A second class of examples involves the
design of optical communication systems where each bit is encoded in a
pulse containing only a small number of photons \cite{Dav77,Ger04}; see
\cite{AASDM02} for a related experimental demonstration.  As a third
application, we mention the use of feedback for cooling of atoms or
nanomechanical devices \cite{DJ99,Asa03,Steck04}, see \cite{ExpCool06} for
a laboratory demonstration (the availability of efficient cooling
techniques is important in experimental physics).

As both the theory and applications are still in their infancy, it is far
from clear at this time what the ultimate impact of quantum feedback
control will be; it is clearly unlikely that future generations of washing
machines or airplanes will be powered by quantum control, but other
application areas, such as the ones mentioned above, might well benefit
from the availability of such technology.  At the very least, however, the
questions that arise in the subject prompt us to explore the fundamental
limits on the control and engineering of physical systems, and a fruitful
interplay between control theory and experimental physics enables the
exploration of these limits using real-life experiments such as the one
displayed in Figure \ref{fig:photo}. In this article we demonstrate the 
underlying theory using a simple but rich toy model; we hope that the 
reader is sufficiently intrigued to come along for the ride!

The study of quantum feedback control was pioneered by V.\,P. Belavkin in
a 1983 paper \cite{Bel83}, long before the experimental realization of
such a setup was realistically feasible.  Belavkin realized that due to
the unavoidable presence of uncertainty in quantum measurements, the
theory of filtering---the extraction of information from noisy
observations---plays a fundamental role in quantum feedback control,
exactly as in the theory of stochastic control with partial observations
\cite{Bensoussan}.  In \cite{Bel83} feedback control was explored in
discrete time using the operational formalism of Davies \cite{Dav76},
which is the precursor of quantum filtering theory.  Belavkin developed
further the theory of quantum filtering in a series of articles
\cite{Bel80,Bel88,Bel92b}, and in \cite{Bel88} the use of dynamic
programming for feedback control in continuous time was sketched. We refer
to \cite{BHJ06} for a modern introduction to quantum filtering theory.

Independently from Belavkin's work, the equations of quantum nonlinear
filtering were developed in the physics literature by Carmichael
\cite{Car93}, based on the work of Davies \cite{Dav76}, under the name of
``quantum trajectory theory''.  Though the connection to filtering theory
was not realized until much later, Carmichael's work nonetheless spawned
an investigation on the use of feedback control in the physics literature
(see e.g.\ \cite{DHJMT00} and the references therein).  It is here that 
most of the practical applications, such as the ones cited above, were 
proposed and developed.

Like its classical counterpart, the mathematical theory of quantum
filtering and stochastic control in continuous time can be quite
technical, requiring a heavy analytic machinery to obtain the
necessary results in a rigorous way.  On the other hand, we believe that
the methods in the field can be understood perfectly well without
being obscured by analytic complications.  What we are going to do in
this paper is to give a complete development of the basic theory
of quantum filtering and feedback control in a simple toy model, which
requires essentially no analytic machinery.  We only assume some 
background in linear algebra and familiarity with elementary probability 
theory with martingales, roughly at the level of the inspiring textbook by 
D.\ Williams \cite{Wil91}.

The model we will investigate likely has a familiar ring to those readers
who are familiar with mathematical finance.  It is in fact little more
than a noncommutative version of the binomial model which has enjoyed wide
popularity since its introduction by Cox, Ross and Rubinstein
\cite{CRR79}.  The model is widely used in the teaching of mathematical
finance \cite{Shr04} and even on Wall Street as a flexible computational
tool for the pricing of financial derivatives.  The model has three key
features, all of which are shared by its noncommutative counterpart, which
make its investigation particularly appealing: (i) on the one hand the
model is as simple as it can possibly be, consisting of a finite number of
two-state random variables (coin flips)---hence analytical complications
are avoided;  (ii) on the other hand, many of the mathematical tools that
are indispensable in the continuous context can be developed within the
framework of this model and admit simple proofs: in particular, change of
measure techniques, the martingale representation theorem, discrete
stochastic integrals and even a (trivial) stochastic calculus can all be
demonstrated in this context; and (iii) the model is in fact an
approximation to the usual continuous time models in the sense that these
are recovered when the limit as the time step tends to zero is taken in a
suitable sense.  For these reasons, we believe that such models are
ideally suited as a first introduction to quantum filtering and feedback
control.

The noncommutative model has an independent history in the quantum
probability literature.  In his book \cite{Mey93} P.-A.\ Meyer uses the
building blocks of the model we use in this article to demonstrate the
basic features of the Fock space, which underlies the continuous quantum
stochastic calculus of Hudson-Parthasarathy \cite{HuP84}, in a discrete
setting (he calls this model, which he attributes to Journ{\'e} \cite[page
18]{Mey93}, the ``toy Fock space'').  Repeated interaction models of a
similar type were considered by K{\"u}mmerer \cite{Kum85}.  Various 
authors have considered the continuous limit of such models
\cite{LP88,Att03,AtP06,Gough04} and have demonstrated convergence of the
associated discretized evolution equations to (Hudson-Parthasarathy)
quantum stochastic differential equations which are routinely used for the
description of realistic physical systems \cite{GaZ04}.  The description
of measurements (in essence, filtering)  in this discrete framework was
deduced from basic physical arguments in Brun \cite{Brun02} and Gough
\cite{GS04}.

In this paper we will not take the simplest and most straightforward route
in treating these discrete models; similar equations can often be deduced
by more elementary means \cite{Brun02,GS04}.  Rather, we will exploit the
richness of the model to mimic as closely as possible the development of
quantum filtering and feedback control in the continuous setting, where
simple arguments based on linear algebra are useless.  Indeed, the
introduction of a ``trivial'' stochastic calculus, martingale
representation, etc., may seem overkill in the context of this model, but
these become powerful and indispensable tools in the treatment of the more
general (and physically relevant) framework of continuous time evolution
and observations.  We hope to provide in this way a (relatively) gentle
introduction to the full technical theory, where all the necessary
ingredients make their appearance in a more accessible setting.

The remainder of the paper is organized as follows.  In \S\ref{sec random
variables} we introduce the basic notions of quantum probability, and
develop in detail the quantum binomial model which we will be studying
throughout the article.  \S\ref{conditional expectation} introduces
conditional expectations in the quantum context and sets up the filtering
problem.  In \S\ref{discrete QSC} we develop a discrete stochastic
calculus and use it to express the model of \S\ref{sec random variables}
in the form of a difference equation.  \S\ref{sec doob} treats filtering
using martingale methods, while \S\ref{reference probability} uses change
of measure techniques.  In \S\ref{sec:fb} we show how to incorporate
feedback into the model of \S\ref{sec random variables}. In
\S\ref{sec:optimal} we develop optimal feedback controls using dynamic
programming techniques, and \S\ref{sec:lyapunov} briefly explores an
alternative control design method using Lyapunov functions.  We conclude
in \S\ref{continuous time} with a set of references for further reading.


\renewcommand{\nompreamble}{ We use the following conventions.  In
general, Hilbert spaces are denoted by the symbol $\mathsf{H}$ and are
endowed with the inner product $\langle\cdot,\cdot\rangle$.  We will
assume all Hilbert spaces are complex and finite-dimensional. Recall that
the adjoint of a linear operator $X:\mathsf{H}\to\mathsf{H}$ is defined as
the unique operator $X^*$ that satisfies $\langle X^*x,y\rangle=\langle
x,Xy\rangle$ $\forall x,y\in\mathsf{H}$, and that $X$ is called
self-adjoint if $X=X^*$.  We denote by $I$ the identity operator.  The
commutator between two linear operators is denoted as $[X,Y]=XY-YX$.  
Calligraphic letters will be used for various purposes, such as classical
$\sigma$-algebras (e.g.\ $\mathcal{Y}$) or functionals on an algebra
(e.g.\ $\mathcal{L}(X)$).  Control sets and related quantities are denoted
by gothic type (e.g.\ $\mathfrak{U}$).  $^*$-algebras (to be defined
below)  will be denoted by script type $\mathscr{A}$, and states on such
an algebra are often denoted by blackboard type (e.g.\ $\mathbb{P}$),
though $\rho$ and $\phi$ will also be used. The generic classical
probability space will be denoted as $(\Omega,\mathcal{F},\mathbf{P})$,
$\mathbf{E_P}$ denotes the expectation with respect to the measure
$\mathbf{P}$, and $\ell^\infty(\mathcal{F})$ denotes the space of
$\mathcal{F}$-measurable (bounded, but trivially so when $\Omega$ is a
finite set) random variables.

As is unavoidable in an article of this length, there is much notation
that is introduced throughout the article and that is constantly reused.  
To help the reader to keep track of the various quantities, we have
provided below a list of commonly used objects, a brief description, and 
the page number on which the object is defined. \hfill\break}

\printglossary[2cm]

\section{The quantum binomial model}\label{sec random variables}

In this section we introduce the basic model that we will be dealing with
throughout: a discretized approximation of the interaction
between an atom and the electromagnetic field.  First, however, we need to 
demonstrate how probability theory fits in the framework of quantum 
mechanics.

\subsection{Random variables in quantum mechanics}\label{sec:rvini}

The basic setting of quantum mechanics, as one would find it in most
textbooks, is something like this.  We start with a Hilbert space
$\mathsf{H}$ and fix some ``state vector'' $\psi\in\mathsf{H}$.  An
``observable'', the physicist's word for random variable, is described by
a self-adjoint operator $X$ on $\mathsf{H}$, and the expectation of $X$
is given by $\langle\psi,X\psi\rangle$.  The set of values that $X$ can
take in a single measurement is its set of eigenvalues, and the
probability of observing the eigenvalue $\lambda_i$ is given by
$\langle\psi,P_i\psi\rangle$ where $P_i$ is the projection operator onto
the eigenspace corresponding to $\lambda_i$.  This is quite unlike the
sort of description we are used to from classical probability theory---or
is it?

In fact, the two theories are not as dissimilar as they may seem, and it
is fruitful to formalize this idea (we will do this in the next section).  
The key result that we need is the following elementary fact from linear
algebra.  This is just the statement that a set of commuting normal
($[X,X^*]=0$) matrices can be simultaneously diagonalized, see e.g.\
\cite[sec.\ 2.5]{HJ85} or \cite[sec.\ 84]{PRH74}.

\begin{theorem}[\bf Spectral theorem]\label{thm spectral theorem}
Let $\mathsf{H}$ be an $n$-dimensional Hilbert space, $n<\infty$. Let
$\mathscr{C}$ be a set of linear transformations from $\mathsf{H} \to
\mathsf{H}$ that is closed under the adjoint {\rm (}i.e.\ if
$C\in\mathscr{C}$ then also $C^*\in\mathscr{C}${\rm )} and such that all
the elements of $\mathscr{C}$ commute {\rm (}i.e.\ $[C,D]=0$ $\forall
C,D\in\mathscr{C}${\rm )}.  Then there exists an orthonormal basis of
$\mathsf{H}$ such that every $C\in\mathscr{C}$ is represented by a
diagonal matrix with respect to this basis.
\end{theorem} 

Let us demonstrate how this works in the simplest case.  Let 
${\rm dim}\,\mathsf{H}=2$, fix some $\psi\in\mathsf{H}$ and let $X=X^*$ 
be a self-adjoint operator on $\mathsf{H}$.  The set $\mathscr{C}=\{X\}$ 
satisfies the conditions of the spectral theorem, so we can find an 
orthonormal basis in $\mathsf{H}$ such that we can express $X$ and 
$\psi$ in this basis as
\begin{equation*}
	X = \begin{pmatrix} x_1 & 0 \\ 0 & x_2\end{pmatrix},
	\qquad
	\psi = \begin{pmatrix}
		\psi_1 \\ \psi_2
	\end{pmatrix}.
\end{equation*}
We can now interpret $X$ as a random variable on some probability space.  
Introduce $\Omega=\{1,2\}$, the map $x:\Omega\to\mathbb{R}$, $x(i)=x_i$,
and the probability measure $\mathbf{P}(\{i\})=|\psi_i|^2$ ($i=1,2$).
Evidently $\langle\psi,X\psi\rangle = \mathbf{E_P}(x)$, i.e.\ the random
variable $x$ has the same expectation under $\mathbf{P}$ as we obtain from
$X$ and $\psi$ using the usual quantum mechanical formalism.  We can also
easily calculate $\mathbf{P}(x=x_i)=|\psi_i|^2$, which is again consistent 
with the rules of quantum mechanics.  As the spectral theorem works for 
more general sets $\mathscr{C}$, we can follow the same procedure to 
represent a set of commuting observables on a classical probability space 
and to calculate the joint statistics.

Up to this point we have done nothing particularly exciting: all we have
done is to express the quantum ``rules'' listed in the first paragraph of
this section in a convenient basis, and we have attached an interpretation
in terms of classical probability theory.  Conceptually, however, this is
an important point: commuting observables {\it are} random variables on a
classical probability space, and should be thought of in that way.  
Formalizing the way we can pass back and forth between the two
descriptions will allow us to utilize classical probabilistic techniques
in the quantum mechanical setting.

What the spectral theorem does not allow us to do is to simultaneously
interpret two noncommuting observables as random variables on a single
probability space.  This is not a shortcoming of the theory, but has an
important physical meaning. Observables that do not commute cannot be
observed in the same realization: their joint probability distribution is
a meaningless quantity (hence they are called {\it incompatible}).  
Strange as this may seem, this idea is a cornerstone of quantum theory and
all the empirical evidence supports the conclusion that this is how nature
works.  We will accept as a physical principle that in a single
realization we can choose to measure at most a commuting set of
observables.  However, we will see later that we can still estimate the
statistics of observables which we did not choose to measure, even if
these do not commute with each other, provided they commute with the
measurements we did choose to perform.  These ideas will be clarified in
due course.

\subsection{Quantum probability spaces}

In this section we are going to introduce and motivate the notion of a
(finite-dimensional) quantum probability space.  The definitions may not
seem completely natural at first sight, but their use will be amply
illustrated in the remainder of the article.

In the example of the previous section, we saw how to construct a set of
quantum observables as maps on the sample space $\Omega$.  To turn this
into a probability space we need to add a $\sigma$-algebra and a
probability measure.  As we are considering a particular set of random
variables $X_1,\ldots,X_k$ it is natural to choose the $\sigma$-algebra
$\sigma(X_1,\ldots,X_k)$, i.e.\ the smallest $\sigma$-algebra with respect
to which $X_1,\ldots,X_k$ are measurable (the $\sigma$-algebra generated
by $X_1,\ldots,X_k$), and the quantum state induces a probability measure
on this $\sigma$-algebra.  We would like to be able to specify the
$\sigma$-algebra directly, however, without having to select some
``preferred'' set of random variables $X_1,\ldots,X_k$---the values taken 
by these random variables are irrelevant, as we are only interested in the 
generated $\sigma$-algebra.

The definition we will give is motivated by the following observation. On
a classical sample space $\Omega$, consider the set $\ell^\infty(\mathcal{F})$
of all random variables measurable with respect to the $\sigma$-algebra
$\mathcal{F}$.  If one is given such a set, then $\mathcal{F}$ can be reconstructed
as $\mathcal{F}=\sigma\{\ell^\infty(\mathcal{F})\}$.  Conversely, however, we can
easily construct $\ell^\infty(\mathcal{F})$ if we are given $\mathcal{F}$.  Hence
these two structures carry the same information; by considering all
measurable random variables at the same time we are no longer giving
preference to a particular set $X_1,\ldots,X_k$.  In quantum probability,
it will be more convenient to characterize a $\sigma$-algebra by its
$\ell^\infty$-space---the quantum analog of this object is a set of
observables $\mathscr{A}$, which can be used directly as input for the
spectral theorem.  To complete the picture we need a way to characterize
sets of observables $\mathscr{A}$ which are mapped by the spectral theorem
to a space of the form $\ell^\infty(\mathcal{F})$.

\begin{definition}
A {\em $^*$-algebra} $\mathscr{A}$ is a set of linear transformations
$\mathsf{H}\to\mathsf{H}$ such that $I,\alpha A+\beta 
B,AB,A^*\in\mathscr{A}$ for any $A,B\in\mathscr{A}$, 
$\alpha,\beta\in\mathbb{C}$.  $\mathscr{A}$ is called {\em commutative} 
if $AB=BA$ for any $A,B\in\mathscr{A}$.  A linear map 
$\rho:\mathscr{A}\to\mathbb{C}$ which is positive $\rho(A^*A)\ge 0$ 
$\forall A\in\mathscr{A}$ and normalized $\rho(I)=1$ is called a {\em 
state} on $\mathscr{A}$.
\end{definition}

\begin{definition}
A $^*$-isomorphism between a commutative $^*$-algebra $\mathscr{A}$ and a 
set of functions $\mathcal{A}$ on some space $\Omega$ is a linear 
bijection $\iota:\mathscr{A}\to\mathcal{A}$ such that 
$\iota(A^*)=\iota(A)^*$ {\rm (}$\iota(A)^*(i)$ is the complex conjugate of 
$\iota(A)(i)$ $\forall i\in\Omega${\rm )} and $\iota(AB)=\iota(A)\iota(B)$ 
{\rm (}$(\iota(A)\iota(B))(i)=\iota(A)(i)\iota(B)(i)$ $\forall 
i\in\Omega${\rm )} for every $A,B\in\mathscr{A}$.
\end{definition}

The reason we want $^*$-isomorphisms is to ensure that we can manipulate
observables and classical random variables in the same way.  That is, if
$X_1,X_2$ are commuting self-adjoint operators that correspond to the
random variables $x_1,x_2$, then $X_1+X_2$ must correspond to $x_1+x_2$
and $X_1X_2$ must correspond to $x_1x_2$.  The notion of a $^*$-algebra
implements exactly the question posed above: every commutative
$^*$-algebra can be mapped to a set of random variables of the form
$\ell^\infty(\mathcal{F})$ for some $\sigma$-algebra $\mathcal{F}$ by using the
spectral theorem and the following Lemma.  The fact that sets of 
measurable functions can be characterized by algebras is well-known;
in fact, Lemma \ref{lem:onetoonespectral} is just a ``trivial'' version 
of the monotone class theorem \cite[Ch.\ I, Thm.\ 8]{Pro04}.

\begin{lemma}
\label{lem:onetoonespectral}
Let $\Omega=\{1,\ldots,n\}$.  Then there is a one-to-one correspondence
between (commutative) $^*$-algebras of $n\times n$ diagonal matrices and 
$\ell^\infty(\mathcal{F})$-spaces.
\end{lemma}

\begin{proof} 
We need two simple facts.  First, let $X_1,\ldots,X_k$ be elements of a 
$^*$-algebra of diagonal matrices.  Then $f(X_1,\ldots,X_k)$ is also an 
element of the $^*$-algebra for any function $f$.  Clearly this is true if 
$f$ is any (complex) polynomial.  But this is sufficient, because given 
a finite number of points $x_1,\ldots,x_n\in\mathbb{C}^k$ we can always 
find for any function $f$ a polynomial $\hat f$ that coincides with $f$ on 
$x_1,\ldots,x_n$.  Second, we claim that for any $\mathcal{F}$ there exists a 
finite set of disjoint sets $S_1,\ldots,S_k\in\mathcal{F}$ such that 
$\bigcup_kS_k=\Omega$ and such that any $x\in\ell^\infty(\mathcal{F})$ can be 
written as $x=\sum_ix_i\chi_{S_i}$ ($\chi_S$ is the indicator function on 
$S$, $x_i\in\mathbb{C}$).  To see this, note that there is only a finite 
possible number of disjoint partitions $\{T_i\}$ of $\Omega$, the latter 
being a finite set.  $\{S_i\}$ is then the (unique) finest such partition 
that is a subset of $\mathcal{F}$.  To show uniqueness: if
$\{S_i\},\{S_i'\}\subset\mathcal{F}$ were two such partitions then $\{S_i\cap
S_j'\}\subset\mathcal{F}$ would be a finer partition unless $\{S_i\}=\{S_i'\}$.

Let us now prove the Lemma.  For a diagonal matrix $X$ we define the map
$x(i)=\iota(X)(i)=X_{ii}$, and similarly for a map $x:\Omega\to\mathbb{C}$
define the diagonal matrix $X_{ii}=\iota^{-1}(x)_{ii}=x(i)$.  This gives a
$^*$-isomorphism $\iota$ between the set $\mathscr{A}$ of all diagonal
matrices and the set $\ell^\infty(\Omega)$ of all maps on $\Omega$. The
claim of the Lemma is that $\iota$ maps any $^*$-subalgebra
$\mathscr{C}\subset\mathscr{A}$ to $\ell^{\infty}(\mathcal{F})$ for some
$\sigma$-algebra $\mathcal{F}$, and conversely that
$\iota^{-1}(\ell^{\infty}(\mathcal{F}))$ is a $^*$-algebra for any $\mathcal{F}$.  
The latter is a straightforward consequence of the definitions, so let us
prove the former statement. Fix the $^*$-algebra $\mathscr{C}$ and define
$\mathcal{C}=\iota(\mathscr{C})$.  Let $\mathcal{F}=\sigma\{\mathcal{C}\}$, and
introduce the smallest $\Omega$-partition $\{S_i\}\subset\mathcal{F}$ as above.  
If $\chi_{S_i}\in\mathcal{C}$ for every $i$, then we must have
$\mathcal{C}=\ell^\infty(\mathcal{F})$: $\mathcal{C}$ is an algebra and hence
contains all random variables of the form $\sum_ix_i\chi_{S_i}$.  Let us
thus assume that there is some $i$ such that $\chi_{S_i}\not\in
\mathcal{C}$; as $\mathcal{F}=\sigma\{\mathcal{C}\}$, however, we must be able
to find functions $x_{(1)},\ldots,x_{(n)}\in\mathcal{C}$ such that
$x_{(1)}^{-1}(x_1)\cap\cdots\cap x_{(n)}^{-1}(x_n)=S_i$.  But if we choose
$f(y_1,\ldots,y_n)=\chi_{\{y_1=x_1,\ldots,y_n=x_n\}}(y_1,\ldots,y_n)$ then
$f(x_{(1)},\ldots,x_{(n)})=\chi_{S_i}\in\mathcal{C}$.  Hence we have a
contradiction, and the Lemma is proved.  
\qquad \end{proof}

We are now ready to introduce the basic mathematical object studied in 
this article: a {\it generalized} or {\it quantum probability space}.
Applying the spectral theorem to this structure gives a fundamental 
result which we will use over and over.

\begin{definition}[\bf Quantum probability space]
The pair $(\mathscr{A},\rho)$, where $\mathscr{A}$ is a $^*$-algebra of 
operators on a finite-dimensional Hilbert space $\mathsf{H}$ and $\rho$ 
is a state on $\mathscr{A}$, is called a (finite-dimensional) {\em quantum 
probability space}.
\end{definition}

\begin{theorem}[\bf Spectral theorem for quantum probability spaces]
\label{thm:specth}
Let $(\mathscr{C},\rho)$ be a {\em commutative} quantum probability space.
Then there is a probability space $(\Omega,\mathcal{F},\mathbf{P})$ and a 
$^*$-isomorphism $\iota:\mathscr{C}\to\ell^\infty(\mathcal{F})$
	\nomenclature{$\iota$}{$^*$-isomorphism between observables
	and random variables}
such that $\rho(X)=\mathbf{E_P}(\iota(X))$ $\forall\,X\in\mathscr{C}$.
\end{theorem}

What is so general about a \emph{generalized} or \emph{quantum}
probability space?  It is the existence of many commutative subalgebras
within $(\mathscr{A},\rho)$. Theorem \ref{thm:specth} does not apply
directly to $(\mathscr{A},\rho)$, as usually such a space will not be
commutative.  On the other hand, in a single realization we can choose to
observe at most a commutative set of observables, which generate a
commutative subalgebra $\mathscr{C}\subset\mathscr{A}$.  The probability
space $(\mathscr{C},\rho|_{\mathscr{C}})$ is commutative and is thus
exactly equivalent to a classical probability space by Theorem
\ref{thm:specth}.  The noncommutative probability space
$(\mathscr{A},\rho)$ describes the statistics of all possible
experiments---it is a collection of many incompatible classical
probability models, each of which coincides with a commutative subalgebra
of $\mathscr{A}$.  The experiment we choose to perform in one single
realization determines which commutative subalgebra of $\mathscr{A}$,
i.e.\ which classical probability space, is needed to describe the random
outcomes of that experiment.

\begin{remark}
In section \ref{sec:rvini} we determined the probability measure 
$\mathbf{P}$ using a state vector $\psi\in\mathsf{H}$; here we have 
replaced this notion by the {\em state} $\rho:\mathscr{C}\to\mathbb{C}$.   
This is in fact a generalization: the vector $\psi$ corresponds to the 
state $\rho(X)=\langle\psi,X\psi\rangle$.  In general we can always 
characterize a state $\rho$ by a {\it density matrix} $\tilde\rho$ as 
follows: $\rho(X)={\rm Tr}[\tilde\rho X]$.  This follows directly from 
linearity of $\rho$.  Positivity and normalization impose the additional 
conditions $\rho\ge 0$ and ${\rm Tr}\,\rho=1$.  A state of the form 
$\rho(X)=\langle\psi,X\psi\rangle={\rm Tr}[(\psi\psi^*)X]$ is known as a 
{\it pure} or {\it vector state}, whereas any other state is known as a 
{\it mixed state}.  Both state vectors and density matrices are commonly 
used in the physics literature, while $\rho:\mathscr{C}\to\mathbb{C}$ is 
usual in quantum probability.
\end{remark}

We will often speak of the $^*$-algebra {\em generated} by a set of 
observables $X_1,\ldots,X_k$.  By this we mean the smallest $^*$-algebra 
${\rm alg}\{X_1,\ldots,X_k\}$ of operators on $\mathsf{H}$ that contains 
$X_1,\ldots,X_k$.  This notion plays exactly the same role as the 
$\sigma$-algebra generated by a set of random variables---indeed, it is 
straightforward to verify that if $X_1,\ldots,X_k$ commute, then
$\iota({\rm alg}\{X_1,\ldots,X_k\})=\ell^\infty(\sigma\{\iota(X_1),\ldots,\iota(X_k)\})$.

In the following we will need to construct quantum probability spaces from
a number of independent smaller probability spaces.  The analogous
classical notion can be illustrated as follows: suppose
$(\Omega_i,\mathcal{F}_i,\mathbf{P}_i)$, $i=1,2$ are two independent copies of
the probability space of a coin flip.  Then
$(\Omega_1\times\Omega_2,\mathcal{F}_1\times\mathcal{F}_2,\mathbf{P}_1\times\mathbf{P}_2)$
is the joint probability space for the two independent coin flips, on
which any random variable $f(\omega_1,\omega_2)=f(\omega_1)$ which depends
only on the first coin and $g(\omega_1,\omega_2)=g(\omega_2)$ which
depends only on the second coin are independent (under the product measure
$\mathbf{P}_1\times\mathbf{P}_2$).

The analogous construction for quantum probability spaces uses the tensor
product; that is, given quantum probability spaces
$(\mathscr{A}_i,\rho_i)$ defined on the Hilbert spaces $\mathsf{H}_i$, we
construct the joint space
$(\mathscr{A}_1\otimes\mathscr{A}_2,\rho_1\otimes\rho_2)$ on the Hilbert
space $\mathsf{H}_1\otimes\mathsf{H}_2$.  Here
$\mathsf{H}_1\otimes\mathsf{H}_2$ and
$\mathscr{A}_1\otimes\mathscr{A}_2$ denote the usual vector and matrix
tensor products, while $\rho_1\otimes\rho_2$ is defined by
$(\rho_1\otimes\rho_2)(A\otimes B)=\rho_1(A)\rho_2(B)$ and is extended to
all of $\mathscr{A}_1\otimes\mathscr{A}_2$ by linearity.  It is not
difficult to verify that also this construction is consistent with the
classical case, as $\ell^\infty(\mathcal{F}_1\times\mathcal{F}_2)=
\ell^\infty(\mathcal{F}_1)\otimes\ell^\infty(\mathcal{F}_2)$; the $^*$-isomorphism
$\iota$ obtained through the spectral theorem then maps observables of the
form $A\otimes I$ to random variables of the form
$f(\omega_1,\omega_2)=f(\omega_1)$, and vice versa.

For the time being this is all the general theory that we need.  We now 
turn to the construction of the class of models that we will consider 
throughout the paper.

\subsection{Two-level systems}

In this article all quantum systems will be built up (by means of tensor
products) from two-level systems, i.e.\ quantum probability spaces
with ${\rm dim}\,\mathsf{H}=2$.  Any observable of a two-level system is
at most a two-state random variable, which makes the theory particularly
simple.  Nonetheless we will find a surprisingly rich structure in the
resulting models.  Moreover, it turns out that such models can
approximate closely the behavior of realistic physical systems, see
section \ref{continuous time}.  The classical counterpart of such models,
i.e.\ the approximation of a continuous system by a sequence of coin
flips, is well known e.g.\ in mathematical finance \cite{CRR79,Shr04}.

We already encountered a two-level system in section \ref{sec:rvini}.  As
we will use these systems so often, however, it will be useful to fix some
notation and to introduce some standard vectors and operators which we 
will use throughout.

Let $\mathscr{M}$
	\nomenclature{$\mathscr{M}$}{Two-level system algebra}
denote the $^*$-algebra of all $2\times 2$ complex 
matrices acting on the Hilbert space $\mathsf{H}=\mathbb{C}^2$. This 
algebra, together with suitable choices of states, will be the building 
block in our model of an atom in interaction with an electromagnetic 
field: the atom will be built on $\mathscr{M}$ (see Example 
\ref{ex:tla}), while the field will be modeled by a tensor product of 
$\mathscr{M}$'s  (see section \ref{sec:repeatedinteraction}).

The canonical two-level system is described by the quantum probability 
space $(\mathscr{M},\rho)$, where $\rho$ is some state on $\mathscr{M}$.  
Implicit in this description is the choice of a standard basis (hence the 
name canonical) which allows us to represent linear operators by the more 
convenient matrices; there is no loss of generality in doing this, and it 
will make our life much easier.  We introduce the standard definitions:
	\nomenclature{$\sigma_\pm,\sigma_z$}{Pauli matrices in
		$\mathscr{M}$}
	\nomenclature{$\Phi$}{Vacuum vector in $\mathbb{C}^2$}
\begin{equation}\label{eq basis}
	\sigma_- = \begin{pmatrix} 0 & 0 \\ 1 & 0 \end{pmatrix}, \qquad
	\sigma_+ = \sigma_-^* = \begin{pmatrix} 0 & 1 \\ 0 & 0 \end{pmatrix}, \qquad
	\Phi = \begin{pmatrix}0 \\1\end{pmatrix}.
\end{equation}
$\Phi$ is called the \emph{vacuum vector} or \emph{ground state vector}
for reasons which will become clear shortly.  Note that $\sigma_-\Phi=0$,
$\sigma_-^2=0$, and $\sigma_-\sigma_+ + \sigma_+\sigma_- = I$.
$\sigma_+$ and $\sigma_-$ are convenient operators, as any matrix 
$M\in\mathscr{M}$ can be written as a linear combination of 
$\sigma_+\sigma_-,\sigma_+,\sigma_-,I$.  We will also sometimes use the 
\emph{Pauli matrix}
\begin{equation}\label{eq Pauli}   
	\sigma_z = \sigma_+\sigma_- - \sigma_-\sigma_+ = 
	\begin{pmatrix}1 & 0 \\ 0 & -1 \end{pmatrix}.
\end{equation}

\begin{example}(\textbf{Two-level atom}).  \label{ex:tla}
In this example we are going to interpret the quantum probability space
$(\mathscr{M}, \rho)$ as the set of observables describing a single atom.  
Introduce the self-adjoint operator $H\in\mathscr{M}$ by
\begin{equation}\label{eq Hamiltonian}
	H = \frac{\hbar\omega_0}{2}\sigma_z = 
	\begin{pmatrix}\frac{\hbar\omega_0}{2} &0 
	\\ 0 &-\frac{\hbar\omega_0}{2}\end{pmatrix}, 
\end{equation}
where $\hbar$ is Planck's constant (a fundamental constant of nature), 
$\omega_0$ is a parameter representing a frequency (the so-called 
\emph{atomic frequency}), and $\sigma_z$ is the Pauli matrix given by Eq.\ 
\eqref{eq Pauli}.  As an observable, the \emph{Hamiltonian} $H$ is 
interpreted as the \emph{energy} of the atom (indeed, $\hbar\omega_0$ 
has units of energy).  If an experiment is done to measure it, then from 
Eq.\ \eqref{eq Hamiltonian} it is clear that the outcome is a random 
variable taking either the value $\hbar\omega_0/2$ or $-\hbar\omega_0/2$ 
(note that $H$ is already diagonal in the given basis).

The fact that the atomic energy takes discrete values is a fundamental
result in quantum mechanics to which the theory owes its name.  No actual
atom has only two energy levels, but in some cases atomic physics is well
approximated by such a description.  Other physical observables, such as
the intrinsic angular momentum (spin) of an electron, are exactly
described by a two-level system.  Throughout this article we will use the
two-level system as a prototypical example of an atom.

We have not yet discussed the state $\rho$; the state depends on the way
the atom was prepared (e.g.\ is the atom at room temperature, has it been
cooled to absolute zero, or perhaps it was excited by a laser?)  A
particularly interesting state, the \emph{ground state}, is the one in
which the energy takes its lowest possible value $-\hbar\omega_0/2$ with
unit probability.  This situation describes an atom that has been cooled
to absolute zero temperature (a good approximation for laser-cooled
atoms).  It is not difficult to verify that the ground state is defined by
$\rho(X)=\langle \Phi, X\Phi\rangle$ where $\Phi$ is the ground state
vector.  In general we can choose any state for the atom, i.e.\ we set
$$
	\rho(X)={\rm Tr}[\tilde\rho X]\qquad
	\mbox{with}\qquad
	\tilde\rho=\begin{pmatrix}
		\rho_{11} & \rho_{12} \\
		\rho_{21} & \rho_{22}
	\end{pmatrix},
	\quad X\in\mathscr{M},
$$
where $\tilde\rho$ is a density matrix ($\tilde\rho\ge 0$, ${\rm 
Tr}\,\tilde\rho=1$).  Then measurement of $H$ yields the outcome 
$-\hbar\omega_0/2$ with probability $\rho_{22}$ and $\hbar\omega_0/2$ with 
probability $\rho_{11}$.  Note that the expectation value of the energy 
$H$ is given by $\rho(H)={\rm Tr}[\tilde\rho 
H]=\frac{1}{2}\hbar\omega_0(\rho_{11}-\rho_{22})$. 

Apart from the energy $H$ the two-level atom may possess an 
\emph{electric dipole moment}, described by the vector observable 
$\boldsymbol{\mu}$.  Its three components are given by
\begin{equation*}
	\mu_a = \alpha_a \sigma_- + \alpha_a^*\sigma_+ = 
	\begin{pmatrix} 0 & \alpha_a^* \\  \alpha_a & 0\end{pmatrix}, 
	\qquad
	a\in \{x,y,z\}, 
\end{equation*}
where $\alpha_x$, $\alpha_y$ and $\alpha_z$ are complex parameters 
(depending on which specific atom we are considering).  Note that 
$\mu_a$, $a\in \{x,y,z\}$ does not commute with $H$. Therefore $\mu_a$ and 
$H$ cannot be simultaneously diagonalized by the spectral theorem, that is 
they cannot be represented on the same probability space. It is impossible 
to devise an experiment that measures both observables simultaneously. 
Note further that $\mu_{x,y,z}$ do not commute with each other unless
$\alpha_{x,y,z}$ have the same phase, i.e.\ it is in general not even 
possible to measure all three components of ${\bf \mu}$ in one realization.
\end{example}

\subsection{The discrete field}\label{sec:discretefield}

Having obtained a simple model for an atom, we proceed by introducing a
simple model for the electromagnetic field, e.g.\ light emitted from a
laser.  In the next section we will couple the atom and the field, which
will allow us to describe measurements of light emitted from or scattered
by the atom.  Before we can do this, however, we must model the field in 
absence of the atom.

\begin{figure}
\centering
\includegraphics[width=\textwidth]{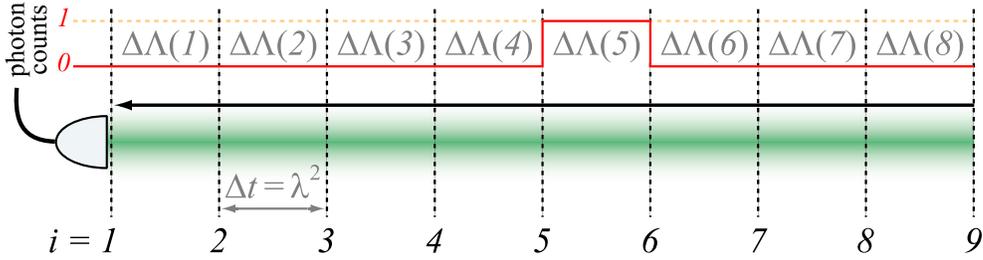}
\caption{\label{fig:noise} Cartoon of a discretized model for the 
electromagnetic field.  The field (e.g.\ a laser beam) propagates to the 
left until it hits a photodetector.  When we discretize time with 
sufficiently small $\Delta t=\lambda^2$ then to good approximation the 
photodetector will detect either zero or one photons in each time step, 
and moreover the number of photons detected in each time interval are 
independent random variables.  The observable $\Delta\Lambda(i)$ 
represents the number of photons detected in the time interval $i$.  A 
possible sample path for the stochastic process $\iota(\Delta\Lambda(i))$ 
is displayed in red.} 
\end{figure}

The model we have in mind is illustrated in Figure \ref{fig:noise}.
Imagine a laser beam which we are measuring with a photodetector.  The
beam propagates towards the photodetector at the speed of light, so at
each time $t$ we will measure a different part of the beam.  We wish to
describe the measurements made by the photodetector at each time $t$.  To
simplify matters we will discretize time into slices.  If these slices are
sufficiently small (as compared to the intensity of the laser), then to
good approximation the photodetector will measure no more than one photon
in each slice.  Moreover, the measurement in each slice will be
independent from the other slices, as we are measuring independent parts
of the field at each time step.  Hence we will model the photodetection of
the field by a collection of independent $\{0,1\}$-valued random
variables, one for each time slice, corresponding to whether a photon was
$(1)$ or was not $(0)$ detected in each slice.

The model being a quantum mechanical model, the above description is not
sufficient.  The events detected by a photodetector are classical events,
corresponding to an entire sample path of observations measured in a
single realization, and hence are necessarily described by a commutative
algebra.  Other observations of the field are possible, using a different
detection setup, which are incompatible with direct photodetection (we
will encounter an example below).  Hence we need a noncommutative model.  
Still, the considerations above are very suggestive: we will build our 
model from a collection of independent two-level systems, one for each 
time slice.

Fix an integer $k$ and finite time interval $[0,T]$, $T>0$. Define 
$\lambda=\sqrt{T/k}$,
	\nomenclature{$k$}{Number of time slices}
	\nomenclature{$\lambda^2$}{Length of each time slice}
i.e.\ the interval $[0,T]$ is divided into $k$ parts 
of equal length $\lambda^2$.  To every time slice $1,\ldots,k$ we want to 
associate an independent copy of the canonical two-level system 
$(\mathscr{M},\rho)$.  But we have already seen how to do this: we must 
use tensor products, i.e.\ we introduce the quantum probability space
	\nomenclature{$\mathscr{W}_k$}{Field algebra}
\begin{equation*}
	(\mathscr{W}_k,\rho_k) = (\mathscr{M}^{\ten k},\rho^{\ten k})=
	(\underbrace{\mathscr{M}\otimes\cdots\otimes\mathscr{M}}_{
	k\mathrm{~times}},~
	\underbrace{\rho\otimes\cdots\otimes\rho}_{k\mathrm{~times}}). 
\end{equation*}
For any $X\in\mathscr{M}$ and integer $1 \le i \le k$, we define an
element $X_i\in\mathscr{W}_k$ by
\begin{equation*}
	X_i = I^{\otimes (i-1)}\otimes X\otimes I^{\otimes (k-i)}.
\end{equation*}
$X_i$ is the observable $X$ in the time slice $i$.  Note the fundamental 
property $[X_i,Y_j]=0$ $\forall\,X,Y\in\mathscr{M}$, $i\ne j$: 
measurements on different slices of the field are always compatible.
For future reference, we also define the algebra $\mathscr{M}_i$ of
observables in time slice $i$:
	\nomenclature{$\mathscr{M}_l$}{Field algebra of time slice $l$}
$$
	\mathscr{M}_i={\rm alg}\{X_i:X\in\mathscr{M}\}.
$$

\subsubsection*{Discrete noises}
For $1 \le i \le k$ we define the standard increments 
	\nomenclature{$A(l),A^*(l),\Lambda(l),t(l)$}{Discrete noises}
\begin{equation}\label{eq def noises}\begin{split}
	  &\Delta A(i) = \lambda (\sigma_-)_i, \ \ \ \ \
	  \Delta\Lambda(i) = (\sigma_+\sigma_-)_i, \\
	  &\Delta A^*(i) = \lambda(\sigma_+)_i,  \ \ \ \ \
	  \Delta t(i) = \lambda^2 I,
\end{split}\end{equation}
and consequently we introduce the discrete noise processes $A$, $A^*$,
$\Lambda$ and the discrete time process $t$ (with the convention 
$A(0)=A^*(0)=\Lambda(0)=t(0)=0$) as follows:
\begin{equation*}\begin{split}
	  &A(l) = \sum_{i=1}^l \Delta A(i), \ \ \ \ \
	  \Lambda(l) = \sum_{i=1}^l \Delta\Lambda(i) \\
	  &A^*(l) = \sum_{i=1}^l \Delta A^*(i), \ \ \ \ \
	  t(l) = \sum_{i=1}^l \Delta t(i).
\end{split}\end{equation*}
The particular scaling of these expressions with $\lambda$ has been 
selected so that the definitions make sense in the continuous limit 
$\lambda\to 0$; see section \ref{continuous time}.  For now, let us 
furnish these objects with a suitable interpretation.  $t(l)=l\lambda^2I$ 
is the easiest: this observable takes the value $l\lambda^2$ with 
probability one.  We interpret $t(l)$ as the time that has elapsed after 
time slice $l$.  Next, let us investigate $\Lambda(l)$.  

\subsubsection*{Photon counting}
For the observable $\sigma_+\sigma_-\in\mathscr{M}$, it is 
easily verified that the spectral theorem maps to a random variable 
$\iota(\sigma_+\sigma_-)$ that takes the value zero or one; indeed, the 
matrix corresponding to $\sigma_+\sigma_-$ (calculated using Eq.\ 
(\ref{eq basis})) is given by
$$
	\sigma_+\sigma_- = 
	\begin{pmatrix}1 & 0 \\ 0 & 0 \end{pmatrix},
$$
which is already in diagonal form.  $\Delta\Lambda(i)$, $i=1,\ldots,k$ 
is a set of independent observables taking the value zero or one.  We 
interpret $\Delta\Lambda(i)$ as the number of photons observed by a 
photodetector in the time slice $i$ (as illustrated in Fig.\ 
\ref{fig:noise}).

The probability of observing a photon in the time slice $i$ depends on the
state $\rho$; in the presence of a light source (a laser or a light bulb)  
this probability will be nonzero.  If we turn off the light source, or if
there was no light source in the first place, then the probability of
measuring a photon will be zero in every time slice.  As we discussed in
the context of a two-level atom, this situation is described by the state
$\phi(X)=\langle\Phi,X\Phi\rangle$ 
	\nomenclature{$\phi$}{Vacuum state on $\mathscr{M}$}
(i.e.\ $\rho_k=\phi^{\ten k}$)---it is 
for this reason that this state is called the \emph{vacuum state}.  Even 
in the absence of an external light source an atom can still interact with
the vacuum: e.g.\ an excited atom can spontaneously emit a photon into the
vacuum, an event which we will want to measure using our photodetector.
For concreteness, \emph{we will always work in the following with the 
vacuum state $\phi^{\ten k}$}.

Now recall that $\Delta\Lambda(i)$, $i=1,\ldots,k$ all commute with each
other, and hence generate a commutative $^*$-algebra $\mathscr{C}_k$. The
spectral theorem maps $\mathscr{C}_k$ to a classical space $(\Omega_k,
\mathcal{F}_k)$, which is simply the probability space corresponding to $k$
independent coin flips (a ``binomial model'' in the terminology of
\cite{Shr04}).  The classical stochastic process
$y_l=\iota(\Delta\Lambda(l))$ is then precisely the signal we observe
coming from the photodetector.  By applying the spectral theorem to the
commutative subalgebras $\mathscr{C}_l=\mathrm{alg}\{\Delta\Lambda(i):
i=1,\ldots,l\}$, $\mathscr{C}_l\subset\mathscr{C}_{l+1}\subset\cdots
\subset\mathscr{C}_k$, we obtain an increasing family of $\sigma$-algebras
$\mathcal{F}_l\subset\mathcal{F}_{l+1}\subset\cdots\subset\mathcal{F}_k$, the {\it
filtration} generated by the stochastic process $y_l$.  As in the
classical theory, the quantum filtration $\{\mathscr{C}_l\}$ allows us to
keep track of what information is available based on observations up to
and including time slice $l$.

It remains to consider $\Lambda(l)\in\mathscr{C}_l$: it is evident from
the definition that this observable represents the number of photons
observed up to and including time slice $l$.  Hence the \emph{number
process} $\{\Lambda(l)\}_{l=1,\ldots,k}$, which is clearly commutative as
$\Lambda(l)\in\mathscr{C}_k$ $\forall l$, maps to a classical counting
process $\iota(\Lambda(l))$.  All of this is not so interesting when we
are working under the measure $\mathbf{P}_k$ induced by the vacuum
$\phi^{\ten k}$ on $(\Omega_k,\mathcal{F}_k)$: under this measure any sample
path of $\iota(\Lambda(l))$ that is not identically zero has zero
probability.  Once we have introduced the interaction with an atom,
however, the number process will allow us to count the number of photons
emitted by the atom.

\subsubsection*{Homodyne detection}
The processes $A(l)$ and $A^*(l)$ are not self-adjoint, but $X_\varphi(l)=
e^{i\varphi}A(l)+e^{-i\varphi}A^*(l)$ is an observable for any phase 
$\varphi\in[0,\pi)$.  Moreover, $\{X_\varphi(l)\}_{l=1,\ldots,k}$ is a 
commutative process for fixed $\varphi$: to see this, note that $\Delta 
X_\varphi(l)=X_\varphi(l)-X_\varphi(l-1)=\lambda(e^{i\varphi}\sigma_- +
e^{-i\varphi}\sigma_+)_l$ commute with each other for different $l$, so 
the same must hold for $X_\varphi(l)$.  Hence, as above, we can use the 
spectral theorem to map $\mathscr{C}_k^\varphi=\mathrm{alg}\{X_\varphi(l):
l=1,\ldots,k\}$ to a classical space $(\Omega_k^\varphi,\mathcal{F}_k^\varphi)$, 
and similarly we obtain the filtration $\{\mathcal{F}_l^\varphi\}_{l=1,\ldots,k}$.
Beware, however: $\{X_\varphi\}$ does not commute with $\{\Lambda\}$, nor 
does it commute with $\{X_{\varphi'}\}$ for $\varphi'\ne\varphi$.  To 
observe each of these processes we need fundamentally different detectors, 
so that in any given realization at most one of these processes is 
available to us.

An optical detector that measures one of the processes $X_\varphi(l)$ is 
known as a \emph{homodyne detection setup}, a term borrowed from radio 
engineering.  What is measured in such a setup is usually interpreted by 
considering the electromagnetic field to be a wave at a certain frequency, 
as in classical electromagnetism.  The homodyne detector then measures 
the component of this wave which is in phase with a certain reference 
wave, known as the local oscillator, and the parameter $\varphi$ 
determines the phase of the local oscillator with respect to an external 
reference.  The interpretation of the light in terms of waves (homodyne) 
or particles (photodetection) are not at odds, because the corresponding 
measurements $\{X_\varphi\}$ and $\{\Lambda\}$ do not commute.  We will 
not dwell further on the physical interpretation of these detection 
schemes; suffice it to say that both homodyne detection and photodetection 
are extremely common techniques, either of which may be more convenient in 
any particular situation.  The interested reader is referred to the 
quantum optics literature \cite{WaM94} for further information (see 
also \cite{Bar03} for a more mathematical perspective).

For concreteness we will always consider the case $\varphi=0$ whenever we
discuss homodyne detection in the following, i.e.\ we will consider the
observation process $X(l)=A(l)+A^*(l)$.  Let us investigate this process a
little more closely.  First, consider the two-level system observable 
$\sigma_++\sigma_-$:
$$
	\sigma_++\sigma_-=
	\begin{pmatrix}
		0 & 1 \\ 1 & 0
	\end{pmatrix} =
	\begin{pmatrix}
		\frac{1}{2} & \frac{1}{2} \\ \frac{1}{2} & \frac{1}{2}
	\end{pmatrix}-
	\begin{pmatrix}
		\frac{1}{2} & -\frac{1}{2} \\ -\frac{1}{2} & \frac{1}{2}
	\end{pmatrix}=P_1-P_{-1}.
$$
$P_{\pm 1}$ are projections, so we can read off that $\sigma_++\sigma_-$ 
has eigenvalues $\pm 1$.  Moreover, note that $\langle\Phi,P_{\pm 
1}\Phi\rangle = 1/2$.  Hence evidently the spectral theorem maps 
$\sigma_++\sigma_-$ to a random variable taking the values $\pm 1$ each 
with probability $1/2$ in the vacuum.  It follows immediately that 
$\iota(\Delta X(l))$, $l=1,\ldots,k$ are independent random variables 
taking values $\pm\lambda$ with equal probability.  Unlike in the case of 
photodetection, the homodyne detector gives a noisy signal even in the 
vacuum: a manifestation of {\it vacuum fluctuations}, as this is called in 
quantum optics.  The process $\iota(X(l))=\sum_{i=1}^l\iota(\Delta X(i))$, 
the {\it integrated photocurrent}, is evidently a symmetric random walk on 
the grid $\lambda\mathbb{Z}$.  It is not difficult to see that 
$\iota(X(\lfloor t/\lambda^2\rfloor))$ converges in law to a Wiener 
process as $\lambda\to 0$ (recall $t(l)=l\lambda^2$, so $l(t)\sim 
\lfloor t/\lambda^2\rfloor$)---see section \ref{continuous time}.

\begin{remark} The reader might wonder at this point why we have chosen
the notation $\iota(\Delta\Lambda(l))$ or $\iota(\Delta X(l))$ for the
observations---why not dispose of the $\Delta$'s?  The reason is that
neither $\Delta X(\lfloor t/\lambda^2\rfloor)$ nor $\Delta\Lambda(\lfloor
t/\lambda^2\rfloor)$ have well-defined limits as $\lambda\to 0$; for
example, $\iota(\Delta X(\lfloor t/\lambda^2\rfloor))$ would have to
converge to white noise, which is not a mathematically well-defined object
(at least in the sense of stochastic processes).  The Wiener process, on
the other hand, has a rigorous definition.  We use the convention that
processes such as $X(\lfloor t/\lambda^2\rfloor)$ have well-defined
limits, whereas $\Delta X(\lfloor t/\lambda^2\rfloor)$ will ``converge''
to $dX_t$ which is only meaningful under the integral sign. \end{remark}

\begin{remark}
It should be emphasized that the model we have introduced here is rather
crude; certainly it was not derived from physical principles!  Its main 
physical justification is that it converges in the limit $\lambda\to 0$
to a model which {\it is} obtained from physical principles, so that for 
small $\lambda$ we can consider it to be a good approximation---see the 
references in section \ref{continuous time}.  The same disclaimer holds 
for what we will discuss in the remainder of section \ref{sec random 
variables}.
\end{remark}

\subsection{Repeated interaction}\label{sec:repeatedinteraction}

Now that we have descriptions of an atom and of the electromagnetic field,
it is time to describe how these objects interact with each other.  The
atom can emit photons into the vacuum, which we subsequently detect using
a photodetector; alternatively, we can measure using a homodyne detection
setup how the presence of the atom perturbs the usual vacuum fluctuations.
The remainder of the article is centered around the theme: how can we put
such measurements to good use?  The development of \emph{filtering theory}
allows us to estimate atomic observables based on the measurements in the
field, and we can subsequently use this information to control the atom
by means of \emph{feedback control}.

\begin{figure}
\centering
\includegraphics[width=\textwidth]{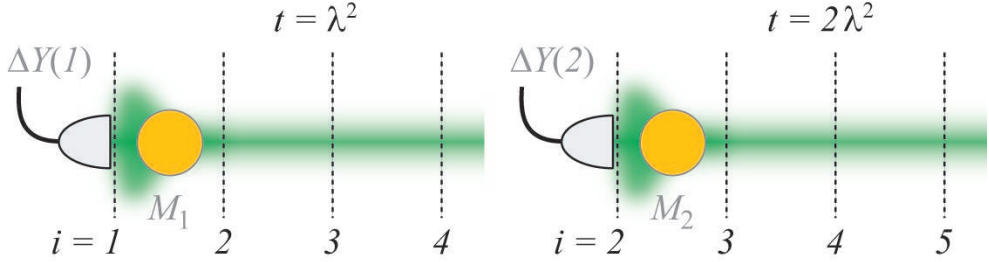}
\caption{\label{fig:qsde} Cartoon of the discretized atom-field 
interaction.  In each time step $i$, the atom and the field at time slice 
$i$ interact through a unitary transformation $M_i$.  The time slice of 
the field is subsequently measured (using a photodetector or homodyne 
detection), which gives rise to the observation $\Delta Y(i)$.  The field 
then propagates to the left by one time slice, and the procedure repeats.
} \end{figure}

The way the interaction works is illustrated in Figure \ref{fig:qsde}. As
before, the field propagates towards the detector, which measures
sequentially time slices of the field.  Now, however, we place the atom
right in front of the detector, so that the current time slice interacts
with the atom prior to being detected.  This sort of model is called a
repeated interaction: the atom is fixed, and interacts identically with
each time slice of the field exactly once before it is detected.
Before we can describe the repeated interaction in detail, we first have
to describe how a single interaction is modeled in quantum mechanics.
Though this article is mostly about quantum \emph{probability}, we need a 
little bit of quantum \emph{dynamics} in order to model physical systems.

\subsubsection*{Interactions in quantum mechanics}

Suppose we have two independent two-level atoms, i.e.\ we are working on 
the quantum probability space $(\mathscr{M}\otimes\mathscr{M},
\rho\otimes\rho)$.  Then, for example, the observable $\mu_x^{1}=
\mu_x\otimes I$ corresponds to the electric dipole moment of the first 
atom, and is independent from any observable $I\otimes X$ of the second 
atom.  If the two atoms subsequently interact, however, then the dipole 
moment of the first atom would likely change; i.e.\ the dipole moment 
after the interaction is described by a different observable 
$\mu_x^{1,{\rm after}}$ than before the interaction, and will likely be 
correlated with some observable of the second atom $(I\otimes X)^{\rm after}$.

In quantum mechanics, every interaction is described by a {\it unitary
operator}.  In the previous example, the interaction between the two atoms
corresponds to a suitably chosen unitary map
$U:\mathbb{C}^2\otimes\mathbb{C}^2\to\mathbb{C}^2\otimes\mathbb{C}^2$.  
This means that any observable $X\in\mathscr{M}\otimes\mathscr{M}$ 
transforms to $U^*XU\in\mathscr{M}\otimes\mathscr{M}$ after the two atoms 
have interacted.  In particular, $\mu_x^{1,{\rm after}}=U^*\mu_x^1U$ will 
generally be different than $\mu_x^1$, and is likely to be correlated with 
$U^*(I\otimes X)U$ for some $X$.  Note that unitary rotation does not 
change the spectrum of an observable; i.e., $\iota(\mu_x^1)$ and 
$\iota(\mu_x^{1,{\rm after}})$ take the same values, but they are defined 
on different probability spaces with different probability measures.

Recall that any unitary matrix $U$ can be written as $U=e^{-iB}$ for 
some self-adjoint matrix $B$.  Hence we can express any unitary 
transformation of $\mathscr{M}\otimes\mathscr{M}$ by
$$
	U=\exp\left(-i\left\{L_1\otimes\sigma_+\sigma_-+L_2\otimes\sigma_++
		L_2^*\otimes\sigma_-+L_3\otimes I\right\}
	\right),
$$
where $L_{1,2,3}\in\mathscr{M}$ and $L_1,L_3$ are self-adjoint.
	\nomenclature{$L_{1,2,3}$}{Interaction operators}

\begin{remark} 
We have described interactions by changing the observables corresponding
to a particular physical quantity.  This is how dynamics is usually
represented in probability theory.  For example, the classical
discrete-time Markov chain $X_{n}=f(X_{n-1},\xi_n)$ where $\xi_n$ is
i.i.d.\ noise is precisely of this form: the random variables $X_n$ and
$X_{n-1}$ represent the same physical quantity $X$, but the interaction
with the noise $\xi$ means we have to represent $X$ by different random
variables $X_n$ at different times.  In quantum mechanics, this 
is called the \emph{Heisenberg picture}.  In the \emph{Schr{\"o}dinger 
picture}, any physical quantity is always represented by the same 
observable, but the underlying state is changed by the interaction.
The two pictures are in some sense dual to each other, as is discussed in 
any quantum mechanics textbook (e.g.\ \cite{EM98}).  However, the 
Schr{\"o}dinger picture is somewhat unnatural if one wants to work with 
(quantum) stochastic processes.  In this article we will always work in 
the Heisenberg picture.
\end{remark}

\subsubsection*{Atom-field interaction}

We now consider the case that a two-level atom is coupled to the vacuum
electromagnetic field, that is, we work on the quantum probability space
$(\mathscr{M}\otimes\mathscr{W}_k,\mathbb{P})$ with
$\mathbb{P}=\rho\otimes\phi^{\ten k}$.
	\nomenclature{$\mathbb{P}$}{Standard state on the atom-field 
	algebra $\mathscr{M}\otimes\mathscr{W}_k$}
The subspace
$(\mathscr{M},\rho)\subset (\mathscr{M}\otimes\mathscr{W}_k,\mathbb{P})$
is known as the \emph{initial system} ($\rho$ is called the 
\emph{initial state}), because observables of the form $X\otimes I$ in 
$\mathscr{M}\otimes\mathscr{W}_k$ represent the atom's physical quantities 
at the initial time $l=0$, i.e., before the atom has interacted with the 
field.  The interaction requires us to modify these observables at each 
time $l>0$, as we will now describe.

Let us begin by considering what happens to the system after one time 
step.  At this time ($l=1$) the atom has interacted only with the 
first time slice of the field, i.e.\ any atomic observable $X\otimes I$ at 
time $l=0$ will evolve to $j_1(X)=M(1)^*(X\otimes I)M(1)$, where $M(1)$ is 
the unitary interaction of the form
$$
	M(1)=
	e^{-i\{L_1\Delta\Lambda(1)
	+L_2\Delta A^*(1)
	+L_2^*\Delta A(1)
	+L_3\otimes\Delta t(1)\}
	}.
$$
In the next time step, we want the atom to interact in an identical 
fashion with the second time slice, so we will define 
$j_2(X)=M(2)^*j_1(X)M(2)$ with
$$
	M(2)=
	e^{-i\{j_1(L_1)\Delta\Lambda(2)
	+j_1(L_2)\Delta A^*(2)
	+j_1(L_2^*)\Delta A(2)
	+j_1(L_3)\Delta t(2)\}
	}.
$$
Note that we do not only need to change the time slice that we are 
interacting with, but we also need to propagate the atomic observables 
that define the interaction.

\begin{remark}
Once again it is useful to have in mind the example of the classical 
Markov chain $X_{n}=f(X_{n-1},\xi_n)$.  Suppose $X_n$ takes values in 
$\mathbb{R}^k$; if we are looking e.g.\ at the first component $X_n^1$, 
then $X_n^1=f(X_{n-1}^1,X_{n-1}^2,\ldots,X_{n-1}^k,\xi_n)$.  Hence the
transformation $f(\cdot,X_{n-1}^2,\ldots,X_{n-1}^k,\xi_n)$ of $X_{n-1}^1$
depends not only on the noise ``time slice'' $\xi_n$, but also on other 
``observables'' of the system $X_{n-1}^2,\ldots,X_{n-1}^k$ which must be 
evaluated at time $n-1$.  Similarly our quantum mechanical interaction $M$ 
depends at each time step $l$ on the $l$th time slice in the field as well 
as on the atomic observables in the previous time step $l-1$.
\end{remark}

Proceeding in a similar fashion, we find that any atomic observable 
$X\in\mathscr{M}$ has evolved to the observable 
$j_l(X)=U(l)^*(X\otimes I)U(l)\in\mathscr{M}\otimes\mathscr{W}_k$ at time 
$l$, where
	\nomenclature{$j_l(X)$}{Time evolution of $X$}
	\nomenclature{$U(l)$}{Interaction unitary}
\begin{equation}
	U(l)=\overrightarrow{\prod}_{i=1}^l M(i)=
		M(1)M(2)\cdots M(l),\qquad
	U(0)=I,
\label{interaction-U}
\end{equation}
is a \emph{repeated interaction} with the one time step interaction
	\nomenclature{$M(l)$}{Single time step interaction}
\begin{equation}
	M(l)=
	e^{-i\{j_{l-1}(L_1)\Delta\Lambda(l)
	+j_{l-1}(L_2)\Delta A^*(l)
	+j_{l-1}(L_2^*)\Delta A(l)
	+j_{l-1}(L_3)\Delta t(l)\}
	}.
\label{interaction-M}
\end{equation}
The map $j_l(\cdot)$ defines the \emph{time evolution} or \emph{flow} of 
the atomic observables.  The connection with dynamical systems theory can 
be made explicit: if we define $U(i,l)=M(i+1)M(i+2)\ldots M(l)$, $i<l$, 
and $J_{i,l}(X)=U(i,l)^*XU(i,l)$ for $X\in\mathscr{M}\otimes\mathscr{W}_k$, 
then $J_{i,l}(\cdot)$ is a two-parameter group of transformations of the 
algebra $\mathscr{M}\otimes\mathscr{W}_k$ (i.e.\ for $i<l<r$  
$J_{l,r}(J_{i,l}(\cdot))=J_{i,r}(\cdot)$, 
$J_{i,r}(J_{i,l}^{-1}(\cdot))=J_{l,r}(\cdot)$, etc.); 
thus we have truly defined a discrete time dynamical system in the 
dynamical systems sense.  We will not need this in the following, however.

There is a different representation of the repeated interaction unitary
$U(l)$ which is often more convenient.  Introduce the unitaries
	\nomenclature{$M_l$}{Single time step interaction, time reversed}
\begin{equation}
\label{eq M}
	M_l=
	e^{-i\{L_1\Delta\Lambda(l)
	+L_2\Delta A^*(l)
	+L_2^*\Delta A(l)
	+L_3\Delta t(l)\}
	}
\end{equation}
which, in contrast to $M(l)$, depend on the \emph{initial} observables 
$L_{1,2,3}$ rather than the time-evolved observable $j_{l-1}(L_{1,2,3})$.
From $j_{l-1}(L_{1,2,3})=U(l-1)^*(L_{1,2,3}\otimes I)U(l-1)$, it follows 
immediately that $M(l)=U(l-1)^*M_lU(l-1)$.  But then
$$
	U(l)= U(l-1)M(l)=U(l-1)U(l-1)^*M_lU(l-1)=M_lU(l-1),
$$
so we obtain the expression
$$
	U(l)=\overleftarrow{\prod}_{i=1}^l M_i=
		M_lM_{l-1}\cdots M_2M_1.
$$
Hence we can use the $M_l$'s rather than $M(l)$'s if we reverse the time 
ordering in the definition of $U(l)$.  This is often convenient because 
$M_l$ is only a function of the initial system and of the $l$th time 
slice, whereas $M(l)$ depends on the entire history of the field up to 
time slice $l$.

The choice of $L_{1,2,3}\in\mathscr{M}$ determines the nature of the 
physical interaction between the atom and the field.  For different atoms 
or for different experimental scenarios this interaction can take very 
different forms; nonetheless the theory that we are about to develop can 
be set up in quite a general way for a large class of interactions.
In section \ref{sec:examplesrepint} we will introduce two specific 
examples of physical interactions which we will use throughout to 
illustrate the general theory.

\subsubsection*{Observations}

Up to this point we have concentrated on the time evolution $j_l(X)$ of an
atomic observable $X\otimes I$; however, the observables corresponding to
time slices of the field are also modified by the interaction with the
atom in a completely analogous way.  Recall, for example, that before the
interaction the number of photons in time slice $l$ is represented by the
observable $I\otimes\Delta\Lambda(l)$.  After the interaction this
quantity is represented by $\Delta Y^\Lambda(l)=U(l)^*(I\otimes
\Delta\Lambda(l))U(l)$, so that $\iota(\Delta Y^\Lambda(l))$ is precisely
what we observe at time $l$ if a photodetector is used to measure the
field.  Similarly, homodyne detection measures
$\Delta Y^X(l)=U(l)^*(I\otimes\Delta A(l)+I\otimes\Delta A^*(l))U(l)$ at
time $l$.  We will use the generic notation $\Delta 
Y(l)=U(l)^*\Delta Z(l)U(l)$, where $\Delta Y=\Delta Y^\Lambda$ if 
$\Delta Z=I\otimes\Delta\Lambda$ and $\Delta Y=\Delta Y^X$ if 
$\Delta Z=I\otimes(\Delta A+\Delta A^*)$.
	\nomenclature{$Y(l)$}{Observation process}
	\nomenclature{$Z(l)$}{Field process with $Y(l)=U(l)^*Z(l)U(l)$}

The first question we should ask ourselves is: does this procedure make
sense?  After all, in the laboratory we measure \emph{classical}
stochastic processes, so that in order for our theory to be consistent the
random variables $\iota(\Delta Y(l))$, $l=1,\ldots,k$ must define a
classical stochastic process on some fixed probability space
$(\Omega_k,\mathcal{F}_k,\mathbf{P}_k)$---in other words, $\mathscr{Y}_k=
{\rm alg}\{\Delta Y(l):l=1,\ldots,k\}$ must be a {\it commutative} 
algebra.  Let us verify that this is indeed the case.  The basic insight 
we need is the following observation: $[M_l,\Delta Z(i)]=0$ for $l>i$ 
(this is easily verified by inspection), so that for $j>l$
\begin{multline*}
	U(j)^*\Delta Z(l)U(j)=
	\left(\overrightarrow{\prod}_{i=1}^j M_i^*\right)
	\Delta Z(l)
	\left(\overleftarrow{\prod}_{i=1}^j M_i\right)
	\\ =
	\left(\overrightarrow{\prod}_{i=1}^l M_i^*\right)
	\Delta Z(l)
	\left(\overleftarrow{\prod}_{i=1}^l M_i\right)=
	U(l)^*\Delta Z(l)U(l).	
\end{multline*}
This is not unexpected: it is simply a statement of the fact that
the $l$th time slice of the field interacts exactly once with the atom, 
viz.\ at the $l$th time step (note that $U(j)^*\Delta Z(l)U(j)=\Delta 
Z(l)$ for $j<l$), for any $l$.  But then for $j>l$
\begin{multline*}
	[\Delta Y(l),\Delta Y(j)]=
	[U(l)^*\Delta Z(l)U(l),U(j)^*\Delta Z(j)U(j)]
	\\ =[U(j)^*\Delta Z(l)U(j),U(j)^*\Delta Z(j)U(j)]=
	U(j)^*[\Delta Z(l),\Delta Z(j)]U(j)=0,
\end{multline*}
so that clearly $\mathscr{Y}_k$ is commutative.  This 
\emph{self-nondemolition} property guarantees that the observations can 
indeed be interpreted as a classical stochastic process through the 
spectral theorem, a crucial requirement for the physical interpretation of 
the theory.

Previously we introduced the filtration $\mathscr{C}_l=
\mathrm{alg}\{\Delta Z(i):i=1,\ldots,l\}$, $l=1,\ldots,k$
	\nomenclature{$\mathscr{C}_l$}{Filtration of $Z$}
corresponding to the photodetection or homodyne detection of the field in 
absence of interaction with the atom.  Similarly we now define the 
filtration $\mathscr{Y}_l=\mathrm{alg}\{\Delta Y(i):i=1,\ldots,l\}$, 
$l=1,\ldots,k$
	\nomenclature{$\mathscr{Y}_l$}{Observation filtration}
which represents the information contained in observations of the field
(after it has interacted with the atom) up to the $l$th time step.  
Applying the spectral theorem to $\mathscr{Y}_k$, we obtain a classical 
probability space with the corresponding observation process $\Delta 
y_l=\iota(\Delta Y(l))$ and filtration $\iota(\mathscr{Y}_l)=
\ell^\infty(\mathcal{Y}_l)$
	\nomenclature{$\mathcal{Y}_l$}{Classical observation filtration}
	\nomenclature{$y_l$}{Classical observation process}
(note that we then have $\mathcal{Y}_l=
\sigma\{\Delta y_1,\ldots,\Delta y_l\}$).  The following fact will be 
useful later on:
$$
	\mathscr{Y}_l=U(j)^*\mathscr{C}_lU(j)=
	\{U(j)^*XU(j):X\in\mathscr{C}_l\}\qquad\forall\,j\ge l.
$$
The proof is identical to the proof of the self-nondemolition property.

We conclude this section with the demonstration of an important property 
of the models under investigation, the \emph{nondemolition property}, to 
which we have already alluded in section \ref{sec:rvini}.  Choose any
$X\in\mathscr{M}$.  Then for any $l\le i$
\begin{multline*}
	[\Delta Y(l),j_i(X)]=
	[U(l)^*\Delta Z(l)U(l),U(i)^*(X\otimes I)U(i)]
	\\ =[U(i)^*\Delta Z(l)U(i),U(i)^*(X\otimes I)U(i)]=
	U(i)^*[\Delta Z(l),X\otimes I]U(i)=0. 
\end{multline*} 
Evidently any atomic observable in time step $i$ commutes with the
entire history of observations up to that time.  Thus in principle we
could decide at any time to stop observing the field and we would still be
able to measure any atomic observable: the joint distribution of $\Delta 
Y(1),\ldots,\Delta Y(l)$ and $j_l(X)$ is defined for any self-adjoint
$X\in\mathscr{M}$, as ${\rm alg}\{\Delta Y(1),\ldots,\Delta Y(l),j_l(X)\}$ 
is a commutative $^*$-algebra, despite the fact that $j_l(X_1)$ and 
$j_l(X_2)$ need not commute for $X_1\ne X_2$.  This enables us to 
meaningfully define the \emph{conditional expectation} of $j_l(X)$ with 
respect to the observations $\Delta Y(i)$, $i=1,\ldots,l$ for {\it any} 
$X\in\mathscr{M}$, which we will do in section \ref{conditional 
expectation}. The nondemolition property provides us with a sensible way
of estimating a whole set observables, despite the fact that these do not
commute with each other, as every observable in this set separately
commutes with the observation history on which we are basing the estimate.  
If the latter were not the case the estimates would have no physical
relevance:  one would not be able to devise an experiment, even in 
principle, which could verify the predictions of such an estimator.
Similarly, it would be meaningless to try to control atomic observables 
which do not commute with the observations on which we have based our 
feedback control law.  The nondemolition property avoids such problems and 
is thus of fundamental importance for quantum filtering and feedback 
control.

\subsection{Examples}\label{sec:examplesrepint}

The interaction matrices $L_{1,2,3}$ determine the nature of the physical 
interaction between the atom and the field.  Though we will set up the 
theory in the following sections for arbitrary $L_{1,2,3}$, we will 
repeatedly demonstrate the theory using the following two examples.  The 
examples illustrate two common experimental scenarios---spontaneous 
emission and dispersive interaction---and we will be able to study these 
examples explicitly through numerical simulations.

\subsubsection*{Spontaneous emission}

The \emph{spontaneous emission} scenario is obtained when an atom sits 
undisturbed in the vacuum at zero temperature.  If the energy of the atom 
is minimal $\mathbb{P}(H)=-\hbar\omega_0/2$, then it will remain this way 
for all time.  Otherwise, the energy of the atom decays to its minimal 
value.  A photodetector measuring the field would see exactly one emitted 
photon at some randomly distributed time: hence the name spontaneous 
emission.  We will reproduce this behavior through numerical simulations 
once we have introduced the corresponding filtering equations.

The spontaneous emission model is characterized by the interaction matrices
$$
	L_1 = L_3 = 0,\qquad
	L_2 = i\sqrt{2\kappa}\,\sigma_-,
$$
where $\kappa$ is the spontaneous emission rate.  For simplicity, we will 
always set $2\kappa=1$.  This is a discretized version of the well known 
Wigner-Weisskopf model.

\begin{remark} In principle, we could calculate the possible observation
sample paths and their probabilities already at this point.  After all, we
would only need to simultaneously diagonalize the matrices $\Delta Y(i)$,
$i=1,\ldots,k$ obtained using this explicit choice for $L_{1,2,3}$.  We
will have a much more efficient way of calculating such sample paths,
however, once we have introduced the filtering equations: these will allow
us to simulate typical sample paths using a Monte Carlo method.  Hence we 
postpone numerical investigation of our examples until then.
\end{remark}

\subsubsection*{Dispersive interaction}

In the spontaneous emission scenario, the atom radiates its excess energy 
by emitting a photon with energy $\hbar\omega_0$ (i.e.\ with 
frequency $\omega_0$)---this is precisely the difference between the two 
values that the energy $H$ can take.  It is possible to suppress the 
spontaneous emission significantly by ``discouraging'' the atom to emit 
photons with frequency $\omega_0$; for example, we can place the atom in 
an optical cavity whose resonant frequency is far detuned from $\omega_0$.  
In this \emph{dispersive regime}, the (effective) interaction between the 
atom and the field is rather different in nature.  An atom which has 
energy $+\hbar\omega_0/2$ will shift the phase of the output light by a 
small amount in one direction, while an atom with energy 
$-\hbar\omega_0/2$ will shift the phase by the same amount in the other 
direction.  Such phase shifts can be observed using a homodyne detection 
setup, and the resulting photocurrent thus carries information on the
energy of the atom.

The dispersive interaction model is characterized by the interaction 
matrices
$$
	L_1 = L_3 = 0,\qquad
	L_2 = i\sqrt{g}\,\sigma_z,
$$
where $g$ is the interaction strength.  For simplicity, we will
always set $g=1$.

\section{Conditional expectations and the filtering problem}
\label{conditional expectation}

Within the context of the repeated interaction models introduced in the
previous section, we have established that the observation history up to
time step $l$ is compatible with any atomic observable at time step $l$.  
In this section we will show how to estimate atomic observables based on
the observation history.  Because only commuting observables are involved,
the corresponding theory can simply be ``lifted'' from classical
probability theory using the spectral theorem.  This is precisely what we
will do.

\subsection{A brief reminder}\label{sec reminder condex}

We begin by briefly recalling how classical conditioning works.  Let 
$(\Omega,\mathcal{F},\mathbf{P})$ be a probability space and let 
$P,Q\in\mathcal{F}$ be a pair of events.  The \emph{conditional probability} of 
$P$ given $Q$ is given by
$$
	\mathbf{P}(P|Q)=\frac{\mathbf{P}(P\cap Q)}{\mathbf{P}(Q)},
$$
provided $\mathbf{P}(Q)\ne 0$.  One could interpret this quantity roughly 
as follows: generate a large number of samples distributed according to 
the measure $\mathbf{P}$, but discard those samples for which the event 
$Q$ is false.  Then $\mathbf{P}(P|Q)$ is the fraction of the remaining 
samples for which the event $P$ is true.

Now suppose $f:\Omega\to\{f_1,\ldots,f_m\}$ and 
$g:\Omega\to\{g_1,\ldots,g_n\}$ are measurable random variables that take 
a finite number of values.  The \emph{conditional expectation} of $f$ given 
$g$ is the random variable defined by
\begin{equation}\label{eq classical conditional expectation}
	\mathbf{E_P}(f|g)(\omega)=
	\sum_{j=1}^n \chi_{g^{-1}(g_j)}(\omega)
	\sum_{i=1}^m f_i\,\mathbf{P}(f^{-1}(f_i)|g^{-1}(g_j)),
	\qquad \omega\in\Omega.	
\end{equation}
In fact, we can consider $\mathbf{E_P}(f|g)(\omega)$ to be a function of 
$g$:
$$
	\mathbf{E_P}(f|g)(\omega)=F(g(\omega)),\qquad
	F:g_j\mapsto
	\sum_{i=1}^m f_i\,\mathbf{P}(f^{-1}(f_i)|g^{-1}(g_j)).
$$
To interpret this quantity, again we generate a large number of samples 
but now we divide these up into disjoint subsets corresponding to the 
value taken by $g$.  Now average $f$ over each of these subsets.  Then 
$F(g_j)$ is the average of $f$ over the subset of samples on which $g$ 
takes the value $g_j$, so that $\mathbf{E_P}(f|g)=F(g)$ is simply 
the expectation of $f$ given that we know $g$.  
Note that we did not define the quantity $\mathbf{P}(P|Q)$ for the case
$\mathbf{P}(Q)=0$;  hence the expressions above for $\mathbf{P}(P|Q)$ are
properly defined everywhere on $\Omega$ except on
$\Omega^0=\bigcup_{j:\mathbf{P}(g^{-1}(g_j))=0}g^{-1}(g_j)\in\mathcal{F}$. We
allow $\mathbf{E_P}(f|g)$ to take arbitrary values on $\Omega^0$ and call
any such random variable a \emph{version} of the conditional expectation
of $f$ with respect to $g$.  There is no loss in doing this: all versions 
of the conditional expectation coincide with unit probability.

The use of these expressions is rather limited, as they do not extend
beyond random variables $g$ that take a finite number of values.  One of
Kolmogorov's fundamental insights in the axiomatization of probability
theory was his abstract definition of the conditional expectation, which
works in the most general settings.  In this article we only use finite
state random variables, but even so we find it significantly easier to use
the abstract definition than these clumsy explicit expressions.  Let us
thus briefly recall the general setting, as can be found in any
probability textbook (e.g.\ \cite{Wil91}).

Let $\mathcal{G}\subset\mathcal{F}$ be a $\sigma$-algebra on which we want to 
condition---i.e., these are the events which we know to be true or false 
due to our observations.  Let $f$ be a $\mathcal{F}$-measurable random 
variable.  Then any $\mathcal{G}$-measurable random variable $\hat f$ that 
satisfies
\begin{equation}\label{eq classical abstract characterization}
	\mathbf{E_P}(\hat fg)=\mathbf{E_P}(fg)\quad
	\forall\,\mathcal{G}\mbox{-measurable (bounded) }g 
\end{equation} 
is a version of the conditional expectation $\hat f=
\mathbf{E_P}(f|\mathcal{G})$.  The conditional expectation is by
construction a function of what we have observed (it is
$\mathcal{G}$-measurable), as it should be.  It is easily verified that
the explicit expression Eq.\ (\ref{eq classical conditional expectation})  
satisfies the abstract definition with $\mathcal{G}=\sigma\{g\}$.  In
fact, in the case where $\Omega$ takes a finite number of values we can
always express $\mathbf{E_P}(f|\mathcal{G})$ as in Eq.\ (\ref{eq classical
conditional expectation}) by constructing $g$ specifically so that it
generates $\mathcal{G}$: expressing $g$ in terms of the smallest
measurable partition of $\Omega$ as in the proof of Lemma
\ref{lem:onetoonespectral}, we only need to make sure that $g$ takes a
different value on each set in that partition.  Hence there always exists
a $\hat f$ that satisfies the abstract definition, and it is not difficult
to prove that all versions coincide with unit probability---the proof of
the latter is the same as the one given below in the quantum case.  
Elementary properties of the conditional expectation are listed in Table
\ref{tab:condex}.

\begin{table} 
\caption{Elementary properties of the conditional expectation in classical
(and quantum) probability; see e.g.\ {\rm\cite{Wil91}}.  Equalities are
always in the sense of versions, e.g.\ $\mathbf{E_P}(f_1|\mathcal{G})+
\mathbf{E_P}(f_2|\mathcal{G})=\mathbf{E_P}(f_1+f_2|\mathcal{G})$ means
that any version of $\mathbf{E_P}(f_1|\mathcal{G})+
\mathbf{E_P}(f_2|\mathcal{G})$ is a version of
$\mathbf{E_P}(f_1+f_2|\mathcal{G})$.}

\begin{center} \footnotesize
\begin{tabular}{|l|l|} 
\hline
Property & Description \\
\hline
Linearity  & $\alpha\,\mathbf{E_P}(f_1|\mathcal{G})
	+\beta\,\mathbf{E_P}(f_2|\mathcal{G})=
	\mathbf{E_P}(\alpha f_1+\beta f_2|\mathcal{G})$. \\
Positivity & If $f\ge 0$ then $\mathbf{E_P}(f|\mathcal{G})\ge 0$ a.s. \\
Invariance of expectation & $\mathbf{E_P}(\mathbf{E_P}(f|\mathcal{G}))=
	\mathbf{E_P}(f)$. \\
Module property & If $g$ is $\mathcal{G}$-measurable,
	$\mathbf{E_P}(fg|\mathcal{G})=g\,\mathbf{E_P}(f|\mathcal{G})$ a.s. 
\\
	& In particular, $\mathbf{E_P}(g|\mathcal{G})=
		g\,\mathbf{E_P}(1|\mathcal{G})=g$ a.s. \\
Tower property & If $\mathcal{H}\subset\mathcal{G}$,
	$\mathbf{E_P}(\mathbf{E_P}(f|\mathcal{G})|\mathcal{H})=
	\mathbf{E_P}(f|\mathcal{H})$. \\
	Independence & If $\sigma\{f,\mathcal{H}\}$ is independent of 
	$\mathcal{G}$, $\mathbf{E_P}(f|\sigma\{\mathcal{H},\mathcal{G}\})=
	\mathbf{E_P}(f|\mathcal{H})$ a.s. \\
\hline
\end{tabular}
\end{center}
\label{tab:condex}
\end{table}

We conclude by recalling the geometric interpretation of the conditional
expectation.  For simplicity, take $\Omega=\{1,\ldots,N\}$ and fix the
$\sigma$-algebra $\mathcal{G}\subset\mathcal{F}$.  Clearly
$\ell^\infty(\mathcal{F})$ is a finite-dimensional linear space and
$\ell^\infty(\mathcal{G})\subset \ell^\infty(\mathcal{F})$ is a linear
subspace.  Moreover $\langle f,g\rangle_\mathbf{P}= \mathbf{E_P}(f^*g)$ is
a pre-inner product: it satisfies all the conditions of the inner product,
except $\langle f,f\rangle_\mathbf{P}=0$ iff $f=0$ $\mathbf{P}$-a.s.\
rather than iff $f(\omega)=0$ $\forall\,\omega\in\Omega$.  Hence under
$\langle\cdot,\cdot\rangle_\mathbf{P}$, $\ell^\infty(\mathcal{F})$ is a
pre-inner product space.  With this interpretation, it is evident from
Eq.\ (\ref{eq classical abstract characterization}) that
$\mathbf{E_P}(f|\mathcal{G})$ is simply the \emph{orthogonal projection}
of $f\in\ell^\infty(\mathcal{F})$ onto $\ell^\infty(\mathcal{G})$.  It follows
from the projection theorem that $\mathbf{E_P}(f|\mathcal{G})$ is the
a.s.\ unique random variable $\hat f\in\ell^\infty(\mathcal{G})$ that
minimizes the mean-square error $\langle f-\hat f,f-\hat
f\rangle_\mathbf{P}= \mathbf{E_P}(|f-\hat f|^2)$. This gives the
conditional expectation, in addition to the \emph{probabilistic}
interpretation above, a firm \emph{statistical} interpretation as the
estimator that minimizes the mean square error.

\subsection{Quantum conditional expectation}

Let $(\mathscr{B},\mathbb{P})$ be a commutative quantum probability space;
by the spectral theorem, it is equivalent to some classical probability
space $(\Omega,\mathcal{F},\mathbf{P})$.  Moreover, the $^*$-subalgebra
$\mathscr{C}\subset\mathscr{B}$ defines through the spectral theorem a
$\sigma$-subalgebra $\mathcal{G}\subset\mathcal{F}$ (as $\iota(\mathscr{C})=
\ell^\infty(\mathcal{G})\subset\ell^\infty(\mathcal{F})$).  Because the 
classical and quantum probability models are completely equivalent, we can 
simply lift the definition of the conditional expectation to the algebraic 
level:
$$
	\mathbb{P}(X|\mathscr{C})=
	\iota^{-1}(\mathbf{E_P}(\iota(X)|\mathcal{G}))
	\qquad
	\forall X\in\mathscr{B}.
$$
This is nothing but a re-expression of the conditional expectation 
$\mathbf{E_P}(f|\mathcal{G})$ in the language of $^*$-algebras.  In fact, 
we can go one step further and directly translate the abstract definition
Eq.\ (\ref{eq classical abstract characterization}) into algebraic 
language: any $\hat X\in\mathscr{C}$ that satisfies
$\mathbb{P}(XC)=\mathbb{P}(\hat XC)$ $\forall\,C\in\mathscr{C}$ is a 
version of $\mathbb{P}(X|\mathscr{C})$.  We emphasize again: by the 
spectral theorem, observables $X\in\mathscr{B}$ {\it are} just random 
variables and $\mathbb{P}(X|\mathscr{C})$ {\it is} just the ordinary 
conditional expectation.  The spectral theorem is a powerful tool indeed!

Usually we do not start with a commutative quantum probability space.  
Rather, we begin with a noncommutative space $(\mathscr{A},\mathbb{P})$
and choose a commutative subalgebra $\mathscr{C}\subset\mathscr{A}$
corresponding to the observations---keep in mind the example of the space
$(\mathscr{M}\otimes\mathscr{W}_k,\rho\otimes\phi^{\ten k})$ with the
observations $\mathscr{Y}_l\subset\mathscr{M}\otimes\mathscr{W}_k$ at time
step $l$.  As we have seen, there could be many elements $X\in\mathscr{A}$
that commute with every element in $\mathscr{C}$, but that do not commute
with each other.  The set $\mathscr{C}'=\{A\in\mathscr{A}:[A,C]=0~\forall
C\in\mathscr{C}\}$ is called the \emph{commutant} of $\mathscr{C}$ in
$\mathscr{A}$.  It is easily verified that $\mathscr{C}'$ is always a
$^*$-algebra, but generally $\mathscr{C}'$ is not commutative.  
Nonetheless we can naturally extend our previous definition of
$\mathbb{P}(\cdot|\mathscr{C})$ from the commutative algebra 
$\mathscr{B}$ to the case $\mathscr{B}=\mathscr{C}'$.

\begin{definition}[\bf Quantum conditional expectation]
\label{de conditional expectation}
Let $(\mathscr{A}, \mathbb{P})$ be a quantum probability space, and let
$\mathscr{C}$ be a commutative subalgebra of $\mathscr{A}$.
Then the map $\mathbb{P}(\cdot |\mathscr{C}):\mathscr{C}' \to \mathscr{C}$ 
is called a version of the conditional expectation from $\mathscr{C}'$ 
onto $\mathscr{C}$ if $\mathbb{P}(\mathbb{P}(A|\mathscr{C})C) = 
\mathbb{P}(AC)$ $\forall\,A\in\mathscr{C}',~C\in \mathscr{C}$.
\end{definition}

How should we interpret this definition?  Let $A=A^*\in\mathscr{C}'$ be an
observable, and suppose we are interested in $\mathbb{P}(A|\mathscr{C})$.
The definition of $\mathbb{P}(A|\mathscr{C})$ only involves the operators
$A$ and $C\in\mathscr{C}$, so for the purposes of evaluating
$\mathbb{P}(A|\mathscr{C})$ we don't need the entire algebra
$\mathscr{C}'$: we might as well confine ourselves to
$\mathscr{C}_A={\rm alg}\{A,C:C\in\mathscr{C}\}\subset\mathscr{C}'$---but
$\mathscr{C}_A$ is commutative!  Hence for any $A=A^*\in\mathscr{C}'$,
Definition \ref{de conditional expectation} simply reduces to the
commutative definition we gave before.  For $A\ne A^*$, we can always
write $A=P+iQ$ where $P=(A+A^*)/2$ and $Q=-i(A-A^*)/2$ are self-adjoint,
and we define $\mathbb{P}(A|\mathscr{C})=\mathbb{P}(P|\mathscr{C})+
i\,\mathbb{P}(Q|\mathscr{C})$.  We do this mostly for computational
convenience (e.g., to ensure that the conditional expectation is
linear); at the end of the day, only conditional expectations of
observables have a meaningful interpretation.

By now it should not come as a surprise that the conditional expectation
can only be defined on the commutant $\mathscr{C}'$.  In section \ref{sec
reminder condex} we gave an ``experimental procedure'' for calculating the
conditional expectation of a random variable $f$: generate a large number
of samples from the joint distribution of $f$ and the observations which
we are conditioning on, separate the samples into subsets corresponding to
the different observations, and average $f$ within each of these subsets.
But quantum mechanics tells us that it is fundamentally impossible to
devise such an experiment when the observable $A$ which we want to
condition is incompatible with the observations!  The commutant
$\mathscr{C}'$ is thus precisely the right notion for a meaningful
definition of the quantum conditional expectation\footnote{
	That is not to say that this is the only definition used in the
	quantum probability literature; in fact, a more general 
	definition, of which our definition is a special case, is very 
	common \cite{Tak71}.  Such a definition is very different in 
	spirit, however, and is not motivated by the probabilistic/statistical
	point of view needed for filtering and feedback control.	
}.

Having introduced an abstract definition, let's do a little ground work.

\begin{lemma}
	$\mathbb{P}(A|\mathscr{C})$, $A\in\mathscr{C}'$ exists and is 
	unique with probability one, i.e.\ any two versions $P$ and $Q$ of 
	$\mathbb{P}(A|\mathscr{C})$ satisfy $\mathbb{P}((P-Q)^*(P-Q))=0$.
\end{lemma}

\begin{proof}
We first introduce some notation.  $\langle F,G\rangle_\mathbb{P}=
\mathbb{P}(F^*G)$ defines a pre-inner product on $\mathscr{C}'$, thus 
turning it into a pre-inner product space.  The associated seminorm is 
denoted by $\|F\|_\mathbb{P}=(\langle F,F\rangle_\mathbb{P})^{1/2}$.

{\it Existence.} For any self-adjoint $A=A^*\in\mathscr{C}'$,
define $\mathscr{C}_A={\rm alg}\{A,C:C\in\mathscr{C}\}$.  The spectral 
theorem then provides a probability space 
$(\Omega_A,\mathcal{F}_A,\mathbf{P}_A)$ and a $^*$-isomorphism 
$\iota_A:\mathscr{C}_A\to \ell^\infty(\mathcal{F}_A)$, and we write
$\iota_A(\mathscr{C})=\ell^\infty(\mathcal{G}_A)$.  For every such $A$, 
define $\mathbb{P}(A|\mathscr{C})=\iota_A^{-1}(a)$ where $a$ is some 
version of $\mathbf{P}_A(\iota_A(A)|\mathcal{G}_A)$.  For non-selfadjoint 
$A\in\mathscr{C}'$, define $\mathbb{P}(A|\mathscr{C})=
\mathbb{P}(P|\mathscr{C})+i\,\mathbb{P}(Q|\mathscr{C})$ with 
$P=(A+A^*)/2$ and $Q=-i(A-A^*)/2$.  It is easily verified that the map 
$\mathbb{P}(\cdot|\mathscr{C})$ thus constructed satisfies Definition 
\ref{de conditional expectation}.

{\it Uniqueness with probability one.} 
Let $P$ and $Q$ be two versions of $\mathbb{P}(A|\mathscr{C})$.  
From Definition \ref{de conditional expectation}, it follows 
that $\langle C,P-Q\rangle_\mathbb{P}=0$ for all $C\in\mathscr{C}$.  But 
$P-Q\in\mathscr{C}$, so $\langle P-Q,P-Q\rangle_\mathbb{P}=
\|P-Q\|_\mathbb{P}^2=0$.
\qquad
\end{proof}

Let us prove that the conditional expectation is the least mean square 
estimate.

\begin{lemma}
$\mathbb{P}(A|\mathscr{C})$ is the least mean 
square estimate of $A$ given $\mathscr{C}$, i.e. for any $C\in\mathscr{C}$
we have $\|A-\mathbb{P}(A|\mathscr{C})\|_\mathbb{P}\le\|A-C\|_\mathbb{P}$.
\end{lemma}

\begin{proof}
For any $C\in\mathscr{C}$ we have $\|A-C\|_\mathbb{P}^2 = 
\|A-\mathbb{P}(A|\mathscr{C})+ 
\mathbb{P}(A|\mathscr{C})-C\|_\mathbb{P}^2$.  
Now note that by Definition \ref{de conditional expectation}, we have
$\langle C,A-\mathbb{P}(A|\mathscr{C})\rangle_\mathbb{P} =
\mathbb{P}(C^*A)-\mathbb{P}(C^*\mathbb{P}(A|\mathscr{C}))=0$
for all $C \in \mathscr{C}$, $A \in \mathscr{C}'$, i.e.,
$A -\mathbb{P}(A|\mathscr{C})$ is orthogonal to $\mathscr{C}$.  In 
particular, $A-\mathbb{P}(A|\mathscr{C})$ is orthogonal to 
$\mathbb{P}(A|\mathscr{C})-C\in \mathscr{C}$, and we obtain
$\|A-C\|_\mathbb{P}^2 =\|A-\mathbb{P}(A|\mathscr{C})\|_\mathbb{P}^2 + 
\|\mathbb{P}(A|\mathscr{C})-C\|_\mathbb{P}^2$.  The result follows 
immediately.
\qquad
\end{proof}

The elementary properties listed in Table \ref{tab:condex} and their
proofs carry over directly to the quantum case.  For example, let us prove
linearity.  It suffices to show that $Z=\alpha\,\mathbb{P}(A|\mathscr{C})
+\beta\,\mathbb{P}(B|\mathscr{C})$ satisfies the definition of
$\mathbb{P}(\alpha A+\beta B|\mathscr{C})$, i.e.\ $\mathbb{P}(ZC)=
\mathbb{P}((\alpha A+\beta B)C)$ for all $C\in\mathscr{C}$.  But this is 
immediate from the linearity of $\mathbb{P}$ and Definition \ref{de 
conditional expectation}.  We encourage the reader to verify the remaining 
properties.

In Section \ref{reference probability} we will need to relate conditional
expectations with respect to different states to each other.  This is done
by the following Bayes-type formula.

\begin{lemma}[\bf Bayes formula]\label{lem Bayes}
Let $(\mathscr{A}, \mathbb{P})$ be a quantum probability space, and let
$\mathscr{C}$ be a commutative subalgebra of $\mathscr{A}$.
Let $V$ be an element in $\mathscr{C}'$ such that $V^*V>0$ and 
$\BB{P}(V^*V)=1$.  We can define a new state $\mathbb{Q}$ on 
$\mathscr{C}'$ by $\mathbb{Q}(X) = \mathbb{P}(V^*XV)$, and
\begin{equation*}
  	\mathbb{Q}(X|\mathscr{C}) 
	= \frac{\mathbb{P}(V^*XV|\mathscr{C})}
  		{\BB{P}(V^*V|\mathscr{C})}
	\qquad\forall\, X\in \mathscr{C}'.
\end{equation*}
\end{lemma}

\begin{proof}
Let $K$ be an element of $\mathscr{C}$. For any $X \in \mathscr{C}'$
we can write
\begin{multline*}
  \mathbb{P}(\mathbb{P}(V^*XV|\mathscr{C})K) = 
  \mathbb{P}(V^*XKV) = 
  \mathbb{Q}(XK) = \mathbb{Q}(\mathbb{Q}(X|\mathscr{C})K) = \\
  \mathbb{P}(V^*V\mathbb{Q}(X|\mathscr{C})K) = 
  \mathbb{P}(\mathbb{P}
  (V^*V\mathbb{Q}(X|\mathscr{C})K|\mathscr{C})) =
  \mathbb{P}(\mathbb{P}(V^*V|\mathscr{C})
  \mathbb{Q}(X|\mathscr{C})K),
\end{multline*}
and the result follows directly.
\qquad
\end{proof}

\subsection{Statement of the filtering problem}

Let us return to the two-level system coupled to the discretized
electromagnetic field, as described in section \ref{sec:repeatedinteraction}.
Recall that we observe the commutative process $\Delta Y(l)$, with either 
$Y=Y^\Lambda$ or $Y^X$ depending on what detector we choose to use, and 
that this process generates the commutative filtration $\mathscr{Y}_l$.  
Suppose we have been observing $Y$ up to and including time $l$.  What can 
we infer about the atom?

The nondemolition condition allows us to make sense of this question.
In the language of this section, nondemolition simply means that
$j_l(X)\in\mathscr{Y}_l'$ for any $l=0,\ldots,k$, $X\in\mathscr{M}$.  
Hence the conditional expectations
	\nomenclature{$\pi_l(X)$}{Conditional state (filtered estimate
		of $X$)}
\begin{equation*}
  \pi_l(X) = \mathbb{P}(j_l(X)|\mathscr{Y}_l), \qquad
	0 \le l \le k,~X \in \mathscr{M}
\end{equation*}
are well defined, and we could e.g.\ interpret $\pi_l(X)$ as the least 
mean square estimate of an atomic observable $X$ at time step $l$ given 
the observation history.

Though it is evident that the conditional expectations $\pi_l(X)$ are well
defined, it is quite another matter to calculate them explicitly.  In
principle we have a recipe which we could follow to perform this
calculation: find a basis in which $j_l(X)$ and $\Delta Y(i)$,
$i=1,\ldots,l$ are simultaneously diagonalized, apply an expression
similar to Eq.\ (\ref{eq classical conditional expectation}), and then
re-express the result in the original basis.  Such a procedure would be
extremely unpleasant, to say the least; clearly we are in need a better 
way of doing this.  Finding a useful explicit representation of the 
conditional expectations $\pi_l(\cdot)$ is known as the \emph{filtering 
problem}.

Sections \ref{sec doob} and \ref{reference probability} are dedicated to
the solution of the filtering problem, which we will do in two different
ways (and obtain two different expressions!)  Both methods give rise to
\emph{recursive} equations, that is, we will find \emph{filtering
equations} of the form $\pi_l(\cdot)=f(\pi_{l-1}(\cdot),\Delta Y(l))$.  
This is very convenient computationally as we don't need to remember the
history of observations: to calculate the current estimates
$\pi_l(\cdot)$, all we need to know are the corresponding estimates in the
previous time step $\pi_{l-1}(\cdot)$ and our observation in the current
time step $\Delta Y(l)$.  In other words, we can update our estimates in 
real time using only the current observation in each time step.

\subsection{A brief look ahead}
\label{sec:lookahead}

The filtering problem is very interesting in itself; its solution provides
us with an efficient way of estimating atomic observables from
the observation data.  We have an ulterior motive, however, for putting
such emphasis on this problem---filtering is a fundamental part of
stochastic control with partial observations, and hence plays a central
role in quantum feedback control.  This is the subject of sections
\ref{sec:fb}--\ref{sec:lyapunov}, and we postpone a detailed discussion 
until then.  Let us briefly sketch, however, why we might expect filtering 
to enter the picture.

In order to design a controller, we first have to specify what goal we are
trying to achieve.  Control goals are generally expressed in terms of
expectations over atomic observables.  For example, we might want to
maximize the expected $x$-dipole moment $\mathbb{P}(j_k(\mu_x))$ at the
terminal time $t(k)=T$, or we might try to steer the dipole moment to a
specific value $\alpha$ (e.g.\ by minimizing
$\mathbb{P}(j_k((\mu_x-\alpha)^2))$).  If we are not sure that we are
going to finish running the control through time $T$, perhaps we would
prefer to maximize the time-average expected dipole moment
$\tfrac{1}{k}\sum_{l=1}^k\mathbb{P}(j_l(\mu_x))$---etc.\ etc.  The problem
is that our control can depend on the entire observation history; the {\it
optimal} control can have quite a complicated dependence on the
observation history!

The key observation we need is that by an elementary property of the
conditional expectation,
$\mathbb{P}(j_l(X))=\mathbb{P}(\pi_l(X))$. Rather than expressing
our control goal directly in terms of the quantum model, we are
free to express the same goal in terms of the filter only, as though it
were the filter we wanted to control rather than the atom itself.  It will
turn out that optimal controls have a very simple dependence on the
filter: the control in time step $l$ can be expressed as some 
(deterministic) function of the filter at time $l-1$.  The filter 
$\pi_l(\cdot)$ is called an \emph{information state} for the control 
problem, as $\pi_l(\cdot)$ contains all the information we need to compute 
the feedback law.  This is discussed in detail in section 
\ref{sec:optimal-IS}.

The fact that optimal controls separate into a filtering step and a simple
feedback step is known as the {\it separation principle}.  The complicated
dependence of the control on the observation history is completely
absorbed by the filter (which is by construction only a function of the
observations).  By solving the filtering problem in a recursive manner, we
obtain an efficient way of calculating the optimal controls.  This amply
motivates, beside its intrinsic interest, our interest in the filtering
problem.

\section{Discrete quantum stochastic calculus}\label{discrete QSC}

In the previous sections we have provided an introduction to quantum
probability, introduced a toy model for the interaction of an atom with
the electromagnetic field, and introduced conditional expectations and the
filtering problem in the quantum setting.  At this point we change gears,
and concentrate on developing machinery that will help us to
actually solve the filtering problem (in two different ways) and 
ultimately to treat control problems associated with our model.

We begin in this section by rewriting the repeated interactions of section
\ref{sec:repeatedinteraction} as difference equations, and we develop a
``stochastic calculus'' to manipulate such equations.  In the discrete
setting this calculus is not a very deep result: it is merely a matrix
multiplication table, as can be seen from the proof below.  In continuous
time, however, difference equations become (quantum stochastic)
differential equations and the stochastic calculus is an indispensable
tool.  While in continuous time linear algebra is replaced by the much
more difficult functional analysis, the quantum stochastic calculus
nonetheless enables one to deal with the corresponding objects using
simple algebraic manipulations.  Bearing in mind our goal to mirror as
closely as possible the continuous theory, we will consistently use the
discrete calculus in what follows.

\subsection{The discrete quantum It\^o calculus} 

Let us fix once and for all the quantum probability space 
$(\mathscr{M}\otimes \mathscr{W}_k,\mathbb{P})$ with $\mathbb{P}=
\rho\otimes\phi^{\otimes k}$, which models an atom together with the 
electromagnetic field as in section \ref{sec:repeatedinteraction}.  The 
algebra $\mathscr{M}\otimes\mathscr{W}_k$ has a natural {\it 
noncommutative} filtration $\{\mathscr{B}_l\}_{l=1,\ldots k}$ defined by
	\nomenclature{$\mathscr{B}_l$}{Filtration in $\mathscr{M}\otimes\mathscr{W}_k$}
$$
   \mathscr{B}_0 = \mathscr{M}\otimes I^{\otimes k} \subset
   \mathscr{B}_1 = \mathscr{M}\otimes\mathscr{M}\otimes I^{\otimes (k-1)}
   \subset\cdots
   \subset \mathscr{B}_k = \mathscr{M}\otimes \mathscr{W}_k.
$$
The noncommutative algebra $\mathscr{B}_l$ contains all observables that 
are a function of the atom and the time slices of the field up to and 
including slice $l$.

A \emph{quantum process} is a map from $\{0,1,\ldots, k\}$ to
$\mathscr{M}\otimes\mathscr{W}_k$. A quantum process $L$ is called
\emph{adapted} if $L(i)\in\mathscr{B}_i$ for every $0 \le i \le k$. It is
called \emph{predictable} if $L(i)\in\mathscr{B}_{i-1}$ for every $1 \le i
\le k$.  For a quantum process $L$ we define $\Delta L$ as $\Delta L(i) = 
L(i) - L(i-1)$ for $1 \le i \le k$.  These definitions are similar to 
those used in classical probability \cite{Wil91}.
	\nomenclature{$\Delta L(l)$}{$L(l)-L(l-1)$}

\begin{definition}{\bf(Discrete quantum stochastic integral)}
Let $L$ be a predictable quantum process and let $M$ be one of the 
processes $A$, $A^*$, $\Lambda$ or $t$.  The \emph{transform of $M$ by 
$L$}, denoted $L\cdot M$ is the adapted quantum process defined by   
\begin{equation*}
  L\cdot M (l) = \sum_{i=1}^l L(i)\Delta M(i), \quad 1 \le l \le k,
	\qquad L\cdot M(0)=0.
\end{equation*}
$L\cdot M$ is also called the \emph{discrete quantum stochastic integral 
of $L$ with respect to $M$}.   
\end{definition}

Note that for any discrete quantum stochastic integral as above, we can 
equivalently write $X = L\cdot M \iff \Delta X = L \Delta M$.  Note also 
that $L(i)$ and $\Delta M(i)$ commute by construction, so that we can 
always exchange the order of the (predictable) integrands and the 
increments $\Delta\Lambda$, $\Delta A$, $\Delta A^*$ and $\Delta t$.  We 
can now state the main theorem of our ``stochastic calculus''.

\begin{theorem}\label{thm Ito rule}{\bf (Discrete quantum It\^o rule)}
Let $L_1$ and $L_2$ be predictable quantum processes, $M_1,M_2\in
\{A,A^*,\Lambda,t\}$, and let $X$ and $Y$ be the transforms $L_1\cdot M_1$ 
and $L_2 \cdot M_2$, respectively.  Then
\begin{equation*}
	\Delta (XY)(l) =   
	X(l-1)\,\Delta Y(l) 
	+ \Delta X(l)\, Y(l-1) 
	+ \Delta X(l)\,\Delta Y(l),
\end{equation*}
where $\Delta X(l)\,\Delta Y(l)$ should be evaluated according 
to the discrete quantum It\^o table:
\vskip.2cm
\begin{center}
  \begin{tabular} {l|llll}
  $\Delta X\backslash\Delta Y$
   & $\Delta A$ & $\Delta \Lambda$ & $\Delta A^*$ & $\Delta t$\\
  \hline 
  $\Delta A$ & $0$ & $\Delta A$ & $\Delta t -\lambda^2\Delta\Lambda$ & $\lambda^2\Delta A$ \\
  $\Delta \Lambda$ & $0$ & $\Delta\Lambda$ & $\Delta A^*$ & $\lambda^2 \Delta \Lambda$  \\
  $\Delta A^*$ & $\lambda^2\Delta\Lambda$ & $0$ & $0$ & $\lambda^2\Delta A^*$ \\
  $\Delta t$ & $\lambda^2\Delta A$ & $\lambda^2\Delta\Lambda$ & $\lambda^2\Delta A^*$ & $\lambda^2\Delta t$
  \end{tabular} 
\end{center}
\vskip.2cm
For example, if $M_1 = A$ and $M_2 = A^*$, then 
$\Delta X\, \Delta Y = L_1L_2 \Delta t - \lambda^2L_1L_2\Delta \Lambda$. 
\end{theorem}

\begin{proof}
The quantum It\^o rule follows from writing out the definition 
of $\Delta(XY)$:
\begin{equation*}
\begin{split}
    \Delta(XY)(l) & = X(l)\,Y(l) - X(l-1)\,Y(l-1) \\ 
    & = (X(l-1)+ \Delta X(l))(Y(l-1) + \Delta Y(l)) - X(l-1)\,Y(l-1)\\
    & = X(l-1)\,\Delta Y(l) + \Delta X(l)\,Y(l-1) + \Delta X(l)\,\Delta Y(l).
\end{split}
\end{equation*}
The It\^o table follows directly from the definition of the increments
$\Delta\Lambda,\Delta A,\Delta A^*,\Delta t$ in section 
\ref{sec:discretefield} by explicit matrix multiplication.
\qquad\end{proof}

The following lemma shows that transforms of $A$, $A^*$ and $\Lambda$ 
have zero expectation. It is a fancy way of saying that 
$\phi(\sigma_-) = \phi(\sigma_+) = \phi(\sigma_+\sigma_-) = 0$ 
and $\phi(I) = 1$.

\begin{lemma}\label{lem vacuum expectation}
For all predictable processes $L$ we have 
$\mathbb{P}(L(l)\Delta t(l)) = \mathbb{P}(L(l))\,\lambda^2$ and
$
	  \mathbb{P}(L(l)\Delta A(l)) = 
	  \mathbb{P}(L(l)\Delta A^*(l)) = 
	  \mathbb{P}(L(l)\Delta \Lambda(l)) = 0
$
for any $1 \le l\le k$.
\end{lemma}

\begin{proof} 
Note that $\phi(\sigma_-) = \phi(\sigma_+) = \phi(\sigma_+\sigma_-) = 0$ 
and $\phi(I) = 1$.  Since $L$ is predictable, $L(i)$ and $\Delta A(i)$ are 
independent under $\mathbb{P}$ and we have for $0 \le i \le k$
\begin{equation*}
  \mathbb{P}(L(i)\Delta A(i)) = \lambda\,\mathbb{P}(L(i))\,
  \phi(\sigma_-) = 0.
\end{equation*} 
Similar reasoning holds for $A^*$, $\Lambda$ and $t$.
\qquad
\end{proof}

\subsection{Quantum stochastic difference equations}\label{sec:qsde}

We are now going to rewrite the repeated interactions of section 
\ref{sec:repeatedinteraction} as difference equations.  Clearly we have
$$
	\Delta U(l)=U(l)-U(l-1)=(M_l-I)\,U(l-1),
	\qquad U(0)=I.
$$
We will need to convert this expression into a more convenient form, 
however, in order for it to be of use in calculations.

Recall that $M_l$, and hence also $M_l-I$, is only a function of the 
initial system and of the $l$th time slice of the field; indeed, this is 
easily read off from Eq.\ (\ref{eq M}).  But any operator 
$X\in\mathscr{M}\otimes\mathscr{M}$ can be written in a unique way as
$$
	X = X^\pm\otimes\sigma_+\sigma_- + X^+\otimes\sigma_+
		+X^-\otimes\sigma_-+X^\circ\otimes I,\qquad
	X^{\circ,+,-,\pm}\in\mathscr{M},
$$
as $\sigma_+\sigma_-,\sigma_+,\sigma_-,I$ is a linearly independent basis 
for $\mathscr{M}$.  By the same reasoning we find that we can write 
without loss of generality
	\nomenclature{$M^{\circ,+,-,\pm}$}{Coefficients in difference
		equation for $U(l)$}
$$
	M_l-I=M^\pm\,\Delta\Lambda(l)+M^+\,\Delta A^*(l)+M^-\,\Delta A(l)
		+M^\circ\,\Delta t(l), 
$$ 
so that $M^{\circ,+,-,\pm}\in\mathscr{B}_0$ are uniquely determined by
$L_{1,2,3}$ and $\lambda$.  The functional dependence of $M^{\circ,+,-,\pm}$
on $L_{1,2,3}$ and $\lambda$ is complicated and we will not attempt to 
give a general expression; in the examples which we consider later it 
will be straightforward to calculate these matrices explicitly.

\begin{remark}
It should be emphasized that whereas $L_{1,2,3}\in\mathscr{B}_0$ do not
depend on the time step size $\lambda$, $M^{\circ,+,-,\pm}$ are in fact
functions of $\lambda$.  As usual, we have chosen the scaling of the
various parameters with $\lambda$ so that the limit as $\lambda\to 0$
gives meaningful results; in this case, it can be shown that the matrix
elements of $M^{\circ,+,-,\pm}$ converge to finite values as $\lambda\to 
0$, and the limiting matrices can be expressed explicitly in terms of 
$L_{1,2,3}$ (see \cite{Gough04} and other references in section 
\ref{continuous time}).
\end{remark}

We now obtain the following {\bf quantum stochastic difference equation}:
\begin{equation}\label{eq U}
	\Delta U(l)=\left\{
		M^\pm\,\Delta\Lambda(l)+M^+\,\Delta A^*(l)+M^-\,\Delta A(l)
		+M^\circ\,\Delta t(l)
	\right\}U(l-1),
\end{equation}
with $U(0)=I$.  It is this form of the equation that will be the most 
useful to us, as we can apply stochastic calculus techniques to manipulate 
it.  We could equivalently express it in terms of discrete quantum 
stochastic integrals:
$$
	U(l)=I+(M^\pm U_-)\cdot\Lambda(l)+
	(M^+U_-)\cdot A^*(l)+(M^-U_-)\cdot A(l)+(M^\circ U_-)\cdot t(l),
$$
where $U_-(l)=U(l-1)$.  This is the way that (quantum) stochastic 
differential equations are defined in continuous time.

To get some familiarity with calculations with the quantum It\^o rule and 
Lemma \ref{lem vacuum expectation}, let us calculate the time evolution of 
the expectation $\mathbb{P}(j_l(X))$ for an atomic observable 
$X \in \mathscr{B}_0$.  By linearity of the state we clearly have
$\Delta\mathbb{P}(j_l(X))=\mathbb{P}(\Delta j_l(X))$, so we are going to 
calculate $\Delta j_l(X)=\Delta(U(l)^*XU(l))$ using the quantum It\^o 
rule and then calculate the expectation of the resulting expression.
First, note that
$$
	\Delta U(l)^*=U(l-1)^*\left\{
		M^{\pm*}\,\Delta\Lambda(l)+M^{+*}\,\Delta A(l)
		+M^{-*}\,\Delta A^*(l)+M^{\circ*}\,\Delta t(l)
	\right\}.
$$
Then we calculate using the quantum It\^o rule
\begin{multline}
\label{eq:lindcalc}
	\Delta(U(l)^*XU(l))=
	(\cdots)\,\Delta\Lambda(l)+(\cdots)\,\Delta A^*(l)
		+(\cdots)\,\Delta A(l) \\
	+U(l-1)^*\left\{
		M^{+*}XM^+ +
		\lambda^2\,M^{\circ*}XM^\circ+M^{\circ*}X+XM^\circ
	\right\}U(l-1)\,\Delta t(l).
\end{multline}
We didn't bother to calculate the $\Delta\Lambda,\Delta A,\Delta A^*$ 
terms: by Lemma \ref{lem vacuum expectation} these vanish anyway when we 
take the expectation with respect to $\mathbb{P}$, so we only need the 
$\Delta t$ term.   Hence we obtain for any $X\in\mathscr{B}_0$
\begin{equation}\label{eq:mastereq}
	\frac{\Delta\,\mathbb{P}(j_l(X))}{\Delta t}=
	\mathbb{P}(j_{l-1}(\mathcal{L}(X))),
\end{equation}
where the \emph{discrete Lindblad generator} $\mathcal{L}(\cdot)$ is 
defined by
\begin{equation}\label{eq Lindblad}
	\mathcal{L}(X) =
		M^{+*}XM^+ +
		\lambda^2\,M^{\circ*}XM^\circ+M^{\circ*}X+XM^\circ.
\end{equation}
	\nomenclature{$\mathcal{L}(X)$}{Discrete Lindblad generator}

\begin{remark}
\label{rem:markov}
Let $X_t$ be a classical, continuous time Markov diffusion, and define the 
semigroup $P_t(f)=\mathbf{E}f(X_t)$.  Then (under sufficient regularity 
conditions) $dP_t(f)/dt=P_t(\mathcal{L}f)$, where $\mathcal{L}$ is the 
infinitesimal generator of the Markov diffusion.  Eq.\ (\ref{eq:mastereq}) 
is strongly reminiscent of this formula, and indeed the continuous time 
version of the formula plays an equivalent role in quantum probability.  
In fact, the semigroup property (expressed in ``differential'' form 
(\ref{eq:mastereq})) suggests that we should interpret our repeated 
interaction model as a ``quantum Markov process'', which can be given a 
precise meaning \cite{Kum85}, but we will not do so here.
\end{remark}

Note that we are not allowed to choose $M^{\pm,+,-,\circ}$ arbitrarily:  
these matrices must be chosen in such a way that the solution $U(l)$ of
the corresponding difference equation is unitary.  Obtaining
$M^{\pm,+,-,\circ}$ directly from $L_{1,2,3}$, as we will do in the
examples below, ensures that this is the case.  One could establish
general conditions on $M^{\pm,+,-,\circ}$ to ensure that $U(l)$ is
unitary, but we will not need these in the following.  There is one
particular necessary condition, however, that we will often use.

\begin{lemma}\label{lem:livanishes}
For any repeated interaction model $\mathcal{L}(I)=0$.
\end{lemma}

\begin{proof}
As $U(l)$ is unitary, $j_l(I)=U(l)^*U(l)=I$ for any $l$.  By Eq.\ 
(\ref{eq:mastereq}) we have
$$
	0=\frac{\mathbb{P}(j_1(I)-j_0(I))}{\Delta t}=
		\mathbb{P}(j_0(\mathcal{L}(I)))=
		\rho(\mathcal{L}(I)).
$$
But this must hold for any initial state $\rho$, so $\mathcal{L}(I)=0$ 
identically.
\qquad
\end{proof}

\subsection{Examples}\label{sec:examplesqsde}

Let us now return to the examples of section \ref{sec:examplesrepint}.
We will calculate the difference equations for the unitary evolution 
$U(l)$ explicitly.

\subsubsection*{Spontaneous emission}

Recall that in this case (with $2\kappa=1$)
$$
	M_l=\exp\left(
		\sigma_-\,\Delta A^*(l) - \sigma_+\,\Delta A(l)
	\right).
$$
As $M_l$ is of the form $P\otimes I^{\ten l-1}\otimes Q\otimes I^{\ten k-l}$,
we can read off the form of $M_l$ by calculating the matrix exponential of 
$$
	\sigma_-\otimes \lambda\,\sigma_+ - \sigma_+\otimes\lambda
		\,\sigma_-
	=
	\begin{pmatrix}
		0 & 0 & 0 & 0 \\
		0 & 0 & \lambda & 0 \\
		0 & -\lambda & 0 & 0 \\
		0 & 0 & 0 & 0
	\end{pmatrix}.
$$
Performing this calculation explicitly, we get
\begin{multline*}
	e^{\sigma_-\otimes \lambda\,\sigma_+ - 
		\sigma_+\otimes\lambda\,\sigma_-} =
	\begin{pmatrix}
		1 & 0 & 0 & 0 \\
		0 & \cos\lambda & \sin\lambda & 0 \\
		0 & -\sin\lambda & \cos\lambda & 0 \\
		0 & 0 & 0 & 1
	\end{pmatrix} \\
	=
	\begin{pmatrix}
		1 & 0 \\
		0 & \cos\lambda
	\end{pmatrix}\otimes \sigma_+\sigma_- +
	\sin\lambda\,\sigma_-\otimes \sigma_+ -
	\sin\lambda\,\sigma_+\otimes \sigma_- +
	\begin{pmatrix}
		\cos\lambda & 0 \\
		0 & 1
	\end{pmatrix}\otimes \sigma_-\sigma_+ \\
	=
	(1-\cos\lambda)\,\sigma_z
	\otimes \sigma_+\sigma_- +
	\sin\lambda\,\sigma_-\otimes \sigma_+ -
	\sin\lambda\,\sigma_+\otimes \sigma_- +
	\begin{pmatrix}
		\cos\lambda & 0 \\
		0 & 1
	\end{pmatrix}\otimes I,
\end{multline*}
where we have used $\sigma_-\sigma_+=I-\sigma_+\sigma_-$.
Hence we obtain
$$
	M_l=
	(1-\cos\lambda)\,\sigma_z\,\Delta\Lambda(l) +
	\frac{\sin\lambda}{\lambda}\,(\sigma_-\,\Delta A^*(l) -
	\sigma_+\,\Delta A(l)) +
	\frac{\cos\lambda-1}{\lambda^2}\,\sigma_+\sigma_-\,\Delta t(l)
	+I.
$$
We can now immediately read off the coefficients in the quantum stochastic 
difference equation for the spontaneous emission model:
\begin{equation}\label{eq:spontemcoeff}
\begin{split}
	& M^\pm=(1-\cos\lambda)\,\sigma_z,
	\qquad
	M^+=\frac{\sin\lambda}{\lambda}\,\sigma_-, \\
	& M^-=-\frac{\sin\lambda}{\lambda}\,\sigma_+,\qquad
	M^\circ = \frac{\cos\lambda-1}{\lambda^2}\,\sigma_+\sigma_-.
\end{split}
\end{equation}

\subsubsection*{Dispersive interaction}

Recall that in this case (with $g=1$)
$$
	M_l=\exp\left(
		\sigma_z(\Delta A^*(l) - \Delta A(l))
	\right).
$$
We proceed as in the case of spontaneous emission.  Starting from
$$
	\sigma_z\otimes\lambda\,(\sigma_+-\sigma_-)=
	\begin{pmatrix}
		0 & 0 & \lambda & 0 \\
		0 & 0 & 0 & -\lambda \\
		-\lambda & 0 & 0 & 0 \\
		0 & \lambda & 0 & 0
	\end{pmatrix},
$$
we calculate the matrix exponential
$$
	e^{\sigma_z\otimes\lambda\,(\sigma_+-\sigma_-)}=
	\begin{pmatrix}
		\cos\lambda & 0 & \sin\lambda & 0 \\
		0 & \cos\lambda & 0 & -\sin\lambda \\
		-\sin\lambda & 0 & \cos\lambda & 0 \\
		0 & \sin\lambda & 0 & \cos\lambda
	\end{pmatrix}=
	\cos\lambda\, I+\sin\lambda\,\sigma_z\otimes(\sigma_+-\sigma_-).
$$
Hence we obtain
$$
	M_l=\frac{\sin\lambda}{\lambda}\,\sigma_z\,(\Delta A^*(l)-\Delta 
		A(l)) + \frac{\cos\lambda-1}{\lambda^2}\,\Delta t(l)+I.
$$
We can now immediately read off the coefficients in the quantum stochastic 
difference equation for the dispersive interaction model:
\begin{equation}\label{eq:dispersivecoeff}
	M^\pm=0,\qquad
	M^+=-M^-=\frac{\sin\lambda}{\lambda}\,\sigma_z,\qquad
	M^\circ = \frac{\cos\lambda-1}{\lambda^2}\,I.
\end{equation}

\subsubsection*{A first simulation}

\begin{figure}
\centering
\includegraphics[width=\textwidth]{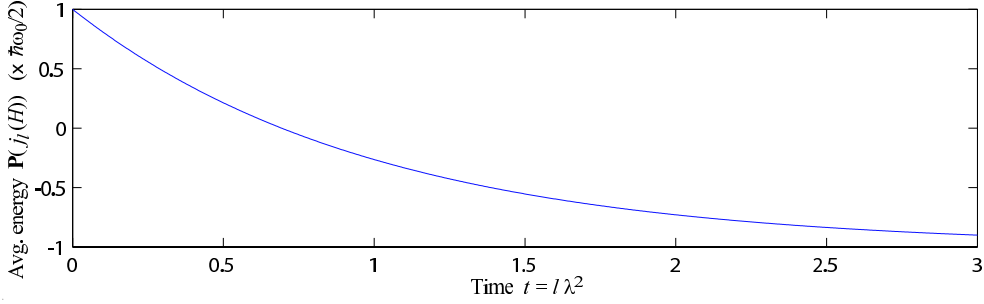}
\caption{\label{fig:lindblad} The expectation of the energy of a 
spontaneously emitting atom as a function of time.  The initial state 
$\rho$ was chosen in such a way that the atomic energy takes its maximal 
value $+\hbar\omega_0/2$ with unit probability at $t=0$, i.e.\ $\rho(x)=
{\rm Tr}[\sigma_+\sigma_- X]$.  The expected energy decays exponentially, 
until it reaches its minimum value $-\hbar\omega_0/2$.  The time scale 
used for the calculation is $\lambda^{-2}=300$.
} \end{figure}

To get some idea for the mean behavior of our two examples, let us 
calculate the expectation of the energy observable $\mathbb{P}(j_l(H))$ as 
a function of time.  The expectation of any atomic observable is given by 
Eq.\ (\ref{eq:mastereq}).  To simulate Eq.\ (\ref{eq:mastereq}) directly, 
we can use a standard technique which will be used several times in the 
following.  Note that $\mathbb{P}(j_l(X))$ is a linear function of $X$.  
Hence we can always find a $2\times 2$ matrix $\tau_l$ such that ${\rm 
Tr}[\tau_lX]=\mathbb{P}(j_l(X))$ for every $X$.  Substituting this 
expression into Eq.\ (\ref{eq:mastereq}), we obtain explicitly a recursive 
equation for $\tau_l$:
$$
	\frac{\Delta\tau_l}{\Delta t}=
		M^{+}\tau_{l-1}M^{+*} +
		\lambda^2\,M^{\circ}\tau_{l-1}M^{\circ*}
		+M^{\circ}\tau_{l-1}+\tau_{l-1}M^{\circ*},
$$
where $\tau_0$ is the density matrix corresponding to the initial state 
$\rho$:  $\rho(X)={\rm Tr}[\tau_0X]$ for every $X$.  This equation is 
called the (discrete) \emph{master equation}, and plays a similar role to 
the forward Kolmogorov (or Fokker-Planck) equation in classical diffusion 
theory.  The recursion is easily implemented in a computer program, and 
the solution allows us to calculate $\mathbb{P}(j_l(X))$ for any $X$.
For the spontaneous emission case, the expected energy $\mathbb{P}(j_l(H))=
{\rm Tr}[\tau_lH]$ is plotted in Figure \ref{fig:lindblad}.

We have used this opportunity as an excuse to introduce the master 
equation; in this simple case, however, we can obtain the result of Figure 
\ref{fig:lindblad} much more directly.  Let us calculate 
$\mathcal{L}(H)=(\hbar\omega_0/2)\,\mathcal{L}(\sigma_z)$ in the case of 
spontaneous emission:
$$
	\mathcal{L}(\sigma_z)=-\frac{2\sin^2\lambda}{\lambda^2}\,
	\sigma_+\sigma_- = -\frac{\sin^2\lambda}{\lambda^2}\,(\sigma_z+I).
$$
Hence we have
$$
	\mathbb{P}(j_l(\sigma_z))+1=
	\Delta\mathbb{P}(j_l(\sigma_z))+\mathbb{P}(j_{l-1}(\sigma_z))+1=
	\cos^2\lambda\,(\mathbb{P}(j_{l-1}(\sigma_z))+1).
$$
Recursing this relation, we find that the expected energy decays 
geometrically:
$$
	\mathbb{P}(j_l(H))=
	-\frac{\hbar\omega_0}{2}+\left(\rho(H)+\frac{\hbar\omega_0}{2}\right)
	\,(\cos^{2}\lambda)^l.
$$
This expression coincides exactly with the plot in Figure \ref{fig:lindblad}.

In the dispersive case, the time dependence of the mean energy is trivial: 
it is not difficult to check that in this case $\mathcal{L}(H)=0$.  This 
is not too surprising, as we have suppressed the atom's ability to radiate 
its excess energy by placing it in an optical cavity.  In fact, using 
$M_l=\exp(\sigma_z(\Delta A^*(l) - \Delta A(l)))$ we find that
$$
	j_l(H) = M_1^*\cdots M_{l-1}^*M_l^*HM_lM_{l-1}\cdots M_1=H.
$$
Evidently not only the expected energy, but even the energy as a random 
variable, is conserved in the dispersive case.

\section{The martingale method}\label{sec doob}

In this section we will provide a solution to the filtering problem using
martingale methods.  The key tools we need are minor variations on two
classical theorems.  First, there is a theorem of Doob which states that
any adapted process can be decomposed as the sum of a predictable process
and a martingale; we extend this theorem slightly to a class of
nondemolition processes.  Applying this theorem to the conditional
expectation $\pi_l(X)$, we can immediately identify the corresponding
predictable part.  Next, the martingale representation theorem tells us
that any martingale can be written as a stochastic integral of some
predictable process.  Consequently we can identify also the martingale
part of $\pi_l(X)$ in terms of a stochastic integral. The recursive nature
of the solution is the motivation behind this procedure: the predictable
parts will exactly turn out to be functions of the filter in the previous
time step, whereas the increment in the discrete stochastic integral will
be directly related to the observation increment.

For excellent introductions to the use of martingale methods in classical 
filtering theory we refer to \cite{DM81,Kr05}; definitive treatments can 
be found in \cite{LiS01,Kal80}.  Martingale methods for quantum filtering
were introduced in \cite{Bel92b} and used in \cite{BGM04,BHJ06}.  The 
treatment below is a discretized version of the latter.

\subsection{The classic theorems}

Recall that $\{\mathscr{Y}_l\}$ is the commutative filtration generated 
by the observations.

\begin{definition}
A quantum process $X$ is called \emph{nondemolished} (by the observations 
$Y$) if $X(l)$ is in the commutant of $\mathscr{Y}_l$ for $0 \le l \le k$.
A nondemolished quantum process $H$ is called an \emph{nd-martingale} 
(with respect to the filtration $\mathscr{Y}_{0\le l \le k}$) if
\begin{equation*} 
  	\mathbb{P}(H(l)|\mathscr{Y}_m) = \mathbb{P}(H(m)|\mathscr{Y}_m)
  	\qquad \forall\,0\le m \le l \le k.
\end{equation*}
A quantum process $A$ is called \emph{$\mathscr{Y}$-predictable} if 
$A(l)$ is in $\mathscr{Y}_{l-1}$ for $1 \le l \le k$, and $A$ is called 
\emph{$\mathscr{Y}$-adapted} if $A(l)$ is in $\mathscr{Y}_{l}$ for $1 \le 
l \le k$.  An nd-martingale that is additionally $\mathscr{Y}$-adapted is 
called a \emph{$\mathscr{Y}$-martingale} or simply a \emph{martingale}.
\end{definition}

The nondemolition requirement ensures that the conditional expectations in
the definition of an nd-martingale are well defined: $\mathscr{Y}_m\subset
\mathscr{Y}_l$ implies $\mathscr{Y}_m'\supset\mathscr{Y}_l'$, so 
$H(l)\in\mathscr{Y}_l'$ implies $H(l)\in\mathscr{Y}_m'$.  Note also that 
by construction, any $\mathscr{Y}$-predictable or $\mathscr{Y}$-adapted
process is a function of the observations and is hence commutative.

\begin{remark}
The concept of a $\mathscr{Y}$-martingale, being $\mathscr{Y}$-adapted and 
hence a classical process, coincides with the classical notion of a 
martingale.  An nd-martingale (``nd'' for nondemolition) is used in a
slightly broader sense; we will encounter an example below of a 
commutative nd-martingale that is not a martingale.  It will be convenient 
to use this terminology, as we will use the following theorem.  
\end{remark}

\begin{theorem}[\bf Doob decomposition]\label{thm doob}
Let $X$ be a quantum process nondemolished by the observations $Y$. 
Then $X$ has the following \emph{Doob decomposition}
\begin{equation*}
  	X(l) = X(0) + A(l) + H(l),\qquad l=1,\ldots,k
\end{equation*}
where $H$ is an nd-martingale (w.r.t.\ $\mathscr{Y}_{0 \le l \le 
k}$) null at $0$, and $A$ is a $\mathscr{Y}$-predictable process null at 
$0$.  Moreover, the decomposition is unique modulo indistinguishability. 
\end{theorem}

The proof is very similar to its classical counterpart \cite{Wil91}.

\begin{proof}
Suppose $X$ has a Doob decomposition as in the Theorem.  Then 
\begin{multline}\label{eq:deltapredictable}
  	\mathbb{P}(X(l)-X(l-1)|\mathscr{Y}_{l-1}) \\ =
   	\mathbb{P}(H(l)-H(l-1)|\mathscr{Y}_{l-1}) +
  	\mathbb{P}(A(l)-A(l-1)|\mathscr{Y}_{l-1}) =
  	A(l)-A(l-1),
\end{multline}
where we have used that $H$ is an nd-martingale and that $A$ is 
predictable.  Hence
\begin{equation}\label{eq previsible}
  	A(l) = 
	\sum_{m=1}^l \mathbb{P}(X(m) - X(m-1)|\mathscr{Y}_{m-1}),
	\qquad A(0)=0.
\end{equation}
For any nondemolished process $X$, define the predictable process $A(l)$ 
as in Eq.\ (\ref{eq previsible}), and define $H(l)=X(l)-X(0)-A(l)$.  Then it 
is easily verified that $H$ is an nd-martingale, hence we have explicitly 
constructed a Doob decomposition.

To prove uniqueness, suppose that $X(l)= X(0)+\tilde{A}(l)+\tilde{H}(l)$ 
is another Doob decomposition. Then $A(l)-\tilde{A}(l) = H(l) -\tilde{H}(l)$ 
for $1 \le l \le k$, and hence
\begin{multline}\label{eq:predictableversion}
	A(l)-\tilde{A}(l)=
	\mathbb{P}(A(l)-\tilde{A}(l)|\mathscr{Y}_{l-1}) \\ 
	= \mathbb{P}(H(l)-\tilde{H}(l)|\mathscr{Y}_{l-1})
	= \mathbb{P}(H(l-1)-\tilde{H}(l-1)|\mathscr{Y}_{l-1}) \\ 
	= \mathbb{P}(A(l-1)-\tilde{A}(l-1)|\mathscr{Y}_{l-1}) =
	A(l-1)-\tilde{A}(l-1)
\end{multline}
where we have used predictability of $A,\tilde{A}$ and the nd-martingale 
property of $H,\tilde H$.  But as $A(0)=\tilde{A}(0)=0$, we obtain by 
induction $0=A(l)-\tilde{A}(l)=H(l)-\tilde{H}(l)$ for $1 \le l \le k$.
Hence $A=\tilde{A}$ and $H=\tilde H$ with probability one (as the 
conditional expectations in Eq.\ (\ref{eq:predictableversion}) are only 
defined up to a choice of version).
\qquad
\end{proof}

\begin{remark}
It should be noted that the Doob decomposition depends crucially on the
choice of filtration $\mathscr{Y}$, which is demonstrated by the following
trivial example.  Consider the (commutative) filtration
$\{\mathscr{I}_l\}$ with $\mathscr{I}_l={\rm alg}\{I\}$ for any $l$.  As
$\mathscr{I}_l$ contains only multiples of the identity, the commutant is
the entire algebra $\mathscr{I}'=\mathscr{M}\otimes\mathscr{W}_k$.  Hence
any process $X$ is nondemolished by $\mathscr{I}$, and as
$\mathbb{P}(X(l)|\mathscr{I}_m)=\mathbb{P}(X(l))$ for any $l,m$ any
process with constant expectation $\mathbb{P}(X(l))=\mathbb{P}(X(m))$ is 
an nd-martingale with respect to $\{\mathscr{I}_l\}$ (but not 
necessarily an $\mathscr{I}$-martingale!)  Using $\mathscr{I}$ as the 
filtration,  we obtain the Doob decomposition $A(l)=\mathbb{P}(X(l)-X(0))$ 
and $H(l)=X(l)-X(0)-\mathbb{P}(X(l)-X(0))$ for any process $X$.  Clearly 
this decomposition is different than the Doob decomposition with respect 
to $\mathscr{Y}$; but note that $H(l)$ is not an nd-martingale with respect 
to $\{\mathscr{Y}_l\}$, so uniqueness is not violated.  The moral of the 
story is that we have to be careful to specify with respect to which 
filtration we decompose a process.  In the following, this will always be 
the filtration $\{\mathscr{Y}_l\}$ generated by the observations. 
\end{remark}

Let $X\in\mathscr{B}_0$.  Applying the Doob decomposition to $\pi_l(X)$ 
and $Y(l)$ gives
	\nomenclature{$\tilde Y(l)$}{Innovations process}
\begin{equation}\label{eq innovations}\begin{split}
  	\pi_l(X) = \rho(X) + B(l) + H(l),\ \ \ \ &\mbox{i.e.} \ \ \ \
  	\Delta \pi_l(X) = \Delta B(l) + \Delta H(l),\\
  	Y(l) = Y(0) + C(l) + \tilde{Y}(l),\ \ \ \  &\mbox{i.e.} \ \ \ \
  	\Delta Y(l) = \Delta C(l) + \Delta\tilde{Y}(l),  
\end{split}\end{equation}
where $B$ and $C$ are predictable processes null at $0$ and $\tilde{Y}$ 
and $H$ are $\mathscr{Y}$-martingales null at $0$. The process $\tilde{Y}$ 
is called the \emph{innovating martingale}. In the next two subsections we 
will investigate the processes $C$ and $\tilde{Y}$ in more detail for both 
the counting and homodyne detection cases. 

\begin{lemma}
The predictable process $B$ in the decomposition of $\pi_l(X)$ is given by 
\begin{equation*}
  	\Delta B(l) =  
	\pi_{l-1}(\mathcal{L}(X))\,\Delta t(l)\,
	\qquad 1 \le l \le k, 
\end{equation*}
where $\mathcal{L}$ is the discrete Lindblad generator
of Eq.\ \eqref{eq Lindblad}.
\end{lemma}

\begin{proof}
By Eq.\ (\ref{eq:deltapredictable}), we have $\Delta 
B(l)=\mathbb{P}(\Delta\pi_l(X)|\mathscr{Y}_{l-1})=
\mathbb{P}(\Delta j_l(X)|\mathscr{Y}_{l-1})$.  To calculate the latter,
let $K$ be an element in $\mathscr{Y}_{l-1}$.  Using Eq.\ 
(\ref{eq:lindcalc}), we obtain
\begin{equation*}
  \mathbb{P}(K\,\Delta j_l(X)) = 
  \mathbb{P}(K\,j_{l-1}(\mathcal{L}(X))\,\Delta t(l)) = 
  \mathbb{P}(K\,\pi_{l-1}(\mathcal{L}(X))\,\Delta t(l)).
\end{equation*} 
As this holds for any $K\in\mathscr{Y}_{l-1}$, and as 
$\pi_{l-1}(\mathcal{L}(X))\in\mathscr{Y}_{l-1}$, the statement of the 
Lemma follows from the definition of the conditional expectation.
\qquad
\end{proof}

At this point, we have the expression
$$
	\Delta\pi_l(X)=\pi_{l-1}(\mathcal{L}(X))\,\Delta t(l) + 
	\Delta H(l). $$ This is almost a recursive equation: what we would
like to do is write something like $\Delta
H(l)=f(\pi_{l-1}(\cdot))\,\Delta Y(l)$, as in that case we could use this
equation to calculate $\pi_l(\cdot)$ using only $\pi_{l-1}(\cdot)$ and
$\Delta Y(l)$.  The problem is that $f(\pi_{l-1}(\cdot))\,\Delta Y(l)$
does not define a martingale; but
$f(\pi_{l-1}(\cdot))\,\Delta\tilde{Y}(l)$ does!  This suggests that we
should try to represent $H$ as a discrete stochastic integral with respect
to the innovating martingale $\tilde{Y}(l)$.  The martingale
representation theorem shows that this is always possible.

\begin{theorem}[\bf Martingale representation]\label{thm martingale rep}
Let $\tilde{Y}$ be the innovating martingale and let $H$ be a 
$\mathscr{Y}$-martingale null at $0$.  Then there exists a 
$\mathscr{Y}$-predictable process $\Xi$ such that $\Delta H(l) = \Xi(l)\, 
\Delta\tilde{Y}(l)$, $l=1,\ldots,k$ modulo indistinguishability.
\end{theorem}

The following proof is reminiscent of \cite[pp.\ 154-155]{Wil91}, but the 
details of the argument are a little more delicate in our case.

\begin{proof}
As all the observables in the theorem are contained in the full 
observation algebra $\mathscr{Y}_k$, which is commutative, this is 
essentially a classical problem.  It will be convenient for the proof to 
treat it as such, i.e., applying the spectral theorem to 
$(\mathscr{Y}_k,\mathbb{P})$ gives the classical probability space 
$(\Omega,\mathcal{F},\mathbf{P})$, the filtration $\{\mathscr{Y}_l\}$ gives 
rise to the classical filtration $\{\mathcal{Y}_l\}$, and we will write 
$y_l=\iota(Y(l))$, $\tilde y_l=\iota(\tilde Y(l))$, and $h_l=\iota(H(l))$.
It will be convenient to write $\omega_l=\iota(\Delta Y(l))=y_l-y_{l-1}$.

We will make fundamental use of the following fact: $\omega_l$ takes one
of two values $\{\omega^+,\omega^-\}$ for every $l=1,\ldots,k$.  To see
this, recall that $\Delta Y(l)=U(l)^*\Delta Z(l)U(l)$ where $\Delta Z(l)$
is one of $\Delta A(l)+\Delta A^*(l)$ or $\Delta\Lambda(l)$ (in fact, any
observable $\Delta Z(l)$ of the form $(Z)_l\in\mathscr{W}_k$,
$Z\in\mathscr{M}$ would do, provided $Z$ is not a multiple of the identity.)
Hence $\Delta Z(l)$ has two distinct eigenvalues, and as unitary rotation 
leaves the spectrum of an operator invariant so does $\Delta Y(l)$.  It 
follows that $\iota(\Delta Y(l))$ is a two-state random variable.  We will 
write $p_l^\pm=\mathbf{P}(\omega_l=\omega^\pm|\mathcal{Y}_{l-1})$ for the 
conditional probability that $\omega_l=\omega^\pm$ given that we have 
observed $\omega_1,\ldots,\omega_{l-1}$.

Now recall that $\mathscr{Y}_l$ is the algebra generated by $\Delta Y(i)$, 
$i=1,\ldots,l$.  Hence every $\mathcal{Y}_l$-measurable random variable 
can be written as a function of $\omega_i$, $i=1,\ldots,l$.  In particular,
$\tilde y_l=\tilde y_l(\omega_1,\ldots,\omega_l)$ and 
$h_l=h_l(\omega_1,\ldots,\omega_l)$.  We would like to find a 
$\xi_l(\omega_1,\ldots,\omega_{l-1})$ (independent of $\omega_l$, hence 
predictable) that satisfies
\begin{equation}\label{eq:mtgtoprove}
	\Delta h_l(\omega_1,\ldots,\omega_l) =
	\xi_l(\omega_1,\ldots,\omega_{l-1})\,
	\Delta \tilde y_l(\omega_1,\ldots,\omega_l).
\end{equation}
To proceed, we split $\Omega$ into three disjoint subsets $\Omega_{1,2,3}$
and define the random variable $\xi_l$ separately on each set.

{\bf Case 1:} $\Omega_1=\{\omega\in\Omega:\Delta\tilde 
y_l(\omega_1,\ldots,\omega_{l-1},\omega^\pm)\ne 0\}$.  
Let us suppose that $\xi_l$ exists.  Then evidently on $\Omega_1$
$$
	\xi_l(\omega_1,\ldots,\omega_{l-1})=
	\frac{\Delta h_l(\omega_1,\ldots,\omega_l)}{
		\Delta \tilde y_l(\omega_1,\ldots,\omega_l)}.
$$
Existence of $\xi_l$ is thus verified by construction if we can show that 
the right hand side is independent of $\omega_l$.  To this end, we express 
the martingale property of $h_l$ as
$$
	\Delta h_l(\omega_1,\ldots,\omega_{l-1},\omega^+)\,
	p_l^+(\omega_1,\ldots,\omega_{l-1})+
	\Delta h_l(\omega_1,\ldots,\omega_{l-1},\omega^-)\,
	p_l^-(\omega_1,\ldots,\omega_{l-1})=0
$$
a.s., where the left hand side is simply the expression for 
$\mathbf{E_P}(\Delta h_l|\mathcal{Y}_{l-1})$.  Similarly, 
$$
	\Delta\tilde y_l(\omega_1,\ldots,\omega_{l-1},\omega^+)\,
	p_l^+(\omega_1,\ldots,\omega_{l-1})+
	\Delta\tilde y_l(\omega_1,\ldots,\omega_{l-1},\omega^-)\,
	p_l^-(\omega_1,\ldots,\omega_{l-1})=0
$$
a.s., as $\tilde y_l$ is a martingale.  Note that necessarily $p_l^\pm\ne 
0$ a.s.\ on $\Omega_1$.  Hence we obtain
\begin{eqnarray*}
	\Delta h_l(\omega_1,\ldots,\omega_{l-1},\omega^+)=-
	\Delta h_l(\omega_1,\ldots,\omega_{l-1},\omega^-)\,
	\frac{p_l^-(\omega_1,\ldots,\omega_{l-1})}{
		p_l^+(\omega_1,\ldots,\omega_{l-1})}, \\
	\Delta \tilde y_l(\omega_1,\ldots,\omega_{l-1},\omega^+)=-
	\Delta \tilde y_l(\omega_1,\ldots,\omega_{l-1},\omega^-)\,
	\frac{p_l^-(\omega_1,\ldots,\omega_{l-1})}{
		p_l^+(\omega_1,\ldots,\omega_{l-1})}.
\end{eqnarray*}
Dividing the first by the second expression, the independence of $\xi_l$ 
from $\omega_l$ follows.

{\bf Case 2:} $\Omega_2=\{\omega\in\Omega:p_l^+\in\{0,1\}\}$.  Using the 
martingale property of $h_l$ and $\tilde y_l$ as above, we conclude that 
on $\Omega_2$ we have $\Delta h_l=\Delta \tilde y_l=0$ a.s.  Hence Eq.\ 
(\ref{eq:mtgtoprove}) holds regardless of the value we assign to $\xi_l$.

{\bf Case 3:} $\Omega_3=\Omega\backslash(\Omega_1\cup\Omega_2)$.  We will 
show that $\mathbf{P}(\Omega_3)=0$, so that we do not need to worry about 
defining $\xi_l$ on this set.  As on $\Omega_3$ we have $p_l^\pm\ne 0$ but 
one of $\Delta\tilde y_l(\omega_1,\ldots,\omega_{l-1},\omega^\pm)=0$, 
using the martingale property as above allows us to conclude that
$\Delta\tilde y_l=0$ a.s.\ on $\Omega_3$.  Recall that $\omega_l=\Delta 
y_l=\Delta c_l+\Delta\tilde y_l$ where $\Delta c_l$ is predictable.  Then
$\omega_l=\Delta c_l(\omega_1,\cdots,\omega_{l-1})$ a.s.\ on $\Omega_3$.
But this would imply that $\omega_l=\mathbf{E_P}(\omega_l|\mathcal{Y}_{l-1})=
\omega^+p^+_l+\omega^-p^-_l$, and as $p^\pm_l\ne 0$ we would be able to 
conclude that $\omega_l\ne \omega^\pm$.  Hence we have a contradiction.

We have now shown how to define $\xi_l$ that satisfies Eq.\
(\ref{eq:mtgtoprove}) except possibly on a set of measure zero.
Setting $\Xi(l)=\iota^{-1}(\xi_l)$, the theorem is proved. 
\qquad
\end{proof}

Though the proof of the discrete martingale representation theorem is in
principle constructive, it is not advisable to follow this complicated
procedure in order to calculate the predictable process $\Xi$.  Instead we
will calculate $\Xi$ using a standard trick of filtering theory, and it
will turn out to depend only on the conditional expectations in the
previous time step. Putting everything together, we obtain a recursive
relation with which we can update our conditional expectations of atomic
operators given the conditional expectations at the previous time step and
the observation result at the present time step. This recursion is called
the \emph{discrete quantum filtering equation}. As the predictable
processes $C$ and $\Xi$ depend on the nature of the observations, we
consider separately the homodyne and photon counting cases.

\subsection{Homodyne detection}\label{sec filter homo}

Let us first consider a homodyne detection setup, i.e.\ an experimental 
setup that allows us to observe 
\begin{equation*}
  Y^X(l) = U(l)^*(A(l) + A^*(l))U(l),\ \ \ \ 0 \le l\le k.
\end{equation*}
We begin by finding the predictable process $C$ in the Doob decomposition 
of $Y^X$.

\begin{lemma}\label{lem pred Y homo} The predictable process $C$ in the 
decomposition of $Y^X$ is given by
\begin{equation*}
	\Delta C(l) = \pi_{l-1}(M^++M^{+*}+\lambda^2
		M^{\circ*}M^++\lambda^2M^{+*}M^\circ)\,
	\Delta t(l),\qquad
	1 \le l \le k.
\end{equation*}
\end{lemma}

\begin{proof}
By Eq.\ (\ref{eq:deltapredictable}), we have $\Delta C(l)=
\mathbb{P}(\Delta Y^X(l)|\mathscr{Y}_{l-1})$.  To calculate the latter, 
let $K$ be an element in $\mathscr{Y}_{l-1}$; we would like to find an 
expression for $\mathbb{P}(K\Delta Y^X(l))$.  To this end, we calculate 
using the discrete quantum It\^o rules
\begin{multline}
\label{eq:homopreditonasty}
	\Delta((A(l)+A^*(l))U(l))=(\cdots)\,\Delta\Lambda(l)+
		(\cdots)\,\Delta A(l) \\
	+\left\{I+(A(l-1)+A^*(l-1))\,M^+ + \lambda^2M^\circ\right\}
		U(l-1)\,\Delta A^*(l) \\
	+ \left\{(A(l-1)+A^*(l-1))\,M^\circ + M^+\right\}U(l-1)\,\Delta t(l),
\end{multline}
where we have retained the relevant terms.  Consequently, we calculate
\begin{multline*}
	\Delta(U(l)^*(A(l)+A^*(l))U(l))=(\cdots)\,\Delta\Lambda(l)+
		(\cdots)\,\Delta A^*(l)+(\cdots)\,\Delta A(l) \\
	+U(l-1)^*\mathcal{L}(I)(A(l-1)+A^*(l-1))U(l-1)\,\Delta t(l) \\
	+j_{l-1}(M^++M^{+*}+\lambda^2
		M^{\circ*}M^++\lambda^2M^{+*}M^\circ)\,\Delta t(l),
\end{multline*}
where $\mathcal{L}(X)$ is the discrete Lindblad generator.  But
$\mathcal{L}(I)$ vanishes by Lemma \ref{lem:livanishes}.  Hence we obtain
$$
	\mathbb{P}(K\Delta Y^X(l))=
	\mathbb{P}(K\,j_{l-1}(M^++M^{+*}+\lambda^2
		M^{\circ*}M^++\lambda^2M^{+*}M^\circ)\,\Delta t(l))
$$
using Lemma \ref{lem vacuum expectation}, or equivalently (using 
$K\in\mathscr{Y}_{l-1}$)
$$
	\mathbb{P}(K\Delta Y^X(l))=
	\mathbb{P}(K\,\pi_{l-1}(M^++M^{+*}+\lambda^2
		M^{\circ*}M^++\lambda^2M^{+*}M^\circ)\,\Delta t(l)).
$$
As this holds for any $K\in\mathscr{Y}_{l-1}$, and as
$\pi_{l-1}(\cdots)\in\mathscr{Y}_{l-1}$, the statement of the Lemma 
follows from the definition of the conditional expectation.
\qquad
\end{proof}

From Theorem \ref{thm martingale rep} we know that $\Delta H = \Xi \Delta
\tilde{Y}$ for some predictable process $\Xi$.  It remains to determine
$\Xi$; we approach this problem using a standard technique.  Since
$\mathbb{P}(j_l(X)Y(l)|\mathscr{Y}_l)  = \pi_l(X)Y(l)$, the uniqueness of
the Doob decomposition ensures that $j_l(X)Y(l)$ and $\pi_l(X)Y(l)$ have
equal predictable parts.  We write $D$ and $E$ for the predictable
processes in the Doob decomposition of $j_l(X)Y(l)$ and $\pi_l(X)Y(l)$,
respectively. Solving the equation $\Delta D = \Delta E$ will then allow
us to determine $\Xi$.

\begin{lemma}\label{lem xi homo}
For any $X\in\mathscr{B}_0$, define
	\nomenclature{$\mathcal{J}(X)$}{$\mathcal{J}$-coefficient of the 
		filtering equations}
$$
	\mathcal{J}(X)=XM^++M^{+*}X+\lambda^2
                M^{\circ*}XM^++\lambda^2M^{+*}XM^\circ,
$$
so that $\Delta C(l)=\pi_{l-1}(\mathcal{J}(I))\,\Delta t(l)$.
Then we can write
$$
	\Xi(l)=
	\frac{\pi_{l-1}(\mathcal{J}(X))
	-\pi_{l-1}(X+\lambda^2\mathcal{L}(X))\,\pi_{l-1}(\mathcal{J}(I))}
	{I-\lambda^2\pi_{l-1}(\mathcal{J}(I))^2}.
$$
\end{lemma}

\begin{proof}
We begin by evaluating $\Delta D$.  For any $K\in\mathscr{Y}_{l-1}$ we
want to calculate $\mathbb{P}(\Delta(j_l(X)Y(l))K)
=\mathbb{P}(\Delta(U(l)^*X(A(l)+A^*(l))U(l))K)$.  We proceed
exactly as in the proof of Lemma \ref{lem pred Y homo}.  Using Eq.\
(\ref{eq:homopreditonasty}) and the quantum It\^o rules, we obtain
\begin{multline*}
	\Delta(U(l)^*X(A(l)+A^*(l))U(l))=(\cdots)\,\Delta\Lambda(l)+
		(\cdots)\,\Delta A^*(l)+(\cdots)\,\Delta A(l) \\
	+U(l-1)^*\mathcal{L}(X)(A(l-1)+A^*(l-1))U(l-1)\,\Delta t(l) \\
	+j_{l-1}(XM^++M^{+*}X+\lambda^2
		M^{\circ*}XM^++\lambda^2M^{+*}XM^\circ)\,\Delta t(l).
\end{multline*}
It follows in the usual way that
\begin{multline*}
	\Delta D(l)=
	\pi_{l-1}(\mathcal{L}(X))\,Y(l-1)\,\Delta t(l) \\ +
	\pi_{l-1}(XM^++M^{+*}X+\lambda^2
		M^{\circ*}XM^++\lambda^2M^{+*}XM^\circ)
	\,\Delta t(l).
\end{multline*}
We now turn our attention to $\Delta E$. First note that the It\^o rule 
gives
\begin{equation*}
  \Delta(\pi_{l}(X)Y(l)) = (\Delta\pi_{l}(X))\,Y(l-1) +
  \pi_{l-1}(X)\,\Delta Y(l) +  \Delta\pi_{l}(X)\, \Delta Y(l). 
\end{equation*}
By uniqueness of the Doob decomposition, $\Delta E$ is the sum of the 
predictable parts $\Delta E_{1,2,3}$ of the three terms on the right hand 
side.  Let us investigate each of these individually.  The first term can 
be written as
$$
	(\Delta\pi_{l}(X))\,Y(l-1) =
	\pi_{l-1}(\mathcal{L}(X))\,Y(l-1)\,\Delta t(l)+
	Y(l-1)\,\Delta H(l).
$$
It is easily verified, however, that $\sum_lY(l-1)\,\Delta H(l)$ inherits 
the martingale property from $H$.  Hence by uniqueness of the Doob 
decomposition, we obtain $\Delta E_1(l)=\pi_{l-1}(\mathcal{L}(X))\,
Y(l-1)\,\Delta t(l)$.  Moving on to the second term, we write
$$
	\pi_{l-1}(X)\,\Delta Y(l)=
	\pi_{l-1}(X)\,\Delta C(l)+\pi_{l-1}(X)\,\Delta\tilde Y(l).
$$
Similarly as above we find that $\sum_l\pi_{l-1}(X)\,\Delta\tilde Y(l)$ 
inherits the martingale property from $\tilde Y$, so that evidently
$\Delta E_2(l)=\pi_{l-1}(X)\,\Delta C(l)$ (where $\Delta C(l)$ is given 
explicitly in Lemma \ref{lem pred Y homo}).  It remains to deal with the 
third term.  To this end, let us write
$$
	\Delta\pi_{l}(X)\,\Delta Y(l) =
	\left\{\pi_{l-1}(\mathcal{L}(X))\,\Delta t(l)+
	\Xi(l)\,\Delta\tilde Y(l)\right\}
	(\Delta C(l)+\Delta\tilde Y(l)).
$$
As before, processes of the form $\sum_lX(l)\,\Delta\tilde Y(l)$ with 
$\mathscr{Y}$-predictable $X$ inherit the martingale property from $\tilde 
Y$.  Thus we only need to retain the predictable terms:
$$
	\Delta\pi_{l}(X)\,\Delta Y(l) =
	(\cdots)\,\Delta\tilde Y(l)+
	\pi_{l-1}(\mathcal{L}(X))\,\Delta t(l)\,\Delta C(l)+
	\Xi(l)\,(\Delta\tilde Y(l))^2.
$$
Similarly, we expand $(\Delta\tilde Y(l))^2$ as
$$
	(\Delta\tilde Y(l))^2=
	(\cdots)\,\Delta\tilde Y(l)
		+(\Delta Y(l))^2-(\Delta C(l))^2,
$$
where we have used $\Delta Y(l)=\Delta C(l)+\Delta\tilde Y(l)$.
But using the It\^o rules we calculate $(\Delta Y(l))^2=
U(l)^*(\Delta A(l)+\Delta A^*(l))^2U(l)=\Delta t(l)$.
Hence we can read off
$$
	\Delta E_3(l)=
	\pi_{l-1}(\mathcal{L}(X))\,\Delta t(l)\,\Delta C(l)+
	\Xi(l)\,(\Delta t(l)-(\Delta C(l))^2).
$$
But recall that $j_l(X)Y(l)$ and $\pi_l(X)Y(l)$ have equal predictable 
parts.  Hence setting $\Delta D(l)=\Delta E_1(l)+\Delta E_2(l)+\Delta 
E_3(l)$ and solving for $\Xi(l)$, the Lemma follows.
\qquad
\end{proof}

Putting everything together we get the following {\bf discrete quantum 
filtering equation for homodyne detection}
\begin{multline*}
	\Delta\pi_l(X)=
	\pi_{l-1}(\mathcal{L}(X))\,\Delta t(l) \\ +
	\frac{\pi_{l-1}(\mathcal{J}(X))
	-\pi_{l-1}(X+\lambda^2\mathcal{L}(X))\,\pi_{l-1}(\mathcal{J}(I))}
	{I-\lambda^2\pi_{l-1}(\mathcal{J}(I))^2}\,
	(\Delta Y(l)-\pi_{l-1}(\mathcal{J}(I))\,\Delta t(l)).
\end{multline*}

\subsection{Photodetection}\label{sec filter count}

We now turn our attention to a setup where we are counting photons in the 
field, i.e.\ we are observing 
\begin{equation*}
  	Y^\Lambda(l) = U(l)^*\Lambda(l)U(l),\qquad 0 \le l\le k.
\end{equation*}
The procedure here is much the same as in the homodyne detection case.

\begin{lemma}\label{lem pred count} 
The predictable process $C$ in the decomposition of $Y^\Lambda$ is given by
\begin{equation*}
	\Delta C(l) = \pi_{l-1}(M^{+*}M^+)\,\Delta t(l),\qquad
	1\le l\le k.
\end{equation*}
\end{lemma}

\begin{proof}
By Eq.\ (\ref{eq:deltapredictable}), we have $\Delta C(l)=
\mathbb{P}(\Delta Y^\Lambda(l)|\mathscr{Y}_{l-1})$.  To calculate the 
latter, let $K$ be an element in $\mathscr{Y}_{l-1}$; we would like find 
an expression for $\mathbb{P}(K\Delta Y^\Lambda(l))$.  To this end, we 
calculate using the discrete quantum It\^o rules
\begin{multline}
\label{eq:countpreditonasty}
	\Delta(\Lambda(l)U(l))=(\cdots)\,\Delta\Lambda(l)+
		(\cdots)\,\Delta A(l) \\
	+(I+\Lambda(l-1))M^+U(l-1)\,\Delta A^*(l)
	+ \Lambda(l-1)M^\circ U(l-1)\,\Delta t(l),
\end{multline}
where we have retained the relevant terms.  Consequently, we calculate
\begin{multline*}
	\Delta(U(l)^*\Lambda(l)U(l))=(\cdots)\,\Delta\Lambda(l)+
		(\cdots)\,\Delta A^*(l)+(\cdots)\,\Delta A(l) \\
	+U(l-1)^*\mathcal{L}(I)\Lambda(l-1)U(l-1)\,\Delta t(l)
	+j_{l-1}(M^{+*}M^+)\,\Delta t(l),
\end{multline*}
where $\mathcal{L}(X)$ is the discrete Lindblad generator.  But 
$\mathcal{L}(I)$ vanishes, and the Lemma follows by the usual argument.
\qquad
\end{proof}

Next, we determine the predictable process $\Xi(l)$ such that
$\Delta H(l)=\Xi(l)\,\Delta\tilde Y(l)$.

\begin{lemma}
The process $\Xi(l)$ is given by
$$
	\Xi(l)=
	(I-\lambda^2\pi_{l-1}(M^{+*}M^+))^{-1}
	\left[
	\frac{\pi_{l-1}(M^{+*}XM^+)}{\pi_{l-1}(M^{+*}M^+)}
	- \pi_{l-1}(X+\lambda^2\mathcal{L}(X))\right].
$$
\end{lemma}

\begin{proof}
We begin by finding the predictable process $\Delta D$ in the Doob 
decomposition of $j_l(X)Y(l)$.  For any $K\in\mathscr{Y}_{l-1}$ we
want to calculate $\mathbb{P}(\Delta(j_l(X)Y(l))K)
=\mathbb{P}(\Delta(U(l)^*X\Lambda(l)U(l))K)$.  Using Eq.\
(\ref{eq:countpreditonasty}) and the quantum It\^o rules, we obtain
\begin{multline*}
	\Delta(U(l)^*X\Lambda(l)U(l))=(\cdots)\,\Delta\Lambda(l)+
		(\cdots)\,\Delta A^*(l)+(\cdots)\,\Delta A(l) \\
	+U(l-1)^*\mathcal{L}(X)\Lambda(l-1)U(l-1)\,\Delta t(l)
	+j_{l-1}(M^{+*}XM^+)\,\Delta t(l).
\end{multline*}
It follows in the usual way that
$$
	\Delta D(l)=
	\pi_{l-1}(\mathcal{L}(X))\,Y(l-1)\,\Delta t(l) +
	\pi_{l-1}(M^{+*}XM^+)\,\Delta t(l).
$$
We now turn our attention to the predictable process $\Delta E$ in the 
Doob decomposition of $\pi_l(X)Y_l$.  First note that the It\^o rule 
gives
\begin{equation*}
  \Delta(\pi_{l}(X)Y(l)) = (\Delta\pi_{l}(X))\,Y(l-1) +
  \pi_{l-1}(X)\,\Delta Y(l) +  \Delta\pi_{l}(X)\, \Delta Y(l). 
\end{equation*}
By uniqueness of the Doob decomposition, $\Delta E$ is the sum of the 
predictable parts $\Delta E_{1,2,3}$ of the three terms on the right hand 
side.  As in the proof of Lemma \ref{lem xi homo}, we find that
$\Delta E_1(l)=\pi_{l-1}(\mathcal{L}(X))\,Y(l-1)\,\Delta t(l)$
and $\Delta E_2(l)=\pi_{l-1}(X)\,\Delta C(l)$ (where $\Delta C(l)$ is 
given explicitly in Lemma \ref{lem pred count}).  To deal with the third 
term, we write
$$
	\Delta\pi_{l}(X)\,\Delta Y(l) =
	\left\{\pi_{l-1}(\mathcal{L}(X))\,\Delta t(l)+
	\Xi(l)\,\Delta\tilde Y(l)\right\}
	(\Delta C(l)+\Delta\tilde Y(l)).
$$
As before, processes of the form $\sum_lX(l)\,\Delta\tilde Y(l)$ with 
$\mathscr{Y}$-predictable $X$ inherit the martingale property from $\tilde 
Y$.  Thus we only need to retain the predictable terms:
$$
	\Delta\pi_{l}(X)\,\Delta Y(l) =
	(\cdots)\,\Delta\tilde Y(l)+
	\pi_{l-1}(\mathcal{L}(X))\,\Delta t(l)\,\Delta C(l)+
	\Xi(l)\,(\Delta\tilde Y(l))^2.
$$
Similarly, we expand $(\Delta\tilde Y(l))^2$ as
$$
	(\Delta\tilde Y(l))^2=
	(\cdots)\,\Delta\tilde Y(l)
		+(\Delta Y(l))^2-(\Delta C(l))^2,
$$
where we have used $\Delta Y(l)=\Delta C(l)+\Delta\tilde Y(l)$.
But using the It\^o rules we calculate $(\Delta Y(l))^2=
U(l)^*\Delta\Lambda(l)^2U(l)=\Delta Y(l)$.
Hence
$$
	(\Delta\tilde Y(l))^2=
	(\cdots)\,\Delta\tilde Y(l)
		+\Delta C(l)-(\Delta C(l))^2,
$$
and we can read off
$$
	\Delta E_3(l)=
	\pi_{l-1}(\mathcal{L}(X))\,\Delta t(l)\,\Delta C(l)+
	\Xi(l)\,(I-\Delta C(l))\Delta C(l).
$$
But recall that $j_l(X)Y(l)$ and $\pi_l(X)Y(l)$ have equal predictable 
parts.  Hence setting $\Delta D(l)=\Delta E_1(l)+\Delta E_2(l)+\Delta 
E_3(l)$ and solving for $\Xi(l)$, the Lemma follows.
\qquad
\end{proof}

Putting everything together we obtain the following 
{\bf discrete quantum filtering equation for photon counting}
\begin{multline*}
	\Delta\pi_l(X)=\pi_{l-1}(\mathcal{L}(X))\,\Delta t(l)+
	(I-\lambda^2\pi_{l-1}(M^{+*}M^+))^{-1}\times\\
	\left[
	\frac{\pi_{l-1}(M^{+*}XM^+)}{\pi_{l-1}(M^{+*}M^+)}
	- \pi_{l-1}(X+\lambda^2\mathcal{L}(X))\right]
	(\Delta Y(l)-\pi_{l-1}(M^{+*}M^+)\,\Delta t(l)).
\end{multline*}

\subsection{How to use the filtering equations}\label{sec:howtouse}

The filtering equations of sections \ref{sec filter homo} and \ref{sec 
filter count} may seem a little abstract at this point; $\pi_l(X)$ is some 
observable in the algebra $\mathscr{M}\otimes\mathscr{W}_k$, and it 
appears that we would need to know $\pi_{l-1}(Z)$ for \emph{every}
$Z\in\mathscr{B}_0$, in addition to the observation increment $\Delta 
Y(l)$, in order to be able to calculate $\pi_l(X)$ for arbitrary $X$.  The 
equations are much less abstract than they might seem, however.
First of all, recall that both $\pi_l(X)$ and $\Delta Y(l)$ are elements 
of the (commutative) observation algebra $\mathscr{Y}_k$; in fact, the 
filtering equations live entirely within this algebra.  Hence these are 
just classical equations in disguise (as they should be!); we could write 
explicitly, e.g.\ in the homodyne detection case,
\begin{multline}\label{eq homo sillyclassical}
	\Delta\iota(\pi_l(X))=
	\iota(\pi_{l-1}(\mathcal{L}(X)))\,\lambda^2 
	+(\iota(\Delta Y(l))-\iota(\pi_{l-1}(\mathcal{J}(I)))\,\lambda^2)
	\times \\
	\frac{\iota(\pi_{l-1}(\mathcal{J}(X)))
	-\iota(\pi_{l-1}(X+\lambda^2\mathcal{L}(X)))\,
	\iota(\pi_{l-1}(\mathcal{J}(I)))}
	{1-\lambda^2\,\iota(\pi_{l-1}(\mathcal{J}(I)))^2}
\end{multline}
using the $^*$-isomorphism $\iota$ obtained by applying the spectral 
theorem to $(\mathscr{Y}_k,\mathbb{P})$.  For any $0\le l\le k$ and 
$X\in\mathscr{M}$, $\iota(\pi_l(X))$ is a random variable that is a 
function of the random process $\iota(\Delta Y(i))$ up to and including 
time $l$.  But an elementary property of the conditional expectation is 
that $\pi_l(X)$ is linear in $X$; hence $\iota(\pi_l(X))$ is also linear 
in $X$.  This means we can always write $\iota(\pi_l(X))={\rm 
Tr}[\rho_lX]$, where $\rho_l$ is a (random) $2\times 2$-matrix (as 
$\mathscr{M}$ is two-dimensional).  We obtain the following recursion for 
$\rho_l$:
\begin{multline}\label{eq homo density form}
	\Delta\rho_l=
	\mathcal{\overline{L}}(\rho_{l-1})\,\lambda^2 +
	\frac{\mathcal{\overline{J}}(\rho_{l-1})
	-{\rm Tr}[\mathcal{\overline{J}}(\rho_{l-1})]\,
	(\rho_{l-1}+\lambda^2\mathcal{\overline{L}}(\rho_{l-1}))}
	{1-\lambda^2\,{\rm Tr}[\mathcal{\overline{J}}(\rho_{l-1})]^2}\times\\
	(\iota(\Delta Y(l))-{\rm Tr}[\mathcal{\overline{J}}(\rho_{l-1})]
		\,\lambda^2),
\end{multline}
where we have written
\begin{equation*}
\begin{split}
	& \mathcal{\overline{L}}(\rho)=
	M^+\rho M^{+*}+\lambda^2\,M^\circ\rho M^{\circ *}+
	M^\circ \rho+\rho M^{\circ *}, \\
	& \mathcal{\overline{J}}(\rho)=
	M^+\rho+\rho M^{+*}+\lambda^2M^+\rho M^{\circ *}+
	\lambda^2M^\circ\rho M^{+*}.
\end{split}
\end{equation*}
$\rho_0$ is simply the density matrix corresponding to the initial state 
$\rho(\cdot)$, i.e.\ $\rho(X)={\rm Tr}[\rho_0X]$ for every 
$X\in\mathscr{M}$.  The matrix $\rho_l$ is called the \emph{conditional 
density matrix} and contains all the information needed to calculate 
$\pi_l(X)$ for every $X\in\mathscr{M}$.  Furthermore, Eq.\ (\ref{eq homo 
density form}) is a simple nonlinear recursion for $2\times 2$-matrices.  
At any time step $l$ we only need to remember the $2\times 2$-matrix 
$\rho_l$; when the $(l+1)$th observation $\iota(\Delta Y(l+1))$ becomes 
available, which takes one of the values $\pm\lambda$, we simply plug this 
value into Eq.\ (\ref{eq homo density form}) and obtain the updated matrix 
$\rho_{l+1}$.  Such a recursion is very efficient and would, if 
necessary, be easily implemented on a digital signal processor.  

As filtering equations are entirely classical, there is no real need to
make the explicit distinction between their representation in terms of
classical random variables on a probability space \emph{vs.}\ elements of
the observation algebra.  We urge the reader to always think of sets of
commuting observables as random variables: this is implied by the spectral
theorem, and is at the heart of quantum mechanics!  Eq.\ (\ref{eq homo
sillyclassical}) is notationally tedious and completely unnecessary, as
it does not add anything to the filtering equation as we have already
written it in section \ref{sec filter homo}.  In some cases, e.g.\ in the 
proof of the martingale representation theorem, it is convenient to use 
explicitly the structure of the underlying probability space; but in much 
of this article we will not make an explicit distinction between random 
variables and observables.

A similar story holds for the photodetection case; we leave it up to the 
reader to calculate the associated recursion for the conditional density 
matrix.

\subsection{The Markov property and Monte Carlo simulations}
\label{sec:filter-markov}

The filtering equations that we have developed take as input the
observation process obtained from the system.  Though this is precisely
how it should be, one would think that further investigation of the
filters can not proceed without the availability of typical sample
paths of the observations from some other source, be it an actual physical
system or a direct computer simulation of the underlying repeated
interaction model.  It is thus somewhat surprising that we can actually
simulate such sample paths using the filtering equation only, without any
auxiliary input.  This is due to the \emph{Markov property} of the filter,
which we will demonstrate shortly.  We can use this property of the
filtering equations to perform Monte Carlo simulations of both the sample
paths of the observation process and sample paths of the filter itself
(called ``quantum trajectories'' in the physics literature).  In addition,
the Markov property is key for the development of feedback controls, as we
will see in sections \ref{sec:fb}--\ref{sec:lyapunov}.

We will consider below the homodyne detection case, but the photodetection
case proceeds identically.  Set $\Delta y_l=\iota(\Delta Y(l))$.  Recall
that the homodyne detection signal $\Delta y_l$ takes one of two values
$\pm\lambda$ for every $l$.  Suppose we have observed $y_i$ up to and
including time $l-1$; we would like to be able to calculate the
probability distribution of $\Delta y_l$ using this information. This
calculation is carried out in the following lemma.

\begin{lemma}  \label{lemma:Delta-y}
We have
\begin{equation}
	\mathbf{P}[ \Delta y_l = \pm \lambda \, \vert \, \mathcal{Y}_{l-1} ]
	= p(\Delta y_l = \pm \lambda;  \rho_{l-1} ) ,
\label{eq:Delta-y-1}
\end{equation}
where
\begin{equation}
	p(\Delta y = \pm \lambda; \rho  ) =
	\frac{1}{2}\pm\frac{\lambda}{2}\,
                {\rm Tr}[\mathcal{\overline{J}}(\rho)   ]
\label{eq:Delta-y-2}
\end{equation}
depends only on the filter in the previous time step.
\end{lemma}

\begin{proof}
Let $p^+_l = \mathbf{P}[ \Delta y_l = \pm \lambda \, \vert \, 
\mathcal{Y}_{l-1} ]$
be the probability that the observation in the next time step $\Delta 
y_l$ takes the value $+\lambda$. Using the Doob decomposition we have
$$
        \Delta Y(l) = \Delta C(l)+\Delta\tilde Y(l)=
        \pi_{l-1}(\mathcal{J}(I))\,\Delta t(l)+\Delta\tilde Y(l),
$$
where $\Delta\tilde Y(l)$ is a martingale increment, so that
$$
        \lambda\,p^+_l-\lambda(1-p^+_l)=
        \iota(\mathbb{P}(\Delta Y(l)|\mathscr{Y}_{l-1}))=
        \iota(\pi_{l-1}(\mathcal{J}(I)))\,\lambda^2.
$$
Thus  
$$
        p^+_l
        =\frac{1}{2}+\frac{\lambda}{2}\,
                {\rm Tr}[\rho_{l-1}\mathcal{J}(I)]
        =\frac{1}{2}+\frac{\lambda}{2}\,
                {\rm Tr}[\mathcal{\overline{J}}(\rho_{l-1})]
$$
depends only on the filter in the previous time step. 
\qquad
\end{proof}

With this distribution in hand, we can also calculate the statistics of
the filter in the next time step, and we can prove the Markov property by
recursing this procedure.

\begin{lemma}\label{lem:markov}
The filter $\rho_l$ satisfies the Markov property:
$$
	\mathbf{E_P}(g(\rho_j)|\sigma\{\rho_0,\ldots,\rho_{l}\})
	= \mathbf{E_P}(g(\rho_j)|\sigma\{\rho_{l}\})\qquad
	\forall\,l\le j\le k.
$$
\end{lemma}

\begin{proof}
Using the recursion (\ref{eq homo density form}), $\rho_j$ can be written 
as a deterministic function $f(\rho_{j-1},\Delta y_j)$ of $\rho_{j-1}$ and 
$\Delta y_j$.  By the martingale property of the innovation process
$$
	\mathbf{P}(\Delta y_j=\pm\lambda|\sigma\{
		y_1,\ldots,y_{j-1}\})
	=\frac{1}{2}\pm\frac{\lambda}{2}\,
		{\rm Tr}[\mathcal{\overline{J}}(\rho_{j-1})],
$$
which is only a function of $\rho_{j-1}$.  As 
$\sigma\{\rho_0,\ldots,\rho_{j-1}\}\subset\sigma\{y_1,\ldots,y_{j-1}\}$,
we obtain
$$
	\mathbf{P}(\Delta y_j=\pm\lambda|\sigma\{
		\rho_0,\ldots,\rho_{j-1}\})
	=\frac{1}{2}\pm\frac{\lambda}{2}\,
		{\rm Tr}[\mathcal{\overline{J}}(\rho_{j-1})].
$$
Hence for any function $g$
\begin{multline*}
	\mathbf{E_P}(g(\rho_j)|\sigma\{\rho_0,\ldots,\rho_{j-1}\})
	=
	\mathbf{E_P}(g'(\rho_{j-1},\Delta y_j)|
		\sigma\{\rho_0,\ldots,\rho_{j-1}\}) 
	\\ =
	\sum_{a\in\{-1,1\}} g'(\rho_{j-1},a\lambda)\left[
	\frac{1}{2}+a\frac{\lambda}{2}\,{\rm Tr}[
        \mathcal{\overline{J}}(\rho_{j-1})]\right],
\end{multline*}
where we have written $g'(\rho,w)=g(f(\rho,w))$ for the function $g$ 
composed with the one-step filter recursion.  As the right hand side is a 
function of $\rho_{j-1}$ only, we have
$$
	\mathbf{E_P}(g(\rho_j)|\sigma\{\rho_0,\ldots,\rho_{j-1}\})=
	\mathbf{E_P}(\mathbf{E_P}(g(\rho_j)|\sigma\{\rho_0,\ldots,\rho_{j-1}\})
	\,|\,\sigma\{\rho_{j-1}\}),
$$
from which we conclude that
$$
	\mathbf{E_P}(g(\rho_j)|\sigma\{\rho_0,\ldots,\rho_{j-1}\})=
	\mathbf{E_P}(g(\rho_j)|\sigma\{\rho_{j-1}\}).
$$
But setting $\mathbf{E_P}(g(\rho_j)|\sigma\{\rho_0,\ldots,\rho_{j-1}\})=
h(\rho_{j-1})$, we can repeat the argument giving
$$
	\mathbf{E_P}(h(\rho_{j-1})|\sigma\{\rho_0,\ldots,\rho_{j-2}\})
	=
	\mathbf{E_P}(h(\rho_{j-1})|\sigma\{\rho_{j-2}\}).
$$
From the definition of $h(\cdot)$ we immediately obtain
$$
	\mathbf{E_P}(g(\rho_j)|\sigma\{\rho_0,\ldots,\rho_{j-2}\})
	=
	\mathbf{E_P}(g(\rho_j)|\sigma\{\rho_{j-2}\}).
$$
Recursing the argument gives the Markov property.
\qquad
\end{proof}

It is now straightforward to turn this procedure into a Monte Carlo 
algorithm.  The following pseudocode generates random sample paths of the 
observations $\Delta y_l$ and filter $\rho_l$, sampled faithfully from the 
probability measure induced by the repeated interaction model on the space 
of observation paths.
\vskip.2cm
{\sc
\begin{enumerate}
\item Initialize $\rho_0$.
\item $l\leftarrow 0$.
\item Repeat
\begin{enumerate}
\item Calculate $p^+_{l+1}(\rho_l)$.
\item Sample $\xi\sim\mbox{Uniform}[0,1]$.
\item If $\xi<p^+_{l+1}$: $\Delta y_{l+1}\leftarrow +\lambda$;
	Else: $\Delta y_{l+1}\leftarrow -\lambda$.
\item $\rho_{l+1}\leftarrow \rho_l+\Delta\rho_{l+1}(\rho_l,\Delta y_{l+1})$.
\item $l\leftarrow l+1$.
\end{enumerate}
\item Until $l=k$.
\end{enumerate}
}

\subsection{Examples}\label{sec:examplesfilter}

Using the Monte Carlo method developed in the previous section, we can now 
simulate the observation and filter sample paths for our usual examples.
In the simulations we have used the initial state $\rho(X)={\rm Tr}[X]/2$,
under which the probabilities of the energy attaining its maximal or 
minimal values are equal.

\subsubsection*{Spontaneous emission}

\begin{figure}
\centering
\includegraphics[width=\textwidth]{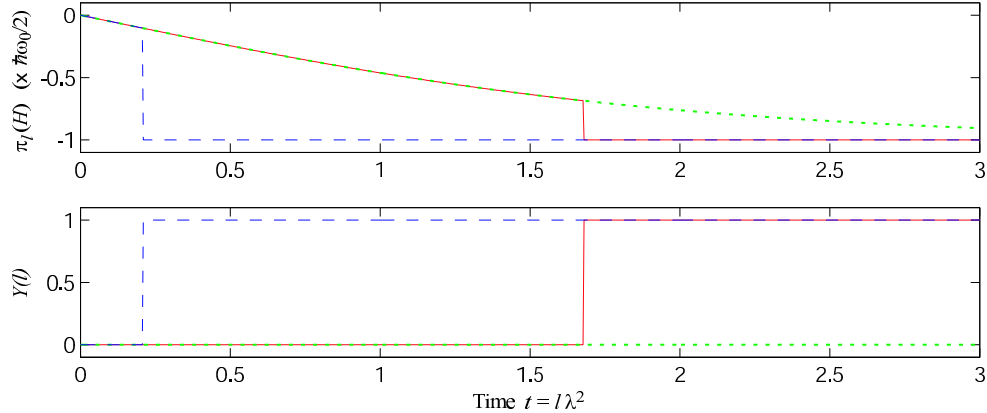}
\caption{\label{fig:photospont}  
{\rm(\bf Spontaneous emission\rm)}  Three typical sample paths of the 
conditional expectation $\pi_l(H)$ of the atomic energy (top plot) with 
respect to the integrated photodetection signal $Y(l)$ (bottom plot).  The 
atom spontaneously emits exactly one photon at a random time, after which 
it attains its lowest energy and remains there.  The initial state chosen 
was $\rho(x)={\rm Tr}[X]/2$; in this state $H=\pm\hbar\omega_0/2$ have 
equal probability.  The time scale used for the calculation is 
$\lambda^{-2}=300$.} 
\end{figure}

In Figure \ref{fig:photospont}, photodetection of a spontaneously emitting
atom is simulated.  The observation process takes a very simple form: if
the atom initially attained its maximal energy (which it does with
probability $\tfrac{1}{2}$ under the state $\rho$), it emits a single
photon at a random time.  If the atom was already at its lowest energy
(also with probability $\tfrac{1}{2}$), the atom never emits a photon. The
conditional expectation of the energy attains its minimal value
immediately after the observation of a photon, as at this point we know
that the atom has attained its lowest energy.  Before the observation of a
photon, the conditional expectation decays: the longer we fail to observe
a photon, the higher the (conditional) probability that the atom started
out with minimal energy to begin with.

Note that the higher the initial expectation of the energy of the atom,
the slower the decay of $\pi_l(H)$: after all, if the probability of the
atom starting out at its minimal energy is very small, then we should fail
to observe a photon for a very long time before concluding that the atom,
against all odds, did start off with minimal energy.  In the extreme case
of unit initial probability that the atom has maximal energy
($\rho(X)={\rm Tr}[\sigma_+\sigma_-X]$), the conditional expectation of
the energy is a step function.  The latter can be verified directly using
the filtering equation for photodetection, which shows that
$\Delta\pi_l(H)=0$ as long as $\Delta Y(l)=0$ and
$\pi_{l-1}(H)=\hbar\omega_0/2$.

\begin{figure}
\centering
\includegraphics[width=\textwidth]{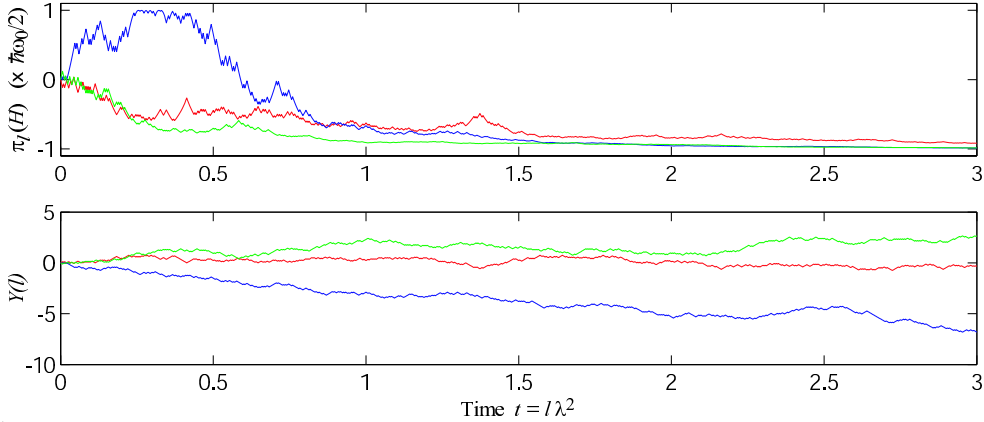}
\caption{\label{fig:homospont} 
{\rm(\bf Spontaneous emission\rm)}  Three typical sample paths of the 
conditional expectation $\pi_l(H)$ of the atomic energy (top plot) with 
respect to the integrated homodyne photocurrent $Y(l)$ (bottom plot).  
It is difficult to infer much about the atomic energy from the 
observation process using the naked eye, but nonetheless the conditional 
expectation is revealing (e.g.\ in the case of the blue sample path, the 
atom initially possessed maximal energy with unit probability).  The 
initial state chosen was $\rho(x)={\rm Tr}[X]/2$, and the time scale used 
for the calculation is $\lambda^{-2}=300$.} 
\end{figure}

Figure \ref{fig:homospont} shows a simulation of the same system, but now 
observed using a homodyne detection setup.  Evidently the way in which
information on the atomic energy is encoded in the homodyne observations 
is very different than in the photodetection case; rather than the sudden 
gain of information when a photon is observed, homodyne detection 
allows us to gradually infer the atomic energy from the noisy 
measurements.  The observation process itself is not very revealing to the 
naked eye, but the filter manages to make sense of it.  In the case of the 
blue curve, for example, we infer that the atom almost certainly attained 
its maximal energy value before time $t\sim 0.5$, whereas after time 
$t\sim 1.5$ all three paths indicate that most likely the atom has 
attained its lowest energy.

\subsubsection*{Dispersive interaction}

The dispersive interaction case is quite different than the spontaneous
emission case.  Recall that in this case the energy observable $j_l(H)=H$
is constant in time.  This does not mean, however, that the conditional
expectation $\pi_l(H)$ is constant.  Whether the energy attains its
maximal or minimal value in a particular realization determines the mean
direction of the phase shift of the outgoing light, which can be measured
using a homodyne detection setup.  As we gain information on the atomic
energy, the conditional expectation gradually attracts to the value
actually taken by the energy $H$ in that realization.

\begin{figure}
\centering
\includegraphics[width=\textwidth]{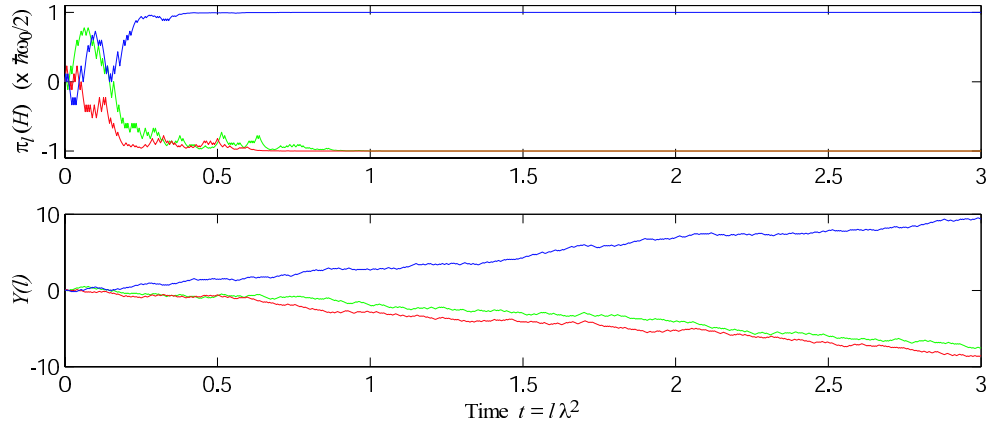}
\caption{\label{fig:homodisp} 
{\rm(\bf Dispersive interaction\rm)}  Three typical sample paths of the 
conditional expectation $\pi_l(H)$ of the atomic energy (top plot) with 
respect to the integrated homodyne photocurrent $Y(l)$ (bottom plot).  
Though the energy observable is constant in time, the observations contain 
information on the atomic energy so that the \emph{conditional} expectation 
converges to the actual value of the energy.  The initial state chosen was 
$\rho(x)={\rm Tr}[X]/2$, and the time scale used is $\lambda^{-2}=300$.} 
\end{figure}

This behavior is demonstrated in the simulation of Figure
\ref{fig:homodisp}.  Each of the filter sample paths $\pi_l(H)$ attracts
to either $+\hbar\omega_0/2$ or $-\hbar\omega_0/2$.  As the probabilities
of maximal and minimal atomic energy are equal under the initial state
$\rho(X)={\rm Tr}[X]/2$, the sample paths of $\pi_l(H)$ are attracted to
$\pm\hbar\omega_0/2$ with equal probability.  It would be presumptious to
conclude this from the simulation of only three sample paths, but if we
are willing to believe that $\pi_l(H)\to\pm\hbar\omega_0/2$ with unit
probability then the result is evident: after all,
$\mathbb{P}(\pi_l(H))=\mathbb{P}(j_l(H))=\rho(H)=0$ for all $l$, so we can
only have $\pi_l(H)\to\pm\hbar\omega_0/2$ if the two possibilities occur
with equal probability (the probabilities change accordingly if we choose
a different initial state $\rho$).  The fact that $\pi_l(H)\to\pm
\hbar\omega_0/2$ with unit probability can also be rigorously proved, but
we will postpone this discussion until section \ref{sec:lyapunov}.

Note that the behavior of the filter can already be seen by inspecting the
observation process using the naked eye: though the observation processes
are still random, the integrated observations have an upward or downward
trend depending on the value of the atomic energy.  This indicates that
the actual observation process $\Delta Y(l)$ is positive (negative) on
average if the atomic energy is positive (negative), i.e., positive atomic
energy leads to an average positive phase shift on the output light,
whereas negative atomic energy gives rise to a negative average phase
shift.

We have not shown a simulation of the dispersively interacting atom under 
photodetection.  Though many photons can be observed, a photodetector 
gives no information on the phase of the output light.  Hence no 
information is gained about the atomic energy, and the conditional 
expectation of the atomic energy is constant.  Evidently the type of 
detector used makes a big difference in this case.

\section{The reference probability approach}\label{reference probability}

In the previous section we obtained explicit nonlinear recursive equations
for the conditional expectations $\pi_l(X)$.  In this section we start
from scratch and solve the filtering problem in an entirely different way.
The key idea here is that the Bayes Lemma \ref{lem Bayes} allows us to
express the filtering problem in terms of an arbitrary (reference) state; 
by choosing the ``Radon-Nikodym'' operator $V$ conveniently, we can reduce 
the calculation to a manageable form.  This gives rise to a \emph{linear} 
recursion for the numerator in the Bayes formula, which we denote by 
$\sigma_l(X)$.  $\pi_l(X)$ is then calculated as $\sigma_l(X)/\sigma_l(I)$.

The reference probability method is widely used in classical filtering 
theory following the work of Duncan \cite{TD68}, Mortensen \cite{RM66}, 
Zakai \cite{Zak69}, and Kallianpur and Striebel \cite{KalStr68,KalStr69}.  
See e.g.\ \cite{AEM95} for a systematic exposition of reference 
probability methods in filtering and control.  The corresponding approach 
in quantum filtering theory, adapted below to the discrete setting, was 
developed in \cite{BH05}.

\subsection{The strategy}

Let us begin by outlining what we are going to do.  We are looking for a
way to calculate $\pi_l(X)=\mathbb{P}(j_l(X)|\mathscr{Y}_l)$.  In
classical filtering theory it has proved advantageous to express the
problem, using the Bayes formula, in terms of a measure under which the
signal (here the atomic observable $j_l(X)$) and the observations (here
$\mathscr{Y}_l$) are independent.  We will aim to do the same in the
quantum case.  Unlike in classical probability, however, we do not have a
suitable quantum version of Girsanov's theorem to help us find such a
change of measure.  We need to use a little intuition to obtain the
required change-of-state operator $V$.

The following simple Lemma will make this task a little easier.

\begin{lemma}\label{lem:rotatedstatecondex}
Let $(\mathscr{A},\mathbb{P})$ be a quantum probability space, 
$\mathscr{C}\subset\mathscr{A}$ a commutative subalgebra, and 
$U\in\mathscr{A}$ be unitary.  Define the rotated state
$\mathbb{Q}(X)=\mathbb{P}(U^*XU)$ on $\mathscr{A}$.  Then
$\mathbb{P}(U^*XU|U^*\mathscr{C}U)=U^*\mathbb{Q}(X|\mathscr{C})U$ for any 
$X\in\mathscr{C}'$.
\end{lemma}

\begin{proof}
This is a simple consequence of the definition of conditional 
expectations.  Let $K\in U^*\mathscr{C}U$.  Then
$\mathbb{P}(U^*\mathbb{Q}(X|\mathscr{C})UK)=
\mathbb{Q}(\mathbb{Q}(X|\mathscr{C})UKU^*)=\mathbb{Q}(XUKU^*)=
\mathbb{P}(U^*XUK)$.  But as this holds for any $K\in U^*\mathscr{C}U$ and 
as $U^*\mathbb{Q}(X|\mathscr{C})U\in U^*\mathscr{C}U$, the Lemma follows 
from the definition of the conditional expectation.
\qquad
\end{proof}

How does this help us in our usual setting $(\mathscr{M}\otimes
\mathscr{W}_k,\mathbb{P})$?  Recall from section \ref{sec:repeatedinteraction}
that $\mathscr{C}_l={\rm alg}\{\Delta Z(i):i=1,\ldots,l\}$ 
($Z=I\otimes\Lambda$ or $Z=I\otimes(A+A^*)$ for photodetection or homodyne 
detection, respectively), and $\mathscr{Y}_l=U(l)^*\mathscr{C}_lU(l)$. For 
the time being, let us fix a time step $l$ and define the state 
$\mathbb{Q}^l(X)=\mathbb{P}(U(l)^*XU(l))$. Then by Lemma 
\ref{lem:rotatedstatecondex}, we can write $\pi_l(X)=
U(l)^*\mathbb{Q}^l(X|\mathscr{C}_l)U(l)$.  Now note that $X\in\mathscr{B}_0$
is of the form $X\otimes I$ in $\mathscr{M}\otimes\mathscr{W}_k$, whereas 
every element of $\mathscr{C}_l$ is of the form $I\otimes C$.  But we 
already know a state under which the initial system and the field are 
independent: this is simply the state $\mathbb{P}$!  Hence if we could 
write $\mathbb{Q}^l(X)=\mathbb{P}(V(l)^*XV(l))$ for some 
$V(l)\in\mathscr{C}_l'$, then we would obtain using the Bayes Lemma 
\ref{lem Bayes}
\begin{equation}\label{eq:kalstrgen}
	\pi_l(X)=U(l)^*\,\mathbb{Q}^l(X|\mathscr{C}_l)\,U(l)=
		\frac{U(l)^*\,\mathbb{P}(V(l)^*XV(l)|\mathscr{C}_l)\,U(l)}
		{U(l)^*\,\mathbb{P}(V(l)^*V(l)|\mathscr{C}_l)\,U(l)}
	\qquad \forall\,X\in\mathscr{B}_0.
\end{equation}
Note that we already have by construction 
$\mathbb{Q}^l(X)=\mathbb{P}(U(l)^*XU(l))$, but $U(l)\not\in\mathscr{C}_l'$.
Hence $V(l)=U(l)$ in the Bayes formula does not work.  As we will 
demonstrate below, however, there is a simple trick which we can use to 
``push'' $U(l)$ into the commutant $\mathscr{C}_l'$ without changing the 
state $\mathbb{Q}^l$.  This gives the desired change of state $V(l)$.

We emphasize that Lemma \ref{lem:rotatedstatecondex} is not an essential
part of the procedure; we could try to find a new reference state
$\mathbb{S}(V^*XV)=\mathbb{P}(X)$ and apply the Bayes lemma directly to
$\mathbb{P}(j_l(X)|\mathscr{Y}_l)$ (in fact, the state $\mathbb{S}$ can be
read off from Eq.\ (\ref{eq:kalstrgen}) by applying Lemma
\ref{lem:rotatedstatecondex} in the reverse direction).  The state
$\mathbb{P}$ is a very convenient reference state, however, as it allows
us to use the natural tensor splitting $\mathscr{M}\otimes\mathscr{W}_k$
and properties of the vacuum state $\phi^{\ten k}$ to determine the
necessary $V(l)$ with minimal effort. It is possible that this procedure
could be streamlined in a more general theory for changes of state (in the
spirit of the techniques widely used in classical probability theory), but
such a theory is not currently available in quantum probability.

\subsection{Homodyne detection}\label{sec zakai homo}

We first consider the homodyne detection case, i.e.\ $Z(l)=A(l)+A^*(l)$.
To use Eq.\ (\ref{eq:kalstrgen}), we are looking for
$V(l)\in\mathscr{C}_l'$ such that
$\mathbb{Q}^l(X)=\mathbb{P}(V(l)^*XV(l))$.  The following trick
\cite{Hol90} allows us to find such a $V(l)$ simply by modifying the
quantum stochastic difference equation for $U(l)$, Eq.\ (\ref{eq U}).

\begin{lemma}\label{lem nondemo}
Let $V$ be the solution of the following difference equation:
	\nomenclature{$V(l)$}{Change of state operator for reference measure}
\begin{equation*}
	\Delta V(l)=
	\left\{
		M^+\,(\Delta A(l)+\Delta A^*(l))+M^\circ\,\Delta t(l)
	\right\}V(l-1),\qquad
	V(0)=I.
\end{equation*}
Then $V(l)\in\mathscr{C}_l'$ and 
$\mathbb{P}(V(l)^*XV(l))=\mathbb{P}(U(l)^*XU(l))$ for any
$X \in \mathscr{M}\otimes\mathscr{W}_k$.  
\end{lemma}

\begin{proof}
The proof relies on the fact that $\mathbb{P}=\rho\otimes\phi^{\ten k}$.
For simplicity, let us assume that $\rho$ is a vector state, i.e.\
$\rho(X)=\langle v,Xv\rangle$ for some $v\in\mathbb{C}^2$ (we will relax 
this requirement later on).  Then $\mathbb{P}(X)=\langle 
v\otimes\Phi^{\ten k},X\,v\otimes\Phi^{\ten k}\rangle$ where $\Phi$ is the 
vacuum vector.  The essential property we need is that $\Delta 
A(l)\Phi^{\ten k}=\Delta\Lambda(l)\Phi^{\ten k}=0$: this follows 
immediately from their definition (as $\Delta A(l)\Phi^{\ten k}=\Phi^{\ten 
l-1}\ten \sigma_-\Phi\ten \Phi^{\ten k-l}=0$, and similarly for 
$\Delta\Lambda$).  Hence using the fact that $U(l-1)\in\mathscr{B}_{l-1}$ 
commutes with $\Delta A(l)$ and $\Delta\Lambda(l)$, Eq.\ (\ref{eq U}) gives
$$
	\Delta U(l)\,v\otimes\Phi^{\ten k}=
	\left\{
		M^+\,\Delta A^*(l)+M^\circ\,\Delta t(l)
	\right\}U(l-1)\,v\otimes\Phi^{\ten k}.
$$
Similarly, \emph{any} difference equation of the form
$$
	\Delta V(l)=\left\{
		N^\pm\,\Delta\Lambda(l)+M^+\,\Delta A^*(l)
		+N^-\,\Delta A(l)+M^\circ\,\Delta t(l)
	\right\}V(l-1)
$$
satisfies
$$
	\Delta V(l)\,v\otimes\Phi^{\ten k}=\left\{
		M^+\,\Delta A^*(l)+M^\circ\,\Delta t(l)
	\right\}V(l-1)\,v\otimes\Phi^{\ten k}.
$$
Hence if $V(l-1)\,v\otimes\Phi^{\ten k}=U(l-1)\,v\otimes\Phi^{\ten k}$,
then $V(l)\,v\otimes\Phi^{\ten k}=U(l)\,v\otimes\Phi^{\ten k}$.  By 
induction $V(l)\,v\otimes\Phi^{\ten k}=
U(l)\,v\otimes\Phi^{\ten k}$ for any $l$ if $V(0)=U(0)=I$.  Thus
\begin{multline*}
	\mathbb{P}(U(l)^*XU(l))=
	\langle U(l)\,v\otimes\Phi^{\ten k},X\,U(l)\,
		v\otimes\Phi^{\ten k}\rangle \\
	=\langle V(l)\,v\otimes\Phi^{\ten k},X\,V(l)\,v\otimes\Phi^{\ten k}\rangle=
	\mathbb{P}(V(l)^*XV(l)),
\end{multline*}
regardless of what we choose for $N^\pm$ and $N^-$.

We are now free to choose $N^\pm$ and $N^-$ so that $V(l)$ satisfies the 
remaining requirement $V(l)\in\mathscr{C}_l'$.  But if we choose $N^\pm=0$ 
and $N^-=M^+$ as in the statement of the Lemma, it follows that
$V(l)\in\mathscr{M}\otimes\mathscr{C}_l\subset\mathscr{C}_l'$ for every 
$l$.  Indeed, suppose that $V(l-1)\in\mathscr{M}\otimes\mathscr{C}_{l-1}$.
Then $V(l)$ is defined by the recursion as a function of $V(l-1)$, $\Delta 
Z(l)=\Delta A(l)+\Delta A^*(l)\in\mathscr{C}_l$ and $M^+,M^\circ\in
\mathscr{M}$, which is obviously contained in $\mathscr{M}\otimes
\mathscr{C}_l$.  The result then follows by induction, and the Lemma is 
proved.

It remains to consider the case that $\rho$ is not a vector state.  
By linearity we can always write $\rho={\rm Tr}[\tilde\rho X]$ for some 
$2\times 2$ density matrix $\tilde\rho$.  But any density matrix can be 
diagonalized (as it is a positive matrix), so that we can write without 
loss of generality $\tilde\rho=\lambda_1\,v_1v_1^*+\lambda_2\,v_2v_2^*$ 
where $v_1$ is the (normalized) eigenvector of $\tilde\rho$ with eigenvalue 
$\lambda_2$, and $v_2$ is the eigenvector with eigenvalue $\lambda_2$.
As the Lemma holds for each of the vector states $\rho_{1,2}(X)=\langle 
v_{1,2},Xv_{1,2}\rangle$, it must hold for arbitrary $\rho$.
\qquad
\end{proof}

The solution of the filtering problem is now remarkably straightforward.

\begin{definition}\label{def unnormalized}
For any atomic observable $X\in\mathscr{B}_0$, the unnormalized 
conditional expectation $\sigma_l(X)\in \mathscr{Y}_l$ is defined as
\begin{equation*}
  \sigma_l(X) = U(l)^*\,\mathbb{P}(V(l)^*XV(l)|\mathscr{C}_l)\,U(l).
\end{equation*}
\end{definition}
	\nomenclature{$\sigma_l(X)$}{Unnormalized conditional state}

\begin{theorem}\label{thm:zakaihomodyne}
The unnormalized conditional expectation $\sigma_l(X)$ 
satisfies the following {\bf linear filtering equation for homodyne 
detection}:
\begin{equation*}
	\Delta\sigma_l(X)=
	\sigma_{l-1}(\mathcal{L}(X))\,\Delta t(l)+
	\sigma_{l-1}(\mathcal{J}(X))\,\Delta Y(l),\qquad
	\sigma_0(X)=\rho(X),
\end{equation*}
where $\mathcal{L}(X)$ is the discrete Lindblad generator and 
$\mathcal{J}(X)$ was defined in Lemma \ref{lem xi homo}.
Furthermore the noncommutative Kallianpur-Striebel formula holds:
\begin{equation*}
  \pi_l(X) = \frac{\sigma_l(X)}{\sigma_l(I)}
	\qquad\forall\, X \in \mathscr{B}_0.
\end{equation*}
\end{theorem}

\begin{proof}
The Kallianpur-Striebel formula is simply Eq.\ (\ref{eq:kalstrgen}).  To 
obtain the linear recursion, we calculate using the discrete It\^o rules
\begin{multline*}
	\Delta(V(l)^*XV(l)) = V(l-1)^*\mathcal{L}(X)V(l-1)\,\Delta t(l)
	\\+V(l-1)^*\mathcal{J}(X)V(l-1)\,(\Delta A(l)+\Delta A^*(l)).
\end{multline*}
Calculating the conditional expectation with respect to $\mathscr{C}_l$, 
we obtain
\begin{multline*}
	\Delta\mathbb{P}(V(l)^*XV(l)|\mathscr{C}_l) 
	= \mathbb{P}(V(l-1)^*\mathcal{L}(X)V(l-1)|\mathscr{C}_l)\,\Delta t(l)
	\\+
	\mathbb{P}(V(l-1)^*\mathcal{J}(X)V(l-1)|\mathscr{C}_l)\,
		(\Delta A(l)+\Delta A^*(l))
\end{multline*}
using $\Delta A(l)+\Delta A^*(l)\in\mathscr{C}_l$.
But note that $\mathscr{C}_l=\mathscr{C}_{l-1}\otimes{\rm alg}\{
\Delta A(l)+\Delta A^*(l)\}$ and that $\Delta A(l)+\Delta A^*(l)$ is 
independent from $\mathscr{M}\otimes\mathscr{C}_{l-1}$ under $\mathbb{P}$.
Hence by the independence property of the conditional expectation (cf.\ 
the last property listed in Table \ref{tab:condex}) and the fact that
$V(l-1)^*XV(l-1)\in\mathscr{M}\otimes\mathscr{C}_{l-1}$ for any 
$X\in\mathscr{B}_0$, we obtain
\begin{multline*}
	\Delta\mathbb{P}(V(l)^*XV(l)|\mathscr{C}_l) 
	= \mathbb{P}(V(l-1)^*\mathcal{L}(X)V(l-1)|\mathscr{C}_{l-1})\,\Delta t(l)
	\\+
	\mathbb{P}(V(l-1)^*\mathcal{J}(X)V(l-1)|\mathscr{C}_{l-1})\,
		(\Delta A(l)+\Delta A^*(l)).
\end{multline*}
Now multiply from the left by $U(l)^*$ and from the right by $U(l)$.
Note that 
$$
	U(l)^*\,\Delta\mathbb{P}(V(l)^*XV(l)|\mathscr{C}_l)\,U(l)=
	\Delta (U(l)^*\,\mathbb{P}(V(l)^*XV(l)|\mathscr{C}_l)\,U(l)),
$$
because $U(l)^*CU(l)=U(l-1)^*CU(l-1)$ for any $C\in\mathscr{C}_{l-1}$ (see 
section \ref{sec:repeatedinteraction}).  Furthermore $U(l)^*(\Delta 
A(l)+\Delta A^*(l))U(l)=\Delta Y(l)$, and the Theorem follows.
\qquad
\end{proof}

Notice how much simpler the linear recursion is compared to the nonlinear 
recursion obtained through the martingale method; for example, this could 
make digital implementation of the linear recursion more straightforward.
Nonetheless the two methods should give the same answer, as they are both 
ultimately expressions for the same quantity $\pi_l(X)$.  Let us verify 
this explicitly.  Note that
$$
	\Delta\pi_l(X)=\Delta\left[\frac{\sigma_l(X)}{\sigma_l(I)}
	\right]=\frac{\sigma_l(X)}{\sigma_l(I)}-
	\frac{\sigma_{l-1}(X)}{\sigma_{l-1}(I)}=
	\frac{\Delta\sigma_l(X)}{\sigma_{l-1}(I)}-
	\sigma_l(X)\,(\sigma_{l-1}(I)^{-1}-\sigma_l(I)^{-1}).
$$
The first term on the right is easily evaluated as
$$
	\frac{\Delta\sigma_l(X)}{\sigma_{l-1}(I)}=
	\pi_{l-1}(\mathcal{L}(X))\,\Delta t(l)+
	\pi_{l-1}(\mathcal{J}(X))\,\Delta Y(l),
$$
whereas we obtain for the second term (using $\mathcal{L}(I)=0$)
$$
	\sigma_l(X)\,(\sigma_{l-1}(I)^{-1}-\sigma_l(I)^{-1})=
	\sigma_l(X)\,\frac{\Delta\sigma_l(I)}{\sigma_l(I)\sigma_{l-1}(I)}
	=\pi_l(X)\,\pi_{l-1}(\mathcal{J}(I))\,\Delta Y(l).
$$
Hence we obtain
$$
	\Delta\pi_l(X)=
	\pi_{l-1}(\mathcal{L}(X))\,\Delta t(l)+
	\left\{\pi_{l-1}(\mathcal{J}(X))
		-\pi_l(X)\,\pi_{l-1}(\mathcal{J}(I))
	\right\}\Delta Y(l).
$$
Writing $\pi_l(X)=\pi_{l-1}(X)+\Delta\pi_l(X)$ and solving for 
$\Delta\pi_l(X)$ gives precisely the nonlinear recursion obtained in 
section \ref{sec filter homo}, taking into account the identity
$$
	(I+\pi_{l-1}(\mathcal{J}(I))\,\Delta Y(l))^{-1}=
	\frac{I-\pi_{l-1}(\mathcal{J}(I))\,\Delta Y(l)}
	{I-\lambda^2\pi_{l-1}(\mathcal{J}(I))^2}
$$
where we have used $(\Delta Y(l))^2=\Delta t(l)$.  This sheds some light 
on the seemingly complicated structure of the discrete nonlinear filtering 
equation.

\subsection{Photodetection}

If we naively try to follow the procedure above in the photodetection case 
$Z(l)=\Lambda(l)$, we run into problems.  Let us see what goes wrong.  
Following the steps in the proof of Lemma \ref{lem nondemo}, we reach the 
point where we have to choose $N^\pm$ and $N^-$ so that 
$V(l)\in\mathscr{C}_l'$.  But this is impossible\footnote{
	This is not surprising for the following reason.  Applying the 
	spectral theorem to the commutative algebra $\mathscr{C}_l$, we 
	obtain a classical measure space $(\Omega,\mathcal{F})$ on which
	$\mathbb{Q}^l$ and $\mathbb{P}$ induce different probability 
	measures $\mathbf{Q}$ and $\mathbf{P}$, respectively.  Suppose 
	there exists a $V\in\mathscr{C}_l'$ such that 
	$\mathbb{Q}^l(X)=\mathbb{P}(V^*XV)$.  Then it is not difficult to 
	verify that $\mathbf{Q}\ll\mathbf{P}$ with 
	$d\mathbf{Q}/d\mathbf{P}=\iota(\mathbb{P}(V^*V|\mathscr{C}_l))$.
	But note that $\iota(\Delta\Lambda(l))$ is distributed under
	$\mathbf{Q}$ in the same way as the observation increment $\Delta 
	Y(l)$ under $\mathbb{P}$, whereas $\iota(\Delta\Lambda(l))=0$ 
	$\mathbf{P}$-a.s.  This contradicts $\mathbf{Q}\ll\mathbf{P}$
	and hence the existence of $V$.
} because we can not get rid of the $\Delta A^*$ term (as $\Delta 
A^*(l)\Phi^{\ten k}\ne 0$), and $\Delta A^*(l)$ does not commute with 
$\Delta\Lambda(l)$.  

We are not restricted, however, to using Lemma \ref{lem:rotatedstatecondex}
with $U=U(l)$.  To deal with the photodetection case, suppose that $R(l)$ 
is some unitary operator of the form $I\otimes R$ in $\mathscr{M}\otimes
\mathscr{W}_k$ (i.e.\ it only acts on the field, not on the atom).  
Define the state $\mathbb{Q}^l(X)=\mathbb{P}(U(l)^*R(l)^*XR(l)U(l))$ and 
the algebra $\mathscr{R}_l=R(l)\mathscr{C}_lR(l)^*$.  Suppose that there 
is some $V(l)\in\mathscr{R}_l'$ such that $\mathbb{Q}^l(X)=
\mathbb{P}(V(l)^*XV(l))$.  Then by Lemma \ref{lem:rotatedstatecondex} and 
the Bayes formula, we have for any $X\in\mathscr{B}_0$ (using
$R(l)^*XR(l)=XR(l)^*R(l)=X$)
$$
	\pi_l(X)=U(l)^*R(l)^*\,\mathbb{Q}^l(X|\mathscr{R}_l)\,R(l)U(l)
	=\frac{
		U(l)^*R(l)^*\,\mathbb{P}(V(l)^*XV(l)|\mathscr{R}_l)\,R(l)U(l)
	}{
		U(l)^*R(l)^*\,\mathbb{P}(V(l)^*V(l)|\mathscr{R}_l)\,R(l)U(l)
	}.
$$
This is precisely as before, except that we have inserted an additional 
rotation $R(l)$.  The idea is that if we choose $R(l)$ appropriately, then
$R(l)\,(\Delta\Lambda(l))\,R(l)^*$ contains a $\Delta A^*$ term so that we 
can proceed as in the previous section to find a suitable $V(l)\in
\mathscr{R}_l'$.  A possible choice of $R(l)$ is given in the following 
Lemma.

\begin{lemma}\label{lem:weylrr}
Define $R(l)=\exp(A(l)-A^*(l))$, $l=1,\ldots,k$.  Then $R(l)$ is unitary,
$$
  \Delta R(l) = \left\{\frac{\sin(\lambda)}{\lambda}\,(\Delta A(l) - 
  \Delta A^*(l)) + 
  \frac{\cos(\lambda)-1}{\lambda^2}\,\Delta t(l)\right\}R(l-1),
$$
and we obtain the expression
$$
	R(l)\,\Delta\Lambda(l)\,R(l)^*=
	 	(\cos^2\lambda-\sin^2\lambda)\, \Delta \Lambda(l) + 
		\frac{\sin\lambda\,\cos\lambda}{\lambda}
			\,(\Delta A(l) + \Delta A^*(l)) + 
	 	\frac{\sin^2\lambda}{\lambda^2}\,\Delta t(l).  
$$
\end{lemma}

\begin{proof}
Unitarity of $R(l)$ is immediate.  To obtain the difference equation,
note that $R(l)=\exp(\Delta A(l)-\Delta A^*(l))\,R(l-1)$, so that we 
essentially need to calculate $\exp(\Delta A(l)-\Delta A^*(l))=
(\exp(\lambda\sigma_--\lambda\sigma_+))_l$.  But it is not difficult to 
evaluate explicitly the matrix exponential of the $2\times 2$ matrix
$\lambda(\sigma_--\sigma_+)$:
$$
	\exp(\lambda(\sigma_--\sigma_+))=
	\exp\left[\begin{pmatrix}
		0 & -\lambda \\
		\lambda & 0
	\end{pmatrix}
	\right]=
	\begin{pmatrix}
		\cos\lambda & -\sin\lambda\\
		\sin\lambda & \cos\lambda
	\end{pmatrix}=
	\cos\lambda\,I+\sin\lambda\,(\sigma_--\sigma_+).
$$
The expression for $\Delta R(l)$ follows directly.  To obtain the 
remaining expression, note that $R(l-1)$ commutes with $\Delta\Lambda(l)$; 
hence 
\begin{multline*}
	R(l)\,(\Delta\Lambda(l))\,R(l)^*=
	e^{\Delta A(l)-\Delta A^*(l)}\,\Delta\Lambda(l)\,
	e^{\Delta A^*(l)-\Delta A(l)} = \\ 
	\left(
	\begin{pmatrix}
		\cos\lambda & -\sin\lambda\\
		\sin\lambda & \cos\lambda
	\end{pmatrix}
	\begin{pmatrix}
		1 & 0 \\
		0 & 0
	\end{pmatrix}
	\begin{pmatrix}
		\cos\lambda & \sin\lambda\\
		-\sin\lambda & \cos\lambda
	\end{pmatrix}		
	\right)_l=
	\left(
	\begin{pmatrix}
		\cos^2\lambda & \sin\lambda\,\cos\lambda\\
		\sin\lambda\,\cos\lambda & \sin^2\lambda
	\end{pmatrix}		
	\right)_l
\end{multline*}
from which the result follows immediately.
\qquad
\end{proof}

For brevity, let us call $R(l)\,(\Delta\Lambda(l))\,R(l)^*=\Delta 
\overline{Z}(l)$.  We are now in a position to repeat Lemma \ref{lem 
nondemo} for the photodetection case.

\begin{lemma}
Let $V$ be the solution to the following difference equation:
\begin{multline*}
	\Delta V(l)=
	\left\{
		\left(\frac{\lambda}{\sin\lambda}\,\,M^+ 
			-\frac{\lambda^2}{\cos\lambda}\,M^\circ
			- \frac{1}{\cos\lambda}\right)
		\left(\Delta\overline{Z}(l)
			-\frac{\sin^2\lambda}{\lambda^2}\,\Delta t(l)
		\right)
		\right. \\
		+ \left.
		\left(\frac{\sin\lambda}{\lambda}\,M^+
			+\cos\lambda\,M^\circ+\frac{\cos\lambda-1}{\lambda^2}
		\right)
		\Delta t(l)
	\right\}V(l-1),\qquad
	V(0)=I.
\end{multline*}
Then $V(l)\in\mathscr{R}_l'$ and $\mathbb{P}(V(l)^*XV(l))=
\mathbb{P}(U(l)^*R(l)^*XR(l)U(l))$ $\forall\,X\in\mathscr{M}\otimes
\mathscr{W}_k$.
\end{lemma}

\begin{proof}
Using the quantum It\^o rules, it is not difficult to calculate
\begin{multline*}
	\Delta(R(l)U(l))=
	\left\{
		(\cdots)\,\Delta\Lambda(l)
		+\left(\cos\lambda\,M^+ 
			-\lambda\sin\lambda\,M^\circ
			- \frac{\sin\lambda}{\lambda}\right)
			\Delta A^*(l)
		\right. \\
		+(\cdots)\,\Delta A(l)
		+ \left.
		\left(\frac{\sin\lambda}{\lambda}\,M^+
			+\cos\lambda\,M^\circ+\frac{\cos\lambda-1}{\lambda^2}
		\right)
		\Delta t(l)
	\right\}R(l-1)U(l-1).
\end{multline*}
We can now follow exactly the procedure in the proof of Lemma \ref{lem 
nondemo}, and using the expression for $\Delta\overline{Z}$ obtained in 
Lemma \ref{lem:weylrr} the result follows.
\qquad
\end{proof}

We can now obtain the linear filter for photodetection as before.

\begin{theorem}\label{thm zakai count}
The unnormalized conditional expectation 
$$
	\sigma_l(X)=U(l)^*R(l)^*\,
	\mathbb{P}(V(l)^*XV(l)|\mathscr{R}_l)\,R(l)U(l),
	\qquad X\in\mathscr{B}_0,
$$ 
satisfies the {\bf linear filtering equation for photodetection}
\begin{equation*}
	\Delta\sigma_l(X)=
	\sigma_{l-1}(\mathcal{L}(X))\,\Delta t(l)+
	\sigma_{l-1}(\mathcal{T}(X))\left(\Delta Y(l)
		-\frac{\sin^2\lambda}{\lambda^2}\,\Delta t(l)
	\right),~~~
	\sigma_0(X)=\rho(X),
\end{equation*}
where $\mathcal{L}(X)$ is the discrete Lindblad generator and 
$\mathcal{T}(X)$ is given by
	\nomenclature{$\mathcal{T}(X)$}{$\mathcal{T}$-coefficient of the 
		filtering equations}
$$
	\mathcal{T}(X)=
	\frac{1}{\cos^2\lambda}\left(
		\frac{\lambda^2}{\sin^2\lambda}\,M^{+*}XM^+ - X
		-\lambda^2\mathcal{L}(X)
	\right).
$$
Furthermore the noncommutative Kallianpur-Striebel formula holds:
\begin{equation*}
  \pi_l(X) = \frac{\sigma_l(X)}{\sigma_l(I)}
	\qquad\forall\, X \in \mathscr{B}_0.
\end{equation*}
\end{theorem}

\begin{proof}
We begin by calculating $\Delta(V(l)^*XV(l))$ using the It\^o rules; this 
is a tedious but straightforward calculation, and we will not repeat it 
here.  The result is
\begin{multline*}
	\Delta(V(l)^*XV(l))=V(l-1)^*\mathcal{L}(X)V(l-1)\,\Delta t(l)
	\\+
	V(l-1)^*\mathcal{T}(X)V(l-1)\left(\Delta\overline{Z}(l)
                -\frac{\sin^2\lambda}{\lambda^2}\,\Delta t(l)
        \right),
\end{multline*}
where $\mathcal{T}(X)$ is given by the expression
$$
	\mathcal{T}(X)=\frac{\lambda^2}{\sin^2\lambda}\,M^{+*}XM^+
		-\frac{\lambda^4}{\cos^2\lambda}\,M^{\circ *}XM^\circ
		-\frac{\lambda^2}{\cos^2\lambda}\,(XM^\circ+M^{\circ *}X)
		-\frac{1}{\cos^2\lambda}\,X.
$$
Using the expression for $\mathcal{L}(X)$, the latter is easily 
transformed into the form given in the Theorem.  The remainder of 
the proof is the same as that of Theorem \ref{thm:zakaihomodyne}.
\qquad
\end{proof}

The presence of the $\sin\lambda$ and $\cos\lambda$ terms in the linear 
filter for photodetection is an artefact of our choice of $R(l)$.  In 
fact, this choice is not unique: many $R(l)$ would work and give rise to 
different linear filters!  However, the Kallianpur-Striebel formula 
guarantees that all these linear filters coincide when normalized with the 
nonlinear filter of section \ref{sec filter count}.  To complete our 
discussion of linear filters, let us show explicitly how normalization of 
the linear filter of Theorem \ref{thm zakai count} gives back, through a 
series of miraculous cancellations, the nonlinear filter for 
photodetection which we obtained earlier through martingale methods.
We begin by writing, as in section \ref{sec zakai homo},
$$
	\Delta\pi_l(X)=\frac{\Delta\sigma_l(X)}{\sigma_{l-1}(I)}
		-\pi_l(X)\,\frac{\Delta\sigma_l(I)}{\sigma_{l-1}(I)}.
$$
Using the fact that $\pi_l(X)=\Delta\pi_l(X)+\pi_{l-1}(X)$,
this gives explicitly
\begin{multline*}
	\Delta\pi_l(X)\left[
		I+\pi_{l-1}(\mathcal{T}(I))\,(\Delta Y(l)-\sin^2\lambda)
	\right] 
\\=
	\pi_{l-1}(\mathcal{L}(X))\,\Delta t(l)
	-\frac{\sin^2\lambda}{\lambda^2}\left[
		\pi_{l-1}(\mathcal{T}(X))-
		\pi_{l-1}(X)\,\pi_{l-1}(\mathcal{T}(I))
	\right]\Delta t(l) 
\\
	+\left[
		\pi_{l-1}(\mathcal{T}(X))-
		\pi_{l-1}(X)\,\pi_{l-1}(\mathcal{T}(I))
	\right]
	\Delta Y(l).
\end{multline*}
Using the expression for $\mathcal{T}(X)$ and $\mathcal{L}(I)=0$, we 
calculate
\begin{multline*}
	\pi_{l-1}(\mathcal{T}(X))-\pi_{l-1}(X)\,\pi_{l-1}(\mathcal{T}(I))	
	=
	-\frac{\lambda^2}{\cos^2\lambda}\,\pi_{l-1}(\mathcal{L}(X))
\\
	+\frac{\lambda^2}{\sin^2\lambda\,\cos^2\lambda}
	\left[
		\pi_{l-1}(M^{+*}XM^+)-\pi_{l-1}(X)\,\pi_{l-1}(M^{+*}M^+)
	\right].
\end{multline*}
Hence we obtain
\begin{multline*}
	\Delta\pi_l(X)
	\left[
		I+\pi_{l-1}(\mathcal{T}(I))\,(\Delta Y(l)-\sin^2\lambda)
	\right] 
\\=
	\frac{1}{\cos^2\lambda}\left[
		\pi_{l-1}(\mathcal{L}(X))
		-\pi_{l-1}(M^{+*}XM^+)+\pi_{l-1}(X)\,\pi_{l-1}(M^{+*}M^+)
	\right]\Delta t(l) 
\\
	+\left[
		\pi_{l-1}(\mathcal{T}(X))-
		\pi_{l-1}(X)\,\pi_{l-1}(\mathcal{T}(I))
	\right]
	\Delta Y(l).
\end{multline*}
Next, we claim that
\begin{multline*}
	\left[
		I+\pi_{l-1}(\mathcal{T}(I))\,(\Delta Y(l)-\sin^2\lambda)
	\right]^{-1}
\\=
	\left[
		I+\cos^2\lambda\,\pi_{l-1}(\mathcal{T}(I))
	\right]^{-1}\,\Delta Y(l)+
	\left[
		I-\sin^2\lambda\,\pi_{l-1}(\mathcal{T}(I))
	\right]^{-1}\,(I-\Delta Y(l)).
\end{multline*}
The easiest way to see this is to consider $\pi_{l-1}(\mathcal{T}(I))$ and 
$\Delta Y(l)$ to be classical random variables through the spectral 
theorem; as $(\Delta Y(l))^2=\Delta Y(l)$ we conclude that $\Delta Y(l)$ 
is a $\{0,1\}$-valued random variable, and the statement follows directly.
Using the explicit expression for $\mathcal{T}(X)$, we find that
$$
	\left[
                I-\sin^2\lambda\,\pi_{l-1}(\mathcal{T}(I))
        \right]^{-1}=
	\frac{\cos^2\lambda}{I-\lambda^2\,\pi_{l-1}(M^{+*}M^+)},
$$
and that
$$
	\left[
                I+\cos^2\lambda\,\pi_{l-1}(\mathcal{T}(I))
        \right]^{-1}=
	\frac{\sin^2\lambda}{\lambda^2\,\pi_{l-1}(M^{+*}M^+)}.
$$
Using these expressions, the remainder of the calculation is a
straightforward exercise and indeed we obtain the expression for 
$\Delta\pi_l(X)$ as in section \ref{sec filter count}.


\section{Feedback control}
\label{sec:fb}

Everything we have done up to this point has been devoid of human
intervention.  An atom sits in free space, emits radiation at its leisure,
and all we have allowed ourselves to do is to observe the output radiation
and to interpret it statistically (filtering).  In this section we will
allow ourselves to manipulate the atom in real time; this provides an
opportunity for feedback control, which we will approach using optimal
control theory (section \ref{sec:optimal}) and Lyapunov methods (section
\ref{sec:lyapunov}).

In this section we carefully introduce the concepts required for feedback
control. As we shall see, some subtleties arise which we explain and
address using the framework developed in \cite{BH05}, suitably adapted to
our discrete setting.  Early work on quantum feedback control appears in
\cite{Bel83,Bel88}, and a large number of applications have been discussed
(albeit not in a mathematically rigorous way) in the physics literature,
see e.g.\ \cite{DHJMT00} and the references therein.

In order to understand feedback, we need to consider how the atomic
dynamics can be influenced by control actions, and what information is
available to determine the control actions.  We develop these ideas in the
following two subsections. In subsection \ref{sec:fb-open} we describe how
a controller can influence the atomic dynamics.  This will allow us to
apply open loop controls, that is, we can apply a deterministic function
to the control inputs of the system.  Though this is not our ultimate
goal, it is a helpful first step; in subsection \ref{sec:fb-closed} we
will show how to replace this deterministic input by some function of the
observation history (feedback). This is likely to be advantageous: the
more information we have, the better we can control!

\begin{remark}
We take a moment at this point to discuss the usage of the term ``quantum 
control'' in the literature.  Often this term is used to refer to 
open-loop control, rather than feedback control.  Such control problems 
can be reduced to deterministic control problems, as we will show in 
subsection \ref{sec:fb-open}.  The classical analog of this concept would 
be deterministic control design for the Fokker-Planck equation.

Our main goal here is to discuss quantum {\it feedback} control, where the
feedback is based on an observations process obtained from the system to
be controlled.  This corresponds to the classical idea of a control system
as consisting of a system to be controlled (the plant); a sensor which
gives rise to an observations process; an actuator which allows one to
modify the dynamics of the system in real time; and a controller, which is
a signal processing device that takes the observations as its input and
produces an actuation signal as its output.

Sometimes the term quantum feedback control is used in a somewhat broader
context.  One could consider the setup described above in absence of the
sensor component.  In this case the controller must be treated as being
itself a physical system, rather than a signal processing device, and the
pursuit of this idea leads to a rather different theory and applications 
(see e.g.\ \cite{WM94b,YK03,JP06a}).  This type of feedback is usually 
called {\it coherent} feedback or {\it all-optical} feedback (in the 
context of quantum optics), to distinguish it from observation-based 
feedback which we consider here.
\end{remark}

\subsection{Open loop control (no feedback)}
\label{sec:fb-open}

To add a control input to our model, recall from section 
\ref{sec:repeatedinteraction} that the time evolution of an observable $X$ 
is given by $j_l(X)=U(l)^*XU(l)$, with repeated interaction unitary
$$
	U(l)=\overrightarrow{\prod}_{i=1}^l M(i)=
		M(1)M(2)\cdots M(l),\qquad
	U(0)=I.
$$
The interaction in every time slice was given by
$$
	M(l)=
	e^{-i\{j_{l-1}(L_1)\Delta\Lambda(l)
	+j_{l-1}(L_2)\Delta A^*(l)
	+j_{l-1}(L_2^*)\Delta A(l)
	+j_{l-1}(L_3)\Delta t(l)\}
	},
$$
where $L_1,L_2,L_3\in\mathscr{B}_0$ are atomic operators.   To add a
control input, we simply allow ourselves a choice of different $M(l)$'s in 
every time step.  To be precise, let us introduce a \emph{control set} 
$\mathfrak{U}$
	\nomenclature{$\mathfrak{U}$}{Control set}
(i.e.\ these are the values the control can take), and for 
each $u\in\mathfrak{U}$ we define a set of atomic operators 
$L_1(u),L_2(u),L_3(u)\in\mathscr{B}_0$.
An \emph{open loop control strategy} is a $k$-tuple
$\mathbf{u}=\{u_1,\ldots,u_k\}$ where $u_i\in\mathfrak{U}$ for all $i$, 
and the corresponding time evolution $j_l^{\bf u}(X)=U(l,{\bf u})^*X
U(l,{\bf u})$ is defined by
	\nomenclature{$j_l^{\bf u}(X)$}{Time evolution with open loop control}
	\nomenclature{$U(l,{\bf u})$}{Interaction unitary with open loop control}
$$
	U(l,\mathbf{u})=\overrightarrow{\prod}_{i=1}^l M^{\bf u}(i,u_i)=
		M^{\bf u}(1,u_1)M^{\bf u}(2,u_2)\cdots M^{\bf u}(l,u_l),\qquad
	U(0,\mathbf{u})=I,
$$
where the single time slice interaction unitary is given by
$$
	M^{\bf u}(l,u)=
	e^{-i\{j_{l-1}^{\bf u}(L_1(u))\Delta\Lambda(l)
	+j_{l-1}^{\bf u}(L_2(u))\Delta A^*(l)
	+j_{l-1}^{\bf u}(L_2(u)^*)\Delta A(l)
	+j_{l-1}^{\bf u}(L_3(u))\Delta t(l)\}
	}.
$$
Note that $M^{\bf u}(l,u)$, for fixed $u$, depends only on ${\bf u}$ 
through $u_1,\ldots,u_{l-1}$, and similarly $U(l,{\bf u})$ depends only on 
$u_1,\ldots,u_l$.  We write $U(l,{\bf u})$ rather than $U(l,u_1,\ldots,u_l)$
purely for notational convenience; the latter would technically be more 
appropriate!

As before, it is convenient to run the definition backwards in time, i.e.
$$
	U(l,\mathbf{u})=\overleftarrow{\prod}_{i=1}^l M_i(u_i)=
		M_l(u_l)M_{l-1}(u_{l-1})\cdots M_1(u_1),\qquad
	U(0,\mathbf{u})=I,
$$
where $M_l(u)$ is given by
        \nomenclature{$M_l(u)$}{Controlled single time step interaction, time reversed}
\begin{equation}\label{eq:openloopmlu}
	M_l(u)=
	e^{-i\{L_1(u)\Delta\Lambda(l)
	+L_2(u)\Delta A^*(l)
	+L_2(u)^*\Delta A(l)
	+L_3(u)\Delta t(l)\}
	}.
\end{equation}
These operators are functions of the current control value, and do not depend on the full sequence $\mathbf{u}$.
The corresponding difference equation is written as
\begin{multline}
	\Delta U(l,\mathbf{u})=
	\left\{
		M^\pm(u_{l})\,\Delta\Lambda(l)
		+M^+(u_{l})\,\Delta A^*(l) 
\right. \\ \left.
		+M^-(u_{l})\,\Delta A(l)
		+M^\circ(u_{l})\,\Delta t(l)
	\right\}U(l-1,\mathbf{u}),\qquad
	U(0,\mathbf{u})=I.
	\label{controlled-U}
\end{multline}
One could imagine the different controls $u\in\mathfrak{U}$ to correspond 
to different values of a magnetic field which is applied to the atom and 
can be changed by the experimenter in each time step.  Another common 
control input is obtained using a laser whose amplitude can be controlled 
by the experimenter.  We will discuss a specific example in the context of 
the discrete model in section \ref{sec:optimal-eg}.

The question of controller design could already be posed at this
deterministic level.  Though this rules out feedback control (as the
observations are a random process), open loop controls for quantum systems
have already generated important applications.  For example, optimal
control theory has allowed the design of time-optimal pulse sequences for
nuclear magnetic resonance (NMR) spectroscopy that significantly
outperform the state-of-the-art for that technology \cite{navin}.  

Open loop control design is beyond the scope of this article, so we will 
not go into detail.  Let us take a moment, however, to make a connection 
with this literature.  Let $X\in\mathscr{B}_0$ be some atomic operator; 
then using the quantum It\^o rules, we easily establish as in section
\ref{sec:qsde} that the expectation of the controlled time evolution 
$j_l^{\bf u}(X)$ in the case of \emph{deterministic} $\bf u$ satisfies
\begin{equation}\label{eq:controlledlindblad}
	\frac{\Delta\,\mathbb{P}(j_l^{\bf u}(X))}{\Delta t}=
	\mathbb{P}(j_{l-1}^{\bf u}(\mathcal{L}(X,u_{l}))),
\end{equation}
where the controlled Lindblad generator is given by
\begin{equation*}
	\mathcal{L}(X,u) =
		M^{+}(u)^*XM^+(u) +
		\lambda^2\,M^{\circ}(u)^*XM^\circ(u)
		+M^{\circ}(u)^*X+XM^\circ(u).
\end{equation*}
As $\mathbb{P}(j_l^{\bf u}(X))$ is linear in $X$, we can introduce a 
(deterministic) $2\times 2$ matrix $\nu_l^{\bf u}$ such that 
$\mathbb{P}(j_l^{\bf u}(X))={\rm Tr}[\nu_l^{\bf u}X]$; then Eq.\ 
(\ref{eq:controlledlindblad}) can be written as a deterministic recursion 
for $\nu_l^{\bf u}$, known as the (controlled, discrete) \emph{master 
equation}:
$$
	\frac{\Delta\nu_l^{\bf u}}{\Delta t}=
		M^{+}(u_{l})\nu_{l-1}^{\bf u}M^+(u_{l})^* +
		\lambda^2\,M^{\circ}(u_{l})\nu_{l-1}^{\bf u}
			M^\circ(u_{l})^*
		+M^{\circ}(u_{l})\nu_{l-1}^{\bf u}
		+\nu_{l-1}^{\bf u}M^\circ(u_{l})^*.
$$
The control goal in such a scenario is generally formulated as the desire 
to choose $\bf u$ so that the expectation of a certain atomic observable, 
or a nonlinear function of such expectations, is maximized.  But these 
expectations can be obtained in closed form by solving the master 
equation, so that the control problem reduces to a \emph{deterministic} 
optimal control problem for the master equation.  This is precisely the 
sort of problem that is solved in \cite{navin} (in a continuous time 
context).

\subsection{Closed loop control  (feedback)}
\label{sec:fb-closed}

We now turn to the issue of letting the control at time step $l$ be a
function of the observation history prior to time $l$.  Mathematically, it
is not entirely clear how to modify the description in the previous
section to allow for feedback in the repeated interaction model; there are
in fact some subtleties.  It is the goal of this section to clear up this
point.  For notational simplicity we assume from this point onwards that
$\mathfrak{U}\subset\mathbb{R}$, i.e.\ that our control input is a real
scalar (this is not essential; everything we will do can be generalized).

Let us first consider what we mean by a controller.  A controller is 
a signal processing device---a black box---that on every time step takes 
an observation $\Delta y$ as its input and generates a control signal in 
$\mathfrak{U}$ as its output.  The most general control strategy $\mu$ is 
thus described by a set of functions, one for each time step:
	\nomenclature{$\mu,\mu^*,\bar\mu$}{Control strategies}
	\nomenclature{$f_l(\cdots)$}{Feedback function}
$$
	\mu=\{f_1,f_2(\Delta y_1),f_3(\Delta y_1,\Delta y_2),\ldots,
		f_k(\Delta y_1,\ldots,\Delta y_{k-1})\}.
$$
Causality is enforced explicitly by making the control in time step $l$ a 
function only of the observations up to that time step.  The feedback 
control $\mu$ is called {\it admissible} if $f_l$ takes values in 
$\mathfrak{U}$ for every $l$.  We call the set of all admissible controls 
$\mathfrak{K}$.
	\nomenclature{$\mathfrak{K}$}{Set of admissible controls}

The functions $f_l$ encode the input-output behavior of the controller 
implementing the admissible strategy $\mu\in\mathfrak{K}$.  We now need 
to hook up the controller to our model of the system.  As the observation 
at time $i$ is described by the observable $\Delta Y_i$, we define the 
output of the controller at time $l$ as the observable
	\nomenclature{$\mathfrak{u}_l$}{Control signal $f_l(\Delta 
		Y(1),\ldots,\Delta Y(l-1))$}
$$
	\mathfrak{u}_l=f_l(\Delta Y(1),\ldots,\Delta Y(l-1))
	\in\mathscr{Y}_{l-1}.
$$
Clearly $\iota(\mathfrak{u}_l)$ is a random variable that takes values in 
$\mathfrak{U}$, precisely as it should be.

It now remains to close the loop, i.e., to make the time evolution a
function of the output $\mathfrak{u}_l$ of the controller.  Formally, we 
can proceed exactly as in the open loop case; i.e., we define the time 
evolution $j_l^\mu(X)=U^\mu(l)^*XU^\mu(l)$ by
	\nomenclature{$j_l^{\mu}(X)$}{Time evolution with feedback control}
	\nomenclature{$U^\mu(l)$}{Interaction unitary with feedback control}
$$
	U^\mu(l)=\overrightarrow{\prod}_{i=1}^l M^\mu(i,\mathfrak{u}_i)=
	M^\mu(1,\mathfrak{u}_1)M^\mu(2,\mathfrak{u}_2)\cdots 
	M^\mu(l,\mathfrak{u}_l),\qquad
	U^\mu(0)=I,
$$
where the single time slice interaction unitary is given by
        \nomenclature{$M^\mu(l,u)$}{Controlled single time step interaction}
$$
	M^\mu(l,u)=j_{l-1}^\mu(M_l(u)).
$$
To make this precise, however, we need to define what we mean by composing 
the unitary operator-valued function $M^\mu(l,u)$ with the observable 
$\mathfrak{u}_l$.  We will thus take a moment to talk about such
compositions.

\subsubsection*{Composition of an operator-valued function and an 
observable}

Let us begin with a simple classical analogy.  Let $\mathfrak{u}$ be a 
random variable in $\ell^\infty(\Omega,\mathcal{F},\mathbf{P})$, where 
$\Omega$ is a finite set.  Then $\mathfrak{u}$ takes a finite number of 
values, and we will suppose these lie in a set $\mathfrak{U}$.  Let 
$M:\mathfrak{U}\to V$ be a map from $\mathfrak{U}$ to some linear space 
$V$.  We would like to define $M(\mathfrak{u})$ as a $V$-valued random 
variable.  We could do this as follows:
\begin{equation} \label{eq:classcomp}
	M(\mathfrak{u})(\omega)=
	\sum_{u\in\mathrm{ran}\,\mathfrak{u}}M(u)\,
		\chi_{\{\mathfrak{u}=u\}}(\omega),
\end{equation}
where $\mbox{ran}\,\mathfrak{u}$ is the range of $\mathfrak{u}$ and
$\chi_{\{\mathfrak{u}=u\}}$ is the indicator function on
$\{\omega:\mathfrak{u}(\omega)=u\}$.  This is particularly convenient if
we would like to think of $\ell^\infty(\Omega,\mathcal{F},\mathbf{P})$ as
being itself a linear space; by elementary linear algebra the set of
$V$-valued random variables is isomorphic to
$V\otimes\ell^\infty(\Omega,\mathcal{F},\mathbf{P})$, and the definition
of $M(\mathfrak{u})$ above only involves sums and tensor products
($M(u)\otimes\chi_{\{\mathfrak{u}=u\}}$) which are naturally defined in
this space.

Though this is obviously not the simplest way of defining composition in 
the classical case, (\ref{eq:classcomp}) looks the most natural in the 
noncommutative context.  Let us consider the algebra 
$\mathscr{A}\otimes\mathscr{C}$, where $\mathscr{C}$ is commutative but 
$\mathscr{A}$ not necessarily so.  We can think of 
$\mathscr{A}\otimes\mathscr{C}$ as a linear space of $\mathscr{A}$-valued 
random variables; indeed, applying the spectral theorem to $\mathscr{C}$ 
we find that $\mathscr{A}\otimes\mathscr{C}\simeq\mathscr{A}\otimes
\ell^{\infty}(\mathcal{F})\simeq\ell^\infty(\mathcal{F};\mathscr{A})$.
Now suppose we are given a map $M:\mathfrak{U}\to\mathscr{A}$ and an 
observable $\mathfrak{u}\in\mathscr{C}$ such that $\iota(\mathfrak{u})$ 
takes values in $\mathfrak{U}$.  Then it is natural to define the 
composition $M(\mathfrak{u})$ as an element in 
$\mathscr{A}\otimes\mathscr{C}$ in the same way as (\ref{eq:classcomp}).
This motivates the following definition.

\begin{definition}\label{def:Masafunctionofu}
Let $\mathscr{A}$, $\mathscr{C}$ be $^*$-algebras where $\mathscr{C}$ is 
commutative.  Let $\mathfrak{U}\subset\mathbb{R}$, 
$M:\mathfrak{U}\to\mathscr{A}$ and let $\mathfrak{u}\in\mathscr{C}$ be 
such that ${\rm sp}\,\mathfrak{u}={\rm ran}\,\iota(\mathfrak{u})$, the 
spectrum of $\mathfrak{u}$, is a subset of $\mathfrak{U}$.  Then the 
composition $M(\mathfrak{u})\in\mathscr{A}\otimes\mathscr{C}$ is defined by
\begin{equation*}
	M(\mathfrak{u}) = \sum_{u \in {\rm sp}(\mathfrak{u})} 
	M(u)\,P_{\mathfrak{u}}(u),
\end{equation*} 
where $P_{\mathfrak{u}}(u)=\iota^{-1}(\chi_{\{\iota(\mathfrak{u})=u\}})$ 
is the eigenspace projector of $\mathfrak{u}$ for eigenvalue $u$.    
\end{definition}

In what follows it will be important to understand how the composition 
$M(\mathfrak{u})$ behaves under unitary transformations.  Consider 
$\mathscr{A}\otimes\mathscr{C}\subset\mathscr{B}$, and let $U$ be a 
unitary operator in $\mathscr{B}$.  Consider the map $M_U:\mathfrak{U}\to 
U^*\mathscr{A}U$ defined by $M_U(u) = U^*M(u)U$ for all $u \in \mathfrak{U}$.   
Then $U^*M(\mathfrak{u})U\in U^*(\mathscr{A}\otimes\mathscr{C})U=
U^*\mathscr{A}U\otimes U^*\mathscr{C}U$ is given by
\begin{multline}\label{eq:transformationU}
  U^*M(\mathfrak{u})U = \sum_{u \in \mathrm{sp}(\mathfrak{u})} 
  U^*M(u)U\,U^*P_{\mathfrak{u}}(u)U = \\
  \sum_{u \in \mathrm{sp}(\mathfrak{u})} 
  U^*M(u)U\,P_{U^*\mathfrak{u}U}(u) =
  M_U(U^*\mathfrak{u}U).
\end{multline}
Hence unitary rotations preserve the compositions defined above, as long 
as we remember to rotate both the observable $\mathfrak{u}$ and the map 
$M(u)$.

\subsubsection*{Controlled quantum flows}

Let us return to the controlled time evolution 
$U^\mu(l)$.  Note that $M^\mu(l,\cdot):\mathfrak{U}\to 
j_{l-1}^\mu(\mathscr{M}\otimes\mathscr{M}_l)$, 
whereas $\mathfrak{u}_l\in \mathscr{Y}_{l-1}=j_{l-1}^\mu(\mathscr{C}_{l-1})$.
Hence according to Definition \ref{def:Masafunctionofu}, the composition 
$M(\mathfrak{u}_l)$ makes sense as an operator in the algebra 
$j_{l-1}^\mu(\mathscr{M}\otimes\mathscr{C}_{l-1}\otimes\mathscr{M})
\subset\mathscr{B}_l$, and by (\ref{eq:transformationU})
$$
	M(l,\mathfrak{u}_l)=
	M(l,f(\Delta Y(1),\ldots,\Delta Y(l-1))) =
	j_{l-1}^\mu(M_l(f(\Delta Z(1),\ldots,\Delta Z(l-1)))).
$$
For brevity, we will write
	\nomenclature{$\mathfrak{\check u}_l$}{Control process $f_l(\Delta 
		Z(1),\ldots,\Delta Z(l-1))$}
$$
	\mathfrak{\check u}_l=f_l(\Delta Z(1),\ldots,\Delta Z(l-1))
	\in\mathscr{C}_{l-1}
$$
so that $M(l,\mathfrak{u}_l)=j_{l-1}^\mu(M_l(\mathfrak{\check u}_l))$.
In particular, we can now express $U^\mu(l)$ as
\begin{equation}\label{eq:umulclosed}
	U^\mu(l)=\overleftarrow{\prod}_{i=1}^l M_i(\mathfrak{\check u}_i)=
		M_l(\mathfrak{\check u}_l)
		M_{l-1}(\mathfrak{\check u}_{l-1})\cdots 
		M_1(\mathfrak{\check u}_1),\qquad
	U^\mu(0)=I,
\end{equation}
which gives the controlled quantum stochastic difference equation
\begin{multline}
	\Delta U^\mu(l)=
	\left\{
		M^\pm(\check{\mathfrak{u}}_{l})\,\Delta\Lambda(l)
		+M^+(\check{\mathfrak{u}}_{l})\,\Delta A^*(l) 
\right. \\ \left.
		+M^-(\check{\mathfrak{u}}_{l})\,\Delta A(l)
		+M^\circ(\check{\mathfrak{u}}_{l})\,\Delta t(l)
	\right\}U^\mu(l-1),\qquad
	U^\mu(0)=I.
	\label{controlled-U-fb}
\end{multline}
The main thing to note is that in order to close the loop in 
(\ref{controlled-U}), the open loop controls $u_l$ should be replaced by
$\mathfrak{\check u}_l\in\mathscr{C}_{l-1}$ rather than $\mathfrak{u}_l\in
\mathscr{Y}_{l-1}$.  This ensures that the corresponding flow $j^\mu_l(X)$ 
depends on the observations history in the right way, by virtue of 
(\ref{eq:transformationU}).  Other than this subtlety, the treatment of 
closed loop time evolution proceeds much along the same lines as in the 
absence of feedback.

The notion of a controlled quantum flow \cite{BH05} summarizes these ideas
in a general context.  The following definition is a discrete version of
this concept and defines a generic repeated interaction model with
scalar feedback.

\begin{definition}[\bf Controlled quantum flow] \label{def:cqf}
The quadruple $(\mathfrak{U},M,Z,\mu)$ s.t.
\begin{enumerate}
\item $\mathfrak{U}\subset\mathbb{R}$,
\item $M_l:\mathfrak{U}\to\mathscr{M}\otimes\mathscr{M}_l$, $M_l(u)$ is a 
unitary operator of the form {\rm (\ref{eq:openloopmlu})} 
$\forall\,u\in\mathfrak{U}$,
\item  $Z$ is an adapted process $Z(l)\in\mathscr{W}_l$, $l=1,\ldots,k$ 
such that $Z(l)$ is self-adjoint and $\mathscr{C}_l={\rm 
alg}\{Z(i):i=1,\ldots,l\}$ is commutative for every $l$,
\item $\mu\in\mathfrak{K}$ is an admissible control strategy,
\end{enumerate}
defines a controlled quantum flow $j^\mu_l(X)=U^\mu(l)^*XU^\mu(l)$. Here
$U^\mu(l)$ is given by {\rm (\ref{eq:umulclosed})} and the corresponding 
observations process $Y^\mu(l)$ is given by
	\nomenclature{$Y^\mu(l)$}{Observations process for control 
		strategy $\mu$}
$$
	\Delta Y^\mu(l)=U^\mu(l)^*\,\Delta Z(l)\,U^\mu(l),\qquad
	l=1,\ldots,k.
$$
\end{definition}

\begin{remark} 
For illustrative purposes, we will concentrate in the following on the
homodyne detection case $Z=A+A^*$; the theory for photodetection 
$Z=\Lambda$ proceeds along the same lines. The notion of a controlled
quantum flow is much more general, however, and even allows for feedback
to the detection apparatus.  For example, recall that a homodyne detection
setup can measure any of the processes $e^{i\varphi}A+e^{-i\varphi}A^*$.  
We could now make $\varphi$ a function of the past observations (i.e.\
$\varphi(l)=\tilde f_l(\Delta Z(1),\ldots, \Delta Z(l-1))$), thus feeding
back to the homodyne detector; this fits within the framework of the
controlled quantum flow as it just requires us to use a ``nonlinear''
$Z$.  Feedback to the detector has proven useful for sensitive 
measurements of the phase of an optical pulse, see
\cite{wisemanadaptive,AASDM02}.  We will not consider this further here,
but everything we will discuss can be adapted to this case as well.  
We encourage the reader to work out the following sections for the
feedback scenario of his choice!
\end{remark}

\subsection{Filtering in the presence of feedback}

Now that we have resolved how to model a quantum system with feedback, the
next question to be resolved is whether filtering still works in this
context. Fortunately this is indeed the case, and in fact little changes
in the proofs.  The resulting filters are completely intuitive: one
obtains the same filter from a controlled quantum flow as one would obtain
by first calculating the filter with an open loop control, then
substituting the feedback law into the filtering equation.  In this
section we will briefly discuss filtering in the presence of feedback
using the reference probability method.  From this point onwards we 
restrict ourselves to the homodyne detection case $Z=A+A^*$.  

Fix an admissible feedback control strategy $\mu\in\mathfrak{K}$ and the
corresponding control observables $\mathfrak{u}_l$, $\mathfrak{\check
u}_l$.  We wish to find an expression for $\pi_l^\mu(X)=
\mathbb{P}(j_l^\mu(X)|\mathscr{Y}^\mu_l)$.  For this to make sense we have
to make sure that the self-nondemolition and nondemolition properties
still hold.  If they do not then we did something wrong (recall that these
properties are essential for a meaningful interpretation of the theory);  
but let us take a moment to verify that everything is as it should be.

\begin{lemma}
The observation algebra $\mathscr{Y}_l^\mu$ is commutative
(self-nondemolition) and $j_l^\mu(X)\in(\mathscr{Y}_l^{\mu})'$
(nondemolition) for every $l=1,\ldots,k$ and 
$X\in\mathscr{M}\otimes\mathscr{C}_l$.
\end{lemma}

\begin{proof}
The unitary $M_l(\mathfrak{\check u}_l)$, by construction, commutes 
with every element of $\mathscr{C}_{l-1}$.  Hence $U^\mu(l)^*\,\Delta 
Z(i)\,U^\mu(l)=\Delta Y^\mu(i)$ for every $i\le l-1$, and we find that
$$
	[\Delta Y^\mu(i),\Delta Y^\mu(l)]=
	U^\mu(l)^*[\Delta Z(i),\Delta Z(l)]U^\mu(l)=0.
$$
This establishes self-nondemolition.  Nondemolition is established similarly.
\qquad
\end{proof}

To apply the reference probability method, we need a suitable 
change of state.  The proof of the following Lemma is omitted as it is 
identical to that of Lemma \ref{lem nondemo}.

\begin{lemma}  \label{lem nondemo u}
Let $V^\mu$ be the solution of the following difference equation:
\begin{equation*}
	\Delta V^\mu(l)=
	\left\{
		M^+(\mathfrak{\check u}_l)\,(\Delta A(l)+\Delta A^*(l))
		+M^\circ(\mathfrak{\check u}_l)\,\Delta t(l)
	\right\}V^\mu(l-1),\qquad
	V^\mu(0)=I.
\end{equation*}
Then $V^\mu(l)\in\mathscr{C}_l'$ and 
$\mathbb{P}(V^\mu(l)^*XV^\mu(l))=\mathbb{P}(U^\mu(l)^*XU^\mu(l))$ for any
$X \in \mathscr{M}\otimes\mathscr{W}_k$.  
\end{lemma}

From Lemma \ref{lem Bayes}, we immediately obtain the Kallianpur-Striebel 
formula,
\begin{equation}
\label{eq:KS-u}
	\pi_l^\mu(X)=\frac{\sigma_l^\mu(X)}{\sigma_l^\mu(I)}
		\quad\forall\,X\in\mathscr{M}\otimes\mathscr{C}_{l},
	\qquad
	\sigma_l^\mu(X)=U^\mu(l)^*\,\mathbb{P}(V^\mu(l)^*XV^\mu(l)|
		\mathscr{C}_l)\,U^\mu(l),
\end{equation}
and following the proof of Theorem \ref{thm:zakaihomodyne} gives the 
unnormalized controlled filtering equation for homodyne detection
\begin{equation*}
	\Delta\sigma_l^\mu(X)=
	\sigma_{l-1}^\mu(\mathcal{L}(X,\mathfrak{\check u}_l))\,\Delta t(l)+
	\sigma_{l-1}^\mu(\mathcal{J}(X,\mathfrak{\check u}_l))\,\Delta Y^\mu(l)
	,\qquad
	\sigma_0^\mu(X)=\rho(X).
\end{equation*}
Here we have used the controlled Lindblad generator
\begin{equation*}
	\mathcal{L}(X,u) =
		M^{+}(u)^*XM^+(u) +
		\lambda^2\,M^{\circ}(u)^*XM^\circ(u)
		+M^{\circ}(u)^*X+XM^\circ(u),
\end{equation*}
and we have written
$$
	\mathcal{J}(X,u)=XM^+(u) + M^{+}(u)^*X + \lambda^2
                M^{\circ}(u)^*XM^+(u) + 
		\lambda^2M^{+}(u)^* X M^\circ(u).
$$
As in section \ref{sec zakai homo}, we can also normalize this equation; 
this gives rise to the nonlinear controlled filtering equation for 
homodyne detection:
\begin{multline*}
	\Delta\pi_l^\mu(X)=
	\pi_{l-1}^\mu(\mathcal{L}(X,\mathfrak{\check u}_l))\,\Delta t(l) \, + \\
	\frac{\pi_{l-1}^\mu(\mathcal{J}(X,\mathfrak{\check u}_l))
	-\pi_{l-1}^\mu(X+\lambda^2\mathcal{L}(X,\mathfrak{\check u}_l))
		\,\pi_{l-1}^\mu(\mathcal{J}(I,\mathfrak{\check u}_l))}
	{I-\lambda^2\pi_{l-1}^\mu(\mathcal{J}(I,\mathfrak{\check u}_l))^2} 
	\times \\
	(\Delta Y^\mu(l)-\pi_{l-1}^\mu(\mathcal{J}(I,\mathfrak{\check u}_l))\,\Delta t(l)). 
\end{multline*}
Finally, we claim that even in the controlled case $\Delta
Y^\mu(l)-\pi_{l-1}^\mu(\mathcal{J}(I,\mathfrak{\check u}_l))\,\Delta t(l)$ 
is a martingale.  To show this, it suffices to demonstrate that
$\pi_{l-1}^\mu(\mathcal{J}(I,\mathfrak{\check u}_l))\,\Delta t(l)$ is the 
predictable part in the Doob decomposition of $\Delta Y^\mu(l)$: but this 
follows exactly as in the proof of the uncontrolled counterpart of this 
result,  Lemma \ref{lem pred Y homo}.

\subsection{The controlled quantum filter}
\label{sec:fb-controlled-filter}

Since $\pi_l^\mu(X)$ is linear in $X$, we can find a random $2\times 2$ 
matrix $\rho_l^\mu$ such that $\iota(\pi_l^\mu(X))={\rm Tr}[\rho_l^\mu X]$ 
for every $X\in\mathscr{M}$.  In fact, the conditional density matrix
$\rho^\mu_l$ satisfies the recursion
\begin{equation}\label{eq homo density form closed loop}
	\rho^\mu_l=\Gamma(\rho^\mu_{l-1}, u^\mu_l, \Delta y^\mu_l )
\end{equation}
where $y_l^\mu=\iota(Y^\mu(l))$,  $u_l^\mu=\iota(\mathfrak{u}_l)=
f_l(\Delta y_1^\mu,\ldots,\Delta y_{l-1}^\mu)$, 
	\nomenclature{$y^\mu_l$}{Classical observations process for 
		control strategy $\mu$}
	\nomenclature{$u^\mu_l$}{Classical feedback signal for
		control strategy $\mu$}
\begin{multline*}
	 \Gamma (\rho, u, \Delta y) =
	\rho+\mathcal{\overline{L}}(\rho,u)\,\lambda^2 +
	\frac{\mathcal{\overline{J}}(\rho,u)
	-{\rm Tr}[\mathcal{\overline{J}}(\rho,u)]\,
	(\rho+\lambda^2\mathcal{\overline{L}}(\rho,u))}
	{1-\lambda^2\,{\rm Tr}[\mathcal{\overline{J}}(\rho,u)]^2} \times\\
	(\Delta y -{\rm Tr}[\mathcal{\overline{J}}(\rho,u)]
		\,\lambda^2),
\end{multline*}
and 
\begin{equation*}
\begin{split}
	& \mathcal{\overline{L}}(\rho,u)=
	M^+(u)\rho M^{+}(u)^*+\lambda^2\,M^\circ(u)\rho M^{\circ}(u)^*+
	M^\circ(u) \rho+\rho M^{\circ}(u)^*, \\
	& \mathcal{\overline{J}}(\rho,u)=
	M^+(u)\rho+\rho M^{+}(u)^*+\lambda^2M^+(u)\rho M^{\circ}(u)^*+
	\lambda^2M^\circ(u)\rho M^{+}(u)^*.
\end{split}
\end{equation*}
To obtain the recursion (\ref{eq homo density form closed loop}) we can 
essentially follow the procedure used in section \ref{sec:howtouse} for 
the uncontrolled case.  The only subtlety here is that we need to deal 
with the presence of feedback in terms such as 
$\pi_{l-1}^\mu(XM^+(\check{\mathfrak{u}}_l))$ that occur in the recursion 
for $\pi_l^\mu(X)$.  The following lemma shows how to do this; it is
comparable to the classical statement $\mathbf{E}[ f(X,Y) \vert 
Y=y]=\mathbf{E}[f(X,y)\vert Y=y]$.

\begin{lemma}\label{lem:pulloutfeedback}
Consider a map $X:\mathfrak{U}\to\mathscr{M}$ and an observable
$\check{\mathfrak{u}}\in\mathscr{C}_{l}$ such 
that $u^\mu=\iota(U^\mu(l)^*\check{\mathfrak{u}}U^\mu(l))$ takes values 
in $\mathfrak{U}$.  Then $\iota(\pi_l^\mu[X(\check{\mathfrak{u}})])=
{\rm Tr}[\rho_l X(u^\mu))]$.
\end{lemma}

\begin{proof}
Using Definition \ref{def:Masafunctionofu}, we have
\begin{eqnarray*}
	\pi_l [ X(\check{\mathfrak{u}}) ] &=& 
	\pi_l \left[ \sum_{u\in\mathrm{sp}(\check{\mathfrak{u}})} 
	X(u)\, P_{\check{\mathfrak{u}}}(u)
	\right]
\\ &=&
	\sum_{u\in\mathrm{sp}(\check{\mathfrak{u}})}
	\mathbb{P}[U^\mu(l)^* X(u)\,P_{\check{\mathfrak{u}}}(u) U^\mu(l)
	|\mathscr{Y}_l]
\\ &=& 
	\sum_{u\in\mathrm{sp}(\check{\mathfrak{u}})}
	\mathbb{P}[U^\mu(l)^* X(u)U^\mu(l)|\mathscr{Y}_l]
	\,U^\mu(l)^*P_{\check{\mathfrak{u}}}(u)U^\mu(l)
\\ &=&
	\sum_{u\in\mathrm{sp}(\check{\mathfrak{u}})}
	\pi_l[X(u)]\,P_{U^\mu(l)^*\check{\mathfrak{u}}U^\mu(l)^*}(u).
\end{eqnarray*}
Applying $\iota$ to both sides, the result follows immediately.
\qquad
\end{proof}

It is sometimes more convenient to use the unnormalized form of the
filter.  As $\sigma_l^\mu(X)$ is linear in $X$, we can proceed exactly as 
before to find a random $2\times 2$ matrix $\varrho_l^\mu$ such that 
$\iota(\sigma_l^\mu(X))={\rm Tr}[\varrho_l^\mu X]$ for every 
$X\in\mathscr{M}$.  The result of Lemma \ref{lem:pulloutfeedback} is 
easily shown to hold also for $\sigma_l^\mu(\cdot)$, and we obtain
\begin{equation}
	\varrho_l^\mu= \Sigma( \varrho_{l-1}^\mu, u^\mu_l, \Delta y^\mu_l),
\label{eq homo density form closed loop unnormalized}
\end{equation}
where
\begin{equation}
	\Sigma( \varrho, u, \Delta y) = \rho+
	\mathcal{\overline{L}}(\varrho ,u)\,\lambda^2 +
	\mathcal{\overline{J}}(\nu,u)\,\Delta y .
\label{eq homo lambda def}
\end{equation}
The conditional density matrix $\rho_l^\mu$ can then be calculated as 
$\rho_l^\mu=\varrho_l^\mu/{\rm Tr}[\varrho_l^\mu]$.

\begin{remark} \label{rmk:markov-controlled-filter} 
The filters (\ref{eq homo density form closed loop}) and (\ref{eq homo 
density form closed loop unnormalized}) define classical controlled Markov 
processes \cite[Chapter 4]{KV86}, where future increments depend 
explicitly on previous control actions and on a commutative driving 
process $\Delta y^\mu_l$.  By analogy with the classical case, one could 
even consider the controlled quantum flow as a kind of ``quantum 
controlled Markov process'', though a precise statement of this concept is 
not yet used in the literature (but see also Remark \ref{rem:markov}).
\end{remark}

\subsection{Separated strategies}

\begin{figure}
\centering
\includegraphics[width=\textwidth]{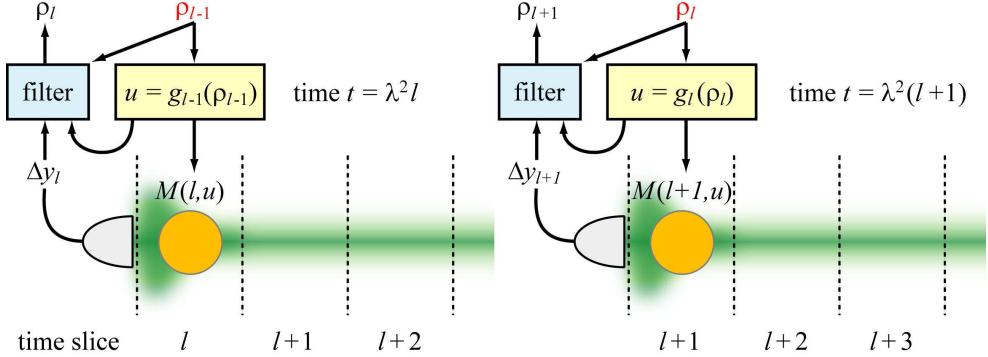}
\caption{\label{fig:control} Cartoon illustrating the operation of a 
separated feedback control.  The only quantity stored by the controller at 
the beginning of time step $l$ is the filter state in the previous time 
step $\rho_{l-1}$.  The control $u_l$ is calculated using a deterministic 
function $g_l$ and is fed back to the atom through an actuator.  
Simultaneously, time slice $l$ of the field interacts with the atom.
We subsequently detect the field, which gives rise to the observation 
$\Delta y_l$.  The filter uses $\Delta y_l$ and $u_l$ to update the 
filter state to $\rho_l$.  The procedure then repeats in time step 
$l+1$ with $\rho_l$ as basis for the control.} 
\end{figure}

Recall that an admissible strategy $\mu\in\mathfrak{K}$ is defined by a set
of feedback functions $\mu = \{ f_1,\ldots,f_{k} \}$ that describe how the
control values depend on the measurement record.  Any causal control
strategy can be written in this form.  We already remarked in section
\ref{sec:lookahead}, however, that in many cases the controls need not
have an arbitrarily complicated dependence on the measurement record---it
is sufficient to make the feedback at every time step a function of the
filter only.  In other words, many control goals can be attained using
only feedback of the form
\begin{equation}
\label{eq:sepdef}
	u_l^\mu=g_{l-1}(\rho_{l-1}^{\mu}),\qquad l=1,\ldots,k,
\end{equation}
where $g_l$ takes values in $\mathfrak{U}$ for any $l$.  Note that such a
strategy is in fact admissible as $\rho_{l-1}^\mu$ is a function of
$\Delta y_1^\mu,\ldots,\Delta y_{l-1}^\mu$ only, so that we could always
write $g_{l-1}(\rho_{l-1}^{\mu})=f_l(\Delta y_1^\mu,\ldots,\Delta
y_{l-1}^\mu)$ for some function $f_l$.

Both in the case of optimal control (section \ref{sec:optimal}) and in
control design using Lyapunov functions (section \ref{sec:lyapunov}), we
will see that the theory leads very naturally to strategies of the form
(\ref{eq:sepdef}).  This is highly desirable; not only does such structure 
significantly simplify the control design procedure, but the resulting 
controls are also much easier to implement in practice.  Note that to 
implement an arbitrary admissible strategy $\mu\in\mathfrak{K}$, the 
controller must have enough internal memory to store the entire 
observation history.  One has no other choice, as the feedback signal may 
be an arbitrary function of the control history.  On the other hand, to 
implement a strategy of the form (\ref{eq:sepdef}) the controller only 
needs to store the current density matrix $\rho_l^\mu$ at any time.
The density matrix is then updated recursively using the filter (\ref{eq 
homo density form closed loop}).

Control strategies of the form (\ref{eq:sepdef}) are called {\it separated
control strategies}, as the controller separates into two parts: a
filtering step, which can be solved recursively as described above, and a
control step which reduces to the simple evaluation of some deterministic
function of the filter state. The structure of the separated controller is
illustrated in Figure \ref{fig:control}.  The ubiquity of this separation
structure highlights the fundamental importance of filtering in (quantum)
stochastic control.

Let us finally fix the notation.  An admissible separated strategy $\mu$ 
is given by
	\nomenclature{$g_l(\rho)$}{Separated feedback function}
$$
	\mu = \{g_0(\rho),g_1(\rho),\ldots,g_{k-1}(\rho)\},
$$
where $g_i:\mathcal{S}\to\mathfrak{U}$ are functions on the set
$\mathcal{S}$ of $2\times 2$ density matrices (self-adjoint, nonnegative,
unit trace $2\times 2$ matrices).  We denote the set of all admissible 
separated strategies by $\mathfrak{K}_S$, and the corresponding feedback 
signal is given by (\ref{eq:sepdef}).
	\nomenclature{$\mathfrak{K}_S$}{Set of admissible separated controls}

Unlike in the absence of control (see section \ref{sec:filter-markov}), 
the controlled filter $\rho_l^\mu$ need not be a Markov process as in 
general $\mu \in \mathfrak{K}$ may have an arbitrary dependence on the 
observation history ($\rho_l^\mu$ is a controlled Markov process in the 
sense of \cite[Chapter 4]{KV86}).  On the other hand, for separated 
controls $\mu \in \mathfrak{K}_S$ the process $\rho_t^\mu$ is in fact 
Markov.  The following Lemma can be proved as in subsection 
\ref{sec:filter-markov}.
\begin{lemma}
For separated $\mu\in\mathfrak{K}_S$, the filter $\rho_l^\mu$ satisfies 
the Markov property:
$$
        \mathbf{E}^\mu(g(\rho_j^\mu)|\sigma\{\rho_0^\mu,\ldots,\rho_{l}^\mu\})
        = \mathbf{E}^\mu(g(\rho_j^\mu)|\sigma\{\rho_{l}^\mu\})\qquad
        \forall\,l\le j\le k.
$$
\end{lemma}

\section{Optimal control}
\label{sec:optimal}

Now that we have described repeated interactions with feedback and have
obtained filtering equations in this context, the remaining question is:  
how do we choose the control strategy $\mu$ to achieve a particular
control goal?  In this section we approach this problem using optimal
control theory, where the control goal is expressed as the desire to
minimize a certain cost function.  Dynamic programming allows us to 
construct a control strategy and to verify its optimality.  For an 
introduction to optimal stochastic control we refer to, e.g., 
\cite{KV86,Kus71,FR75,KA70,DB95}.

In optimal control theory, the control objective is expressed in terms of 
a cost function $J(\mu)$.  In this paper we consider mainly 
the cost function
        \nomenclature{$J(\mu)$}{Control cost function}
\begin{equation}
	J(\mu) = \mathbb{P}\left[ 
		\sum_{l=1}^{k} j^\mu_{l-1}(Q(\check{\mathfrak{u}}_{l}))  
		\,\Delta t(l)
		+ j^\mu_k(K) \right].
\label{eq:cost}
\end{equation}
This cost is defined in terms of the quantities $Q:\mathfrak{U}\to\mathscr{M}$
and $K\in\mathscr{M}$, with $K$ and $Q(u)$ nonnegative and self-adjoint 
for all $u\in\mathfrak{U}$.  Large values of these penalty quantities 
correspond to undesirable behavior, while small values correspond to the 
desirable behavior being sought. The penalty $Q(u)$ will be accumulated at 
each time step, contributing to a running cost.  For example, it may take 
the form $Q(u)=Q + c(u)$ where $Q\in\mathscr{M}$ penalizes deviation from 
a desired state, and the function $c:\mathfrak{U}\to[0,\infty)$ penalizes 
the control effort. $K$ is a terminal cost incurred at the end of 
the time horizon.  Ultimately the goal is to find, if possible, an 
admissible control strategy $\mu$ that minimizes the cost $J(\mu)$

\begin{remark}
Note that the cost $J(\mu)$ is defined on a fixed time interval of length 
$k$. We could also consider this problem on an infinite time horizon, and 
use either discounted or average cost per unit time extensions of
(\ref{eq:cost}) (without the terminal cost $K$).  Though the theory of 
this section extends readily to this scenario, we will restrict ourselves 
to the cost $J(\mu)$ for concreteness.
\end{remark}

The key step needed to apply dynamic programming is to express the cost 
function in terms of a filtered quantity. This will be discussed further 
in subsection \ref{sec:optimal-IS} where we describe the concept of 
information state. The following calculation, which makes use of 
properties of conditional expectations, the spectral theorem, and Lemma 
\ref{lem:pulloutfeedback}, expresses the cost (\ref{eq:cost}) in terms of 
the filter (\ref{eq homo density form closed loop}):
{\allowdisplaybreaks
\begin{eqnarray}
	J ( \mu ) & = &  \mathbb{P} \left[ \sum_{l=1}^{k} 
		\mathbb{P}[ j^\mu_{l-1}( Q( \check{\mathfrak{u}}_{l} )) \vert 
		\mathscr{Y}_{l-1} ]\,\Delta t(l)
		+  \mathbb{P}[ j^\mu_k( K ) \vert \mathscr{Y}_k ]
		 \right]
	\nonumber   \\ &=& 
		\mathbb{P} \left[ \sum_{l=1}^{k}  
		\pi^\mu_{l-1}[ Q( \check{\mathfrak{u}}_{l} ) ]\,\Delta t(l)
		+  \pi^\mu_k[   K ]  \right]
	\nonumber   \\ &=&  
		\mathbf{E}^\mu \left[ \sum_{l=1}^{k}  
		\iota(\pi^\mu_{l-1}[ Q( \check{\mathfrak{u}}_{l} ) ])
		\,\Delta t
		+ \iota( \pi^\mu_k[   K ] )
		\right]
  	\nonumber \\ &=&  
		\mathbf{E}^\mu \left[ 
		\sum_{l=1}^{k}  
		\mathrm{Tr}[\rho^\mu_{l-1} Q( u^\mu_{l} ) ]
		\,\Delta t
		+ \mathrm{Tr}[ \rho^\mu_k   K ] 
		\right].
\label{eq:cost-2}
\end{eqnarray}
}
Here $\mathbf{E}^\mu$ is the $\mathbf{P}^\mu$-expectation on the 
probability space $(\Omega^\mu,\mathcal{F}^\mu,\mathbf{P}^\mu)$ obtained 
by applying the spectral theorem to $\mathscr{Y}_k^\mu$.  We retain the 
superscript $^\mu$ to emphasize that the observation algebra is 
control-dependent; $Y^\mu(l)$ and $Y^\nu(l)$ need not commute!

The final form of the cost is now a classical cost function for the
classical recursion (\ref{eq homo density form closed loop}); in other
words, we have reduced the control problem to a problem of classical 
optimal control for the classical controlled Markov process (\ref{eq homo 
density form closed loop}).  Hence we can apply almost directly the 
ordinary dynamic programming methodology, which we will do in the next 
subsection.  We only need to take care to use different probability spaces 
$\Omega^{\mu}$, $\Omega^\nu$ for different control strategies, as the 
corresponding observations need not commute; this does not present any 
additional difficulties, however.

In order to implement the dynamic programming method, we will use the 
fact that the conditional distribution of $\Delta y^\mu_l$ given $\Delta 
y_1^\mu,\ldots, \Delta y_{l-1}^\mu$ can be evaluated explicitly.  This 
proceeds in a manner similar to the proof of Lemma \ref{lemma:Delta-y}, 
but taking into account the dependence of the coefficients on the controls 
using Lemma \ref{lem:pulloutfeedback}.

\begin{lemma}  \label{lemma:Delta-y-fb}
For any admissible strategy $\mu \in \mathfrak{K}$ we have
\begin{equation}
	\mathbf{P}^\mu[ \Delta y_l = \pm \lambda \, \vert \, 
	\mathcal{Y}^\mu_{l-1} ]
	= p(\Delta y = \pm \lambda; \rho^\mu_{l-1}, u^\mu_l ) ,
\label{eq:Delta-y-fb-1}
\end{equation}
where
\begin{equation}
 	p(\Delta y = \pm \lambda; \rho, u ) =
 	\frac{1}{2}\pm\frac{\lambda}{2}\,
		{\rm Tr}[\mathcal{\overline{J}}(\rho,u)   ]
\label{eq:Delta-y-fb-2}
\end{equation}
depends only on the previous filter state and the control value being 
applied, and in particular is independent of the feedback strategy $\mu$.
\end{lemma}

Before we proceed, let us fix some notation.  In the following we will 
encounter separated control strategies $\mu,\mu^*,\bar\mu\in\mathfrak{K}_S$.
We will denote the corresponding control functions by $\mu=\{g_0(\rho),
\ldots,g_{k-1}(\rho)\}$, $\mu^*=\{g_0^*(\rho),\ldots,g_{k-1}^*(\rho)\}$, 
and $\bar\mu=\{\bar g_0(\rho),\ldots,\bar g_{k-1}(\rho)\}$.  We will also 
denote by $\mu\in\mathfrak{K}$ an arbitrary admissible strategy, in which 
case the corresponding feedback process is always denoted by $u_l^\mu$.

\subsection{Dynamic Programming}
\label{sec:optimal-DP}

The goal of this section is to give a brief introduction to dynamic
programming.  The main results we will need are the dynamic programming
equation or Bellman equation (\ref{eq:bellman}), which allows us to
construct explicitly a candidate optimal control strategy in separated
form, and the verification Lemma \ref{lem:separation}, which verifies that
the strategy constructed through the dynamic programming equation is
indeed optimal.  As motivation for the Bellman equation we will first
prove the ``converse'' of the verification lemma: if an optimal separated
strategy exists, then the Bellman equation follows.  The Bellman equation
can also be introduced in a more general context without this assumption,
and without any reference to separated strategies; this will be discussed
briefly at the end of this section.

Let us begin by considering a separated control strategy 
$\mu\in\mathfrak{K}_S$.  A central idea in dynamic programming is that one 
should consider separately the cost incurred over different time intervals 
within the total control interval $0,\ldots,k$.  To this end, let us 
introduce the {\it cost-to-go}
        \nomenclature{$W_l(\mu,\rho)$}{Cost-to-go}
$$
	W_l(\mu,\rho)=
	\mathbf{E}^\mu \left[\left.
		\sum_{i={l+1}}^{k}  
		\mathrm{Tr}[\rho^\mu_{i-1} Q(u^\mu_{i})]\,\Delta t
		+ \mathrm{Tr}[\rho^\mu_k K] 
	\right|
	\rho_l^\mu=\rho\right].
$$
The quantity $W_l(\mu,\rho)$ is the cost incurred by the control $\mu$ 
over the interval $l,\ldots,k$, given that $\rho_l^\mu=\rho$.  Note in 
particular that $W_0(\mu,\rho_0)=J(\mu)$.  It is straightforward to obtain 
a recursion for $W_l(\mu,\rho)$ using the Markov property of the 
controlled filter:
\begin{equation}
\label{eq:costtogo}
\begin{split}
	W_l(\mu,\rho) &=
	\mathbf{E}^\mu[\mathrm{Tr}[\rho^\mu_l Q(u^\mu_{l+1})]\,\Delta t
	  + W_{l+1}(\mu,\rho_{l+1}^\mu)|
	\rho_l^\mu=\rho]
\\ &= 
	\mathrm{Tr}[\rho Q(g_l(\rho))]\,\Delta t + 
	\sum_{\Delta y =\pm\lambda}
	p(\Delta y; \rho,g_l(\rho))\,
	W_{l+1}(\mu,\Gamma(\rho,g_l(\rho),\Delta y)),
\end{split}
\end{equation}
with the terminal condition $W_k(\mu,\rho)={\rm Tr}[\rho K]$.  

The goal of optimal control theory is to find, if possible, a control 
strategy $\mu^*$ that minimizes $J(\mu)$.  For the time being, let us 
suppose that such a strategy exists within the class of separated 
controls, i.e.\ that $\mu^*\in\mathfrak{K}_S$ is such that $J(\mu^*)\le 
J(\mu)$ for any $\mu\in\mathfrak{K}_S$.  It is not unnatural to expect 
that the same strategy $\mu^*$ also minimizes $W_l(\mu,\rho)$ for any $l$ 
and $\rho$: if the strategy $\mu^*$ is optimal over the full time interval 
$0,\ldots,k$, then it should also be optimal over any subinterval 
$l,\ldots,k$.  If we assume that this is indeed the case, then we can 
simplify the expression for $W_l(\mu^*,\rho)$:
\begin{equation*}
\begin{split}
	W_l&(\mu^*,\rho)  =  \inf_{\mu\in\mathfrak{K}_S}W_l(\mu,\rho) \\
	& =  \inf_{g_l,\ldots,g_k}\left[
		\mathrm{Tr}[\rho Q(g_l(\rho))]\,\Delta t + 
		\sum_{\Delta y =\pm\lambda}
		p(\Delta y; \rho,g_l(\rho))\,
		W_{l+1}(\mu,\Gamma(\rho,g_l(\rho),\Delta y))
	\right] \\
	& =  \inf_{g_l}\left[
		\mathrm{Tr}[\rho Q(g_l(\rho))]\,\Delta t + 
		\inf_{g_{l+1},\ldots,g_k}
		\sum_{\Delta y =\pm\lambda}
		p(\Delta y; \rho,g_l(\rho))\,
		W_{l+1}(\mu,\Gamma(\rho,g_l(\rho),\Delta y))
	\right] \\
	& =  \inf_{g_l}\left[
		\mathrm{Tr}[\rho Q(g_l(\rho))]\,\Delta t + 
		\sum_{\Delta y =\pm\lambda}
		p(\Delta y; \rho,g_l(\rho))\,
		\inf_{\mu\in\mathfrak{K}_S}
		W_{l+1}(\mu,\Gamma(\rho,g_l(\rho),\Delta y))
	\right] \\
	& =  \inf_{g_l}\left[
		\mathrm{Tr}[\rho Q(g_l(\rho))]\,\Delta t + 
		\sum_{\Delta y =\pm\lambda}
		p(\Delta y; \rho,g_l(\rho))\,
		W_{l+1}(\mu^*,\Gamma(\rho,g_l(\rho),\Delta y))
	\right].
\end{split}
\end{equation*}
The following Lemma makes these ideas precise.

\begin{lemma} {\rm ({\bf Dynamic programming equation})}.
\label{lem:bellman}
Suppose $\mu^*\in\mathfrak{K}_S$ is optimal in $\mathfrak{K}_S$, i.e.\   
$J(\mu^*)\le J(\mu)$ for any $\mu\in\mathfrak{K}_S$, and define
$V_l(\rho)=W_l(\mu^*,\rho)$.  Then
        \nomenclature{$V_l(\rho)$}{Value function of dynamic programming}
\begin{eqnarray}
\nonumber
	V_l(\rho)
	&=& \inf_{u\in\mathfrak{U}}\left[
	\mathrm{Tr}[\rho Q(u)]\,\Delta t + 
	\sum_{\Delta y =\pm\lambda}
	p(\Delta y; \rho,u)\,
	V_{l+1}(\Gamma(\rho,u,\Delta y))
	\right],\quad l<k,
	\\
	V_k(\rho) &=& {\rm Tr}[\rho K],
\label{eq:bellman}
\end{eqnarray}
for all $l$ and any $\rho$ that is reachable under $\mu^*$ in the 
sense that $\mathbf{P}^{\mu^*}(\rho_l^{\mu^*}=\rho)>0$.  Moreover, the 
infimum in {\rm (\ref{eq:bellman})} is attained at $u=g^{*}_l(\rho)$.
\end{lemma}

\begin{proof}
In view of the discussion above, it suffices to show that the conditions 
of the Lemma imply that $\mu^*$ minimizes $W_l(\mu,\rho)$, i.e., that 
$W_l(\mu^*,\rho)\le W_l(\mu,\rho)$ for any $\mu\in\mathfrak{K}_S$, 
$l=0,\ldots,k$, and any reachable $\rho$.  To this end, let us suppose
that this statement holds for time step $l+1$.  Then by 
(\ref{eq:costtogo})
$$
	\mathrm{Tr}[\rho Q(g_l(\rho))]\,\Delta t +
	\sum_{\Delta y =\pm\lambda}
	p(\Delta y; \rho,g_l(\rho))\,
	W_{l+1}(\mu^*,\Gamma(\rho,g_l(\rho),\Delta y))
	\le W_l(\mu,\rho)
$$
for any $\mu\in\mathfrak{K}_S$ and reachable $\rho$.  Note that the 
left hand side is precisely $W_l(\mu',\rho)$, where $\mu'$ is the control 
strategy that coincides with $\mu$ at time $l$ and with $\mu^*$ at all 
other times (this follows as $W_{l+1}(\mu^*,\rho)$ only depends on the 
control functions at times $l+1,\ldots,k$).  We would like to show that 
the left hand side is bounded from below by 
$$
	W_l(\mu^*,\rho) = 
	\mathrm{Tr}[\rho Q(g_l^{*}(\rho))]\,\Delta t +
	\sum_{\Delta y =\pm\lambda}
	p(\Delta y; \rho,g_l^{*}(\rho))\,
	W_{l+1}(\mu^*,\Gamma(\rho,g_l^{*}(\rho),\Delta y))
$$
for all reachable $\rho$.  Suppose that this is not the case.  We define a 
new control strategy $\bar\mu\in\mathfrak{K}_S$ as follows.  At time $l$, 
we set
$$
	\bar g_l(\rho)=\left\{\begin{array}{ll}
		g_l^{*}(\rho) &\qquad   W_l(\mu^*,\rho) \le W_l(\mu',\rho),
			\\
		g_l(\rho) &\qquad \mbox{all other $\rho$}.
	\end{array}
	\right.
$$
For any time $i\ne l$, we let $\bar\mu$ coincide with $\mu^*$, i.e.\
$\bar g_i(\rho)=g_i^*(\rho)$ for $i\ne l$.  Clearly
\begin{multline*}
	W_l(\bar\mu,\rho) = 
	\mathrm{Tr}[\rho Q(\bar g_l(\rho))]\,\Delta t + \\
	\sum_{\Delta y=\pm\lambda}
	p(\Delta y; \rho,\bar g_l(\rho))\,
	W_{l+1}(\mu^*,\Gamma(\rho,\bar g_l(\rho),\Delta y)) \le
	W_l(\mu^*,\rho)
\end{multline*}
with strict inequality for some reachable $\rho$.  But then
$$
	J(\mu^*)=
	\mathbf{E}^{\mu^*} \left[
		\sum_{i=1}^{l}  
		\mathrm{Tr}[\rho^{\mu^*}_{i-1} 
		Q(g^{*}_{i-1}(\rho_{i-1}))]\,\Delta t
		+ 
		W_l(\mu^*,\rho^{\mu^*}_l)
	\right]\ge J(\bar\mu),
$$
which contradicts optimality of $\mu^*$.  Hence evidently 
$W_l(\mu^*,\rho)\le W_l(\mu,\rho)$ also for time step $l$.  It remains to 
notice that the statement holds trivially for $l=k$ as $W_k(\mu,\rho)={\rm 
Tr}[\rho K]$ is independent of $\mu$, and hence for any $l=0,\ldots,k$
by induction.

Returning to the statement of the Lemma, (\ref{eq:bellman}) can now be 
obtained as above, and the fact that the infimum is attained at 
$g^{*}_l(\rho)$ follows directly if we compare (\ref{eq:bellman}) with the 
expression for $W_l(\mu^*,\rho)$ given by (\ref{eq:costtogo}).  Hence the 
proof is complete.
\qquad
\end{proof}

Let us now consider (\ref{eq:bellman}) {\it without} assuming that an
optimal control $\mu^*$ exists.  As a backwards in time recursion, this
expression makes sense without any reference to $\mu^*$.  $V_l(\rho)$ is
uniquely defined for any $l=0,\ldots,k$ and $\rho$: after all, the cost
quantities $\mathrm{Tr}[\rho Q(u)]$ and ${\rm Tr}[\rho K]$ are bounded
from below, so that the infimum exists and is uniquely defined at every
step in the recursion starting from the terminal time $k$.

Eq.\ (\ref{eq:bellman}) is called the {\bf dynamic programming equation}
or {\bf Bellman equation}.  By Lemma \ref{lem:bellman}, we
find that if an optimal control $\mu^*\in\mathfrak{K}_S$ exists, then it
must be the case that $V_l(\rho)=W_l(\mu^*,\rho)$ and that all the infima
in the recursion are attained; moreover, in this case an optimal control
can be recovered by choosing the control functions $g_l(\rho)$ to be 
minimizers at every step in the Bellman recursion.

One could now ask whether the converse also holds true.  Suppose that in 
the Bellman recursion all the infima turn out to be attained.  Then we can 
{\it define} a separated control strategy $\mu^*$ by setting
\begin{equation}
	g^{*}_l(\rho) \in \mathop{\mathrm{argmin}}_{u\in\mathfrak{U}}  
	\left[
		\mathrm{Tr}[\rho Q(u)]\,\Delta t + 
		\sum_{\Delta y =\pm\lambda}
		p(\Delta y; \rho,u)\,
		V_{l+1}(\Gamma(\rho,u,\Delta y))
	\right].
\label{eq:dynprog}
\end{equation}
The question is, does this imply that $\mu^*$ is optimal?  This is indeed 
the case, as we will show in the following Lemma.  Note that this 
immediately gives a constructive way of finding optimal controls, which is 
precisely the main idea of {\it dynamic programming}.

\begin{lemma} {\rm ({\bf Verification and separation})}.
\label{lem:separation}
	Suppose that all the infima in the Bellman equation are
	attained.  Then any control $\mu^*$ defined by {\rm 
	(\ref{eq:dynprog})} is optimal in $\mathfrak{K}$, i.e.\ 
	$J(\mu^*)\le J(\mu)$ for any $\mu\in\mathfrak{K}$.  Moreover
	$V_l(\rho)=W_l(\mu^*,\rho)$ for any $l$ and all reachable $\rho$, 
	and in particular $V_0(\rho_0)=J(\mu^*)$.
\end{lemma}

\begin{proof}
By substituting (\ref{eq:dynprog}) in (\ref{eq:bellman}) and comparing 
with (\ref{eq:costtogo}), the last statement is evident.  It remains to 
show that $J(\mu^*)\le J(\mu)$ for any $\mu\in\mathfrak{K}$.  Note that
$$
	V_l(\rho_l^\mu)\le
	\mathrm{Tr}[\rho_l^\mu Q(u_{l+1}^\mu)]\,\Delta t + 
	\sum_{\Delta y =\pm\lambda}
	p(\Delta y; \rho_l^\mu,u_{l+1}^\mu)\,
	V_{l+1}(\Gamma(\rho_l^\mu,u_{l+1}^\mu,\Delta y))	
$$
by (\ref{eq:bellman}), where we have chosen an arbitrary 
$\mu\in\mathfrak{K}$.  But by Lemma \ref{lemma:Delta-y-fb}, this is
$$
	V_l(\rho_l^\mu)\le
	\mathrm{Tr}[\rho_l^\mu Q(u_{l+1}^\mu)]\,\Delta t + 
	\mathbf{E}^\mu(V_{l+1}(\rho_{l+1}^\mu)|\mathcal{Y}_l^\mu).
$$
By recursing the relation backwards from $l=k-1$, we obtain
$$
	V_0(\rho_0)\le\mathbf{E}^\mu\left[
		\sum_{l=1}^{k}  
		\mathrm{Tr}[\rho^\mu_{l-1} Q( u^\mu_{l} ) ]\,\Delta t
		+ \mathrm{Tr}[ \rho^\mu_k   K ] 
	\right] = J(\mu).
$$
But we had already established that $V_0(\rho_0)=J(\mu^*)$.  Hence
the result follows. \qquad
\end{proof}

Let us reflect for a moment on what we have achieved.  We began by showing
that the cost-to-go for any strategy $\mu^*$ that is optimal in the class
$\mathfrak{K}_S$ of separated controls must satisfy the Bellman recursion.  
Conversely, if the infima in the Bellman equation are all attained, we can
{\it construct} a control strategy $\mu^*$ by solving the Bellman
recursion ({\it dynamic programming}).  We then verified in Lemma
\ref{lem:separation} (the {\it verification lemma}) that the strategy thus
constructed is indeed optimal, not only in the class of separated controls
but in the class $\mathfrak{K}$ of all admissible controls.  Combining
these two results, we conclude that any strategy that is optimal in
$\mathfrak{K}_S$ is necessarily optimal within the larger class
$\mathfrak{K}$ (Lemma \ref{lem:separation} is also called the {\it
separation lemma} for this reason).  This shows that the idea of separated
controls is very natural, and indeed universal, for this type of control
problem.

\begin{remark}
We have not given conditions under which existence of the infima in the 
Bellman recursion is guaranteed.  This can be a delicate issue, see e.g.\
\cite{BS78}.  However, if $\mathfrak{U}$ is a finite set the infima are 
trivially attained.  This particularly simple case is also the most 
straightforward to implement on a computer. 
\end{remark}

\begin{remark}
The results above are formulated only for reachable $\rho$.  This is as it 
should be, as for non-reachable $\rho$ the cost-to-go is defined as a 
conditional expectation on a set of measure zero and is thus not unique.  
These issues are irrelevant, however, to the application of the dynamic 
programming algorithm and the verification lemma, which are the most 
important results of this section.  
\end{remark}

Finally, we briefly remark on the case where the infima in the Bellman
equation are not attained.  It follows from Lemma \ref{lem:bellman}
that there cannot exist an optimal separated control in this case.  
Nonetheless the solution $V_l(\rho)$ of the Bellman equation, called the
{\it value function}, is still a relevant quantity; it can be shown to be
the infimum (not minimum!)\ of the cost-to-go over a suitable class of
(not necessarily separated) controls.  This characterization can be 
useful in practice, for example, if we wish to quantify the performance of 
sub-optimal controls.  

This approach also provides a different entry point into the theory of
dynamic programming than the one we have chosen, providing additional
insight into the structure of the theory.  One could begin by proving
directly that the value function $V_l(\rho)$, now {\it defined} as the
infimum of the cost over the time horizon $l,\ldots,k$, satisfies the
dynamic programming equation.  This does not require us to assume the
existence of an optimal control, and in particular places no a priori
preference on the class of separated controls.  The only thing that is
needed is the fact that the filtering equation is a controlled Markov
process, which is key to the entire procedure.  Indeed, we have used this
property in an essential way in the form of Lemma \ref{lemma:Delta-y-fb}.
If the infima in the Bellman recursion are attained then we can construct
an explicit separated control strategy as was demonstrated above;
verification then proceeds as we have indicated.

We will not further detail this approach here, as we do not need these
results in the following.  We refer to \cite{BS78} for an extensive study
of discrete time dynamic programming, or to e.g.\ \cite{FR75} for the
continuous time case.

\subsection{Information States}
\label{sec:optimal-IS}

The key idea that facilitated the dynamic programming solution to the
optimal control problem discussed above was the representation of the cost
$J(\mu)$ in terms of the filtered quantity $\rho^\mu_l$. This is an 
instance of a general methodology involving information state 
representations \cite[Chapter 6]{KV86}. An {\em information state} is a 
quantity that is causally computable from information available to the 
controller, i.e.\ the observation record and previous choices of control 
values. The dynamic evolution of the information state is given by an {\em 
information state filter}. Solving a given optimal control problem reduces 
to expressing the cost function exclusively in terms of a suitable
information state, then applying the dynamic programming method. The
resulting optimal feedback control will be a separated control relative to
the information state filter.  In the previous sections we used the 
density operator $\rho^\mu_l$ as an information state to represent the 
cost $J(\mu)$, as indicated in the calculation (\ref{eq:cost-2}), and the 
filter for this information state was given by (\ref{eq homo density form 
closed loop}).

The choice of information state for a given problem is not unique. In the 
case of the cost function $J(\mu)$ defined by (\ref{eq:cost}), we could 
also use the unnormalized conditional density operator $\varrho^\mu_l$ 
discussed in subsection \ref{sec:fb-controlled-filter} as an information 
state, with the corresponding filter (\ref{eq homo density form closed 
loop unnormalized}). To see this, we use the reference probability method 
and condition on $\mathscr{C}_l$ as follows. We begin by defining the state 
$$
	\mathbb{P}^{0\mu}[X] = \mathbb{P}[ U^\mu(k) X U^\mu(k)^\ast ], 
	\ \ X \in \mathscr{M} \otimes \mathscr{W}_k,
$$
and we denote the associated classical state, as obtained through the 
spectral theorem, by $\mathbf{P}^{0\mu}$.
Under $\mathbf{P}^{0\mu}$, $\Delta y_1, \ldots, \Delta y_k$ are i.i.d.\ 
random variables taking values $\pm\lambda$ with equal probability, and in 
particular the law of the process $\Delta y$ under $\mathbf{P}^{0\mu}$ is 
independent of the feedback control $\mu \in \mathfrak{K}$. Then (cf.\
(\ref{eq:cost-2})) we have
\begin{eqnarray*}
	J ( \mu ) & = &  \mathbb{P} \left[ \sum_{l=1}^{k}  
			V^\mu(l-1)^\ast Q( \check{\mathfrak{u}}_{l} 
			)V^\mu(l-1)\,\Delta t(l)
			  +   V^\mu(k)^\ast K V^\mu(k)  ]
			 \right]
		\nonumber  \\
			&=& \mathbb{P} \left[ \sum_{l=1}^{k} 
			\mathbb{P}[  V^\mu(l-1)^\ast Q( 
			\check{\mathfrak{u}}_{l} )V^\mu(l-1) \vert 
			\mathscr{C}_{l-1} ]  \,\Delta t(l)
			+  \mathbb{P}[ V^\mu(k)^\ast K 
			V^\mu(k)\vert \mathscr{C}_k ]
			 \right]
		\nonumber   \\
			&=& \mathbb{P} \left[ \sum_{l=1}^{k} U^\mu(k) 
			U^\mu(k)^\ast \mathbb{P}[  V^\mu(l-1)^\ast Q( 
			\check{\mathfrak{u}}_{l} )V^\mu(l-1) \vert 
			\mathscr{C}_{l-1} ]U^\mu(k) U^\mu(k)^\ast  
			\,\Delta t(l)
			\right.
		\nonumber \\
			&& \hspace{1.3cm} \left.
			\phantom{\sum_{l=1}^{k}}
			 +  U^\mu(k) 
			U^\mu(k)^\ast \mathbb{P}[ V^\mu(k)^\ast K 
			V^\mu(k) \vert \mathscr{C}_k ] U^\mu(k) 
			U^\mu(k)^\ast
			 \right],
		\nonumber   
\end{eqnarray*}
where the change of state operator $V^\mu(l)$ was defined in Lemma 
\ref{lem nondemo u}.  Now recall that $U^\mu(k)^\ast 
\mathbb{P}[X|\mathscr{C}_l ]U^\mu(k) =
U^\mu(l)^\ast \mathbb{P}[X|\mathscr{C}_l ]U^\mu(l)$, see the proof of the 
nondemolition property in section \ref{sec:repeatedinteraction}.
Changing to the state $\mathbb{P}^{0\mu}$, we therefore obtain
\begin{eqnarray*}
	 J(\mu) &=&  \mathbb{P}^{0\mu} \left[ \sum_{l=1}^{k}  
			\sigma^\mu_{l-1}[ Q( \check{\mathfrak{u}}_{l} ) 
			]\,\Delta t(l)
			+  \sigma^\mu_k[   K ]
			  \right]
		 \nonumber   \\
			  &=&  \mathbf{E}^{0\mu} \left[ \sum_{l=1}^{k}  
			\iota(\sigma^\mu_{l-1}[ Q( \check{\mathfrak{u}}_{l} ) ])
			\,\Delta t
			+ \iota( \sigma^\mu_k[   K ] )
			  \right]
		  \nonumber \\
			  &=&  \mathbf{E}^{0\mu} \left[ \sum_{l=1}^{k}  
			\mathrm{Tr}[\varrho^\mu_{l-1} Q( u^\mu_{l} ) ]
			\,\Delta t
			+ \mathrm{Tr}[ \varrho^\mu_k   K ] 
			  \right],
\end{eqnarray*}
where the unnormalized conditional state $\sigma^\mu_l$ was defined by 
(\ref{eq:KS-u}), and $\varrho^\mu_l$ is the associated density matrix.
Using this representation, we can define a value function $S_l(\varrho)$ 
and find an optimal control using the alternate dynamic programming 
equation
\begin{eqnarray*}
	S_l(\varrho) &=&  \inf_{u \in \mathfrak{U}} \left[
		\mathrm{Tr}[\varrho Q(u) ]\,\Delta t +
		 \sum_{\Delta y = \pm \lambda} 
		 p^0(\Delta y )\,
		 S_{l+1}(\Sigma(\varrho,u,\Delta y))
 	\right], \quad l<k,
 \nonumber \\
	S_k(\varrho) &=& \mathrm{Tr}[ \varrho K ],
\end{eqnarray*}
where $p^0(\Delta y = \pm \lambda )=0.5$. In fact, the optimal control 
is given by
$$
	h^*_l(\varrho) \in \mathop{\mathrm{argmin}}_{u \in \mathfrak{U}}  
	\left[\mathrm{Tr}[\varrho Q(u) ]\,\Delta t+
	\sum_{\Delta y = \pm \lambda} p^0(\Delta y  ) \, 
	S_{l+1}(\Sigma(\varrho,u,\Delta y )) 
	\right]  .
$$
This is a separated feedback control relative to the information state
filter (\ref{eq homo density form closed loop unnormalized}).

The conditional state (either normalized or unnormalized) is not the 
correct choice for every cost function. In 1981 Whittle \cite{W81} 
discovered a different type of information state for optimal control 
problems with {\it exponential} cost functions (risk-sensitive control), 
which are not solvable using the standard conditional states. Instead, the 
filter for the corresponding information state depends explicitly on 
quantities defining the cost function, and the optimal feedback control
is separated relative to this filter. Following \cite{Jam05} we now 
explain this briefly in the context of this paper.

In the quantum setting, a risk-sensitive cost can be defined by
$$
	J^\theta ( \mu ) = \mathbb{P} \left[ 
	C(k)^\ast e^{\theta j^\mu_k(K)} C(k)
	\right] ,
$$
where $\Delta C(l) = \frac{\theta}{2}j^\mu_{l-1}(Q(\check{\mathfrak{u}}_{l-1}))
\,C(l-1)\,\Delta t(l)$, $C(0)=I$ defines the ``exponential'' running cost
and $\theta > 0$ is a fixed real parameter (the risk parameter).   Let 
us now define ${\tilde U}^\mu(l)= U^\mu(l)C(l)$ (which is not unitary in 
general), so that
$$
	J^\theta ( \mu ) = \mathbb{P} \left[ 
	\tilde U^\mu(k)^\ast e^{\theta K } \tilde U^\mu(k)
	\right].
$$
We can now proceed as in the previous part of this section to express the 
control cost in terms of an unnormalized filtered quantity.  The 
corresponding filter is not obtained from the usual unitary $U^\mu(l)$, 
however, but from the modified operator $\tilde U^\mu(l)$.  We can obtain
a change-of-state operator $\tilde V^\mu(l)$ as in Lemma \ref{lem nondemo 
u}, which gives rise to an information state filter that depends 
explicitly on the running cost $Q(u)$.  We can subsequently express
$J^\theta(\mu)$ in terms of this filter, and the optimal control problem 
can then be solved using dynamic programming.  We leave the details as an 
exercise.

In classical stochastic control, the risk-sensitive cost is known to 
possess improved robustness properties compared to the usual cost 
$J(\mu)$; in particular, as the risk parameter $\theta$ increases, the 
optimal performance becomes less sensitive to the details of the 
underlying model.  To what extent these advantages carry over to the 
quantum case remains to be explored (but see \cite{JP05}).

\subsection{Example}
\label{sec:optimal-eg}

We will give a numerical example of dynamic programming for a particularly 
simple system---a controlled version of the dispersive interaction model.

\subsubsection*{The controlled quantum flow}

We consider again the dispersive interaction model of section 
\ref{sec:examplesrepint}, but now we add a control input.  The controlled 
repeated interaction matrices are now given by
$$
	L_1(u) = 0,\qquad
	L_2(u) = i\sigma_z,\qquad
	L_3(u) = iu(\sigma_+-\sigma_-).
$$
Such an interaction can be realized in certain systems by applying a 
magnetic field of strength $u$ to the atom, see e.g.\ Figure 
\ref{fig:photo}, in addition to the usual dispersive interaction with the 
electromagnetic field.  In principle any $u\in\mathbb{R}$ is admissible, 
i.e.\ we should take $\mathfrak{U}=\mathbb{R}$.  As we will be evaluating 
the dynamic programming recursion numerically, however, it is more 
convenient to discretize the admissible controls, i.e.\ we choose 
$\mathfrak{U}=\{-K\varepsilon,-(K-1)\varepsilon,\ldots,K\varepsilon\}\subset
\mathbb{R}$ where $\varepsilon$ is the discretization step size and $2K+1$ 
is the total number of admissible controls.

Our first task is to evaluate the matrices $M^{\pm,+,-,\circ}(u)$ in the 
controlled difference equation.  Let us calculate explicitly
$$
	M_l(u) = \exp(
		\sigma_z(\Delta A^*(l)-\Delta A(l))+
		u(\sigma_+-\sigma_-)\Delta t(l)
	).
$$
Writing
$$
	B(u)=\sigma_z\otimes\lambda\,(\sigma_+-\sigma_-)
		+u(\sigma_+-\sigma_-)\otimes\lambda^2I
	=
	\begin{pmatrix}
		0 & \lambda^2u & \lambda & 0 \\
		-\lambda^2u & 0 & 0 & -\lambda \\
		-\lambda & 0 & 0 & \lambda^2u \\
		0 & \lambda & -\lambda^2u & 0
	\end{pmatrix},
$$
we calculate the matrix exponential
$$
	e^{B(u)}=
	\begin{pmatrix}
		\cos(\lambda w(u)) & u\lambda\,\frac{\sin(\lambda w(u))}{w(u)}
			& \frac{\sin(\lambda w(u))}{w(u)} & 0 \\
		-u\lambda\,\frac{\sin(\lambda w(u))}{w(u)} & \cos(\lambda w(u))
			& 0 & -\frac{\sin(\lambda w(u))}{w(u)} \\
		-\frac{\sin(\lambda w(u))}{w(u)} & 0 
			& \cos(\lambda w(u)) & u\lambda\,
				\frac{\sin(\lambda w(u))}{w(u)} \\
		0 & \frac{\sin(\lambda w(u))}{w(u)}
			& -u\lambda\,
		\frac{\sin(\lambda w(u))}{w(u)} & \cos(\lambda w(u))
	\end{pmatrix},
$$
where we have written $w(u)=\sqrt{1+\lambda^2u^2}$.  Hence we obtain
\begin{multline*}
	M_l(u)=\frac{\sin(\lambda w(u))}{\lambda w(u)}\,\sigma_z\,
	(\Delta A^*(l)-\Delta A(l)) \\ 
	+ \left[
		\frac{\cos(\lambda w(u))-1}{\lambda^2}
		+\frac{\sin(\lambda w(u))}{\lambda w(u)}\,u(\sigma_+-\sigma_-)
	\right]\Delta t(l)+I.
\end{multline*}
We can now immediately read off the coefficients in the quantum stochastic 
difference equation for the controlled dispersive interaction model:
\begin{equation*}
\begin{array}{ll}
	M^\pm(u)=0, &\qquad
	M^+(u)=\frac{\sin(\lambda w(u))}{\lambda w(u)}\,\sigma_z, \\
	~~~ & ~~~\\
	M^-(u)=-\frac{\sin(\lambda w(u))}{\lambda w(u)}\,\sigma_z,
	&\qquad
	M^\circ(u) = \frac{\cos(\lambda w(u))-1}{\lambda^2}
                +\frac{\sin(\lambda w(u))}{\lambda w(u)}\,
	u\,(\sigma_+-\sigma_-).
\end{array}
\end{equation*}

\subsubsection*{An invariant set}

The solution $\rho_l^\mu$ of the filtering recursion is always a density
matrix: recall that $\iota(\pi_l^\mu(X))={\rm Tr}[\rho_l^\mu X]$ is a 
positive map with $\pi_l^\mu(I)=1$, so that $\rho_l^\mu$ must be a 
positive matrix and have unit trace.  The the goal of dynamic programming
is to calculate the feedback function $g^*_l(\rho)$ for any time step $l$
and density matrix $\rho$.  Unfortunately, even the space of $2\times 2$
density matrices is rather large; discretization of this space, as one
would need for computer implementation of the dynamic programming
recursion, would require a tremendous number of discretization points.  
For this reason dynamic programming is computationally expensive,
prohibitively so in moderate- to high-dimensional systems where optimal
control theory often plays the role of a benchmark rather than a practical
solution.  In such cases, one is compelled to settle for control designs
that are sub-optimal, i.e.\ they do not minimize a cost function.  We will 
briefly explore one such approach in section \ref{sec:lyapunov}.

The simple example treated in this section, however, has a feature 
which significantly simplifies the implementation of dynamic programming.
As we will demonstrate shortly, there is an invariant set of density 
matrices which is parametrized by points on a circle: i.e., if we start 
the filter somewhere on this circle, it will always remain there.  
This reduces the dynamic programming algorithm to the calculation of the 
feedback function $g^*_l(\rho)$ on the circle.  With a sufficiently 
accurate discretization of the circle, this problem can be solved numerically 
without too many complications.

To find the desired invariant set, note that
\begin{equation}\label{eq:invsetcalc}
\begin{split}
	&
	\mathcal{L}^\mu(\sigma_z,u)=
	\frac{\sin(2\lambda w(u))}{\lambda w(u)}
	\,u\,(\sigma_++\sigma_-) - \frac{2u^2\sin^2(\lambda w(u))}{w(u)^2}\,
	\sigma_z, \\
	&
	\mathcal{J}^\mu(\sigma_z,u)=\frac{\sin(2\lambda w(u))}{\lambda w(u)}
	\,I, \\
	&
	\mathcal{L}^\mu(\sigma_++\sigma_-,u)=
	-\frac{\sin(2\lambda w(u))}{\lambda w(u)}
	\,u\,\sigma_z 
	- \frac{2\sin^2(\lambda w(u))}{\lambda^2}\,
	(\sigma_++\sigma_-), \\
	&
	\mathcal{J}^\mu(\sigma_++\sigma_-,u)=
	-\frac{2u\sin^2(\lambda w(u))}{w(u)^2}
	\,I, \\
	&
	\mathcal{J}^\mu(I,u)=
	\frac{\sin(2\lambda w(u))}{\lambda w(u)}\,\sigma_z+
	\frac{2u\sin^2(\lambda w(u))}{w(u)^2}\,(\sigma_++\sigma_-).
\end{split}
\end{equation}
Hence evidently the recursion for $\pi_l^\mu(\sigma_z)$ and 
$\pi_l^\mu(\sigma_++\sigma_-)$ forms a closed set of equations.  We claim 
furthermore that $\pi_l^\mu(\sigma_z)^2+\pi_l^\mu(\sigma_++\sigma_-)^2=I$ 
for all $l$, if this is true for $l=0$.  The algebra is a little easier if 
we consider the unnormalized version; using the discrete It\^o rule, we 
calculate
\begin{multline*}
\Delta(\sigma^\mu_l(X)^2)= \\
\left\{
	2\sigma^\mu_{l-1}(X)\,
	\sigma_{l-1}^\mu(\mathcal{L}^\mu(X,u_l^\mu))
	+\lambda^2\sigma_{l-1}^\mu(\mathcal{L}^\mu(X,u_l^\mu))^2
	+\sigma_{l-1}^\mu(\mathcal{J}^\mu(X,u_l^\mu))^2
\right\}\Delta t(l) \\
+
	2\sigma^\mu_{l-1}(X+\lambda^2\mathcal{L}^\mu(X,u_l^\mu))\,
	\sigma_{l-1}^\mu(\mathcal{J}^\mu(X,u_l^\mu))\,
\Delta Y^\mu(l).
\end{multline*}
A tedious but straightforward calculation shows that
$$
	\Delta[\sigma_l^\mu(\sigma_z)^2+\sigma_l^\mu(\sigma_++\sigma_-)^2
		-\sigma_l^\mu(I)^2]=
	\vartheta_{u,\lambda}\,
	[\sigma_{l-1}^\mu(\sigma_z)^2+\sigma_{l-1}^\mu(\sigma_++\sigma_-)^2
                -\sigma_{l-1}^\mu(I)^2]\,\Delta t(l)
$$
where $\vartheta_{u,\lambda}$ is a complicated function of $u$ and 
$\lambda$.  Hence if 
$\pi_0^\mu(\sigma_z)^2+\pi_0^\mu(\sigma_++\sigma_-)^2=I$, it
follows that $\pi_l^\mu(\sigma_z)^2+\pi_l^\mu(\sigma_++\sigma_-)^2=I$
for all $l$.

Now let $\tilde\rho$ be any density matrix.  The remaining insight we need 
is that if we are given $x,z\in\mathbb{R}$ such that ${\rm 
Tr}[\tilde\rho(\sigma_++\sigma_-)]=x$, ${\rm Tr}[\tilde\rho\sigma_z]=z$,
and $x^2+z^2=1$, then this uniquely determines $\tilde\rho$.  To see this, 
let us write without loss of generality $x=\sin\theta$ and $z=\cos\theta$.
Using the constraints ${\rm Tr}[\tilde\rho(\sigma_++\sigma_-)]=\sin\theta$, 
${\rm Tr}[\tilde\rho\sigma_z]=\cos\theta$, ${\rm Tr}\rho=1$ and 
$\rho=\rho^*$, we easily find that
$$
	\tilde\rho=
	\begin{pmatrix}
		\frac{1}{2}+\frac{1}{2}\cos\theta & \frac{1}{2}\sin\theta 
			+i\beta \\
		 \frac{1}{2}\sin\theta-i\beta & \frac{1}{2}
			-\frac{1}{2}\cos\theta
	\end{pmatrix}\quad\mbox{for some}\quad \beta\in\mathbb{R}.
$$
But we can explicitly calculate the eigenvalues of this matrix as
$\tfrac{1}{2}(1\pm (1+4\beta^2)^{1/2})$, so that the remaining requirement 
for the density matrix $\rho\ge 0$ implies that $\beta=0$.  Hence we
conclude that the ``circle'' of density matrices
$$
	S^1=
	\left\{
	\begin{pmatrix}
		\frac{1}{2}+\frac{1}{2}\cos\theta & \frac{1}{2}\sin\theta 
			\\
		 \frac{1}{2}\sin\theta & \frac{1}{2}
			-\frac{1}{2}\cos\theta
	\end{pmatrix}
	:
	\theta\in [0,2\pi)
	\right\},
$$
parametrized by the angle $\theta$, is invariant under the filtering 
equation for our controlled quantum flow in the sense that $\rho_l^\mu\in 
S^1$ for all $l$ if $\rho_0\in S^1$.  We can thus restrict the dynamic 
programming recursion to this set, which yields a feedback control law on 
the circle of the form $g_l(\theta)$.

\subsubsection*{Dynamic programming}

We are now in a position to solve the dynamic programming algorithm 
numerically.  To this end, we have discretized the circle into a set of 
$10^5$ equidistant points, and we have discretized the control set 
$\mathfrak{U}$ into $400$ equidistant points in the interval $[-10,10]$.  
As in the previous simulations, we have chosen $\lambda^{-2}=300$ and a 
terminal time of $3$ (i.e.\ $k=900$).

\begin{figure}
\centering
\includegraphics[width=\textwidth]{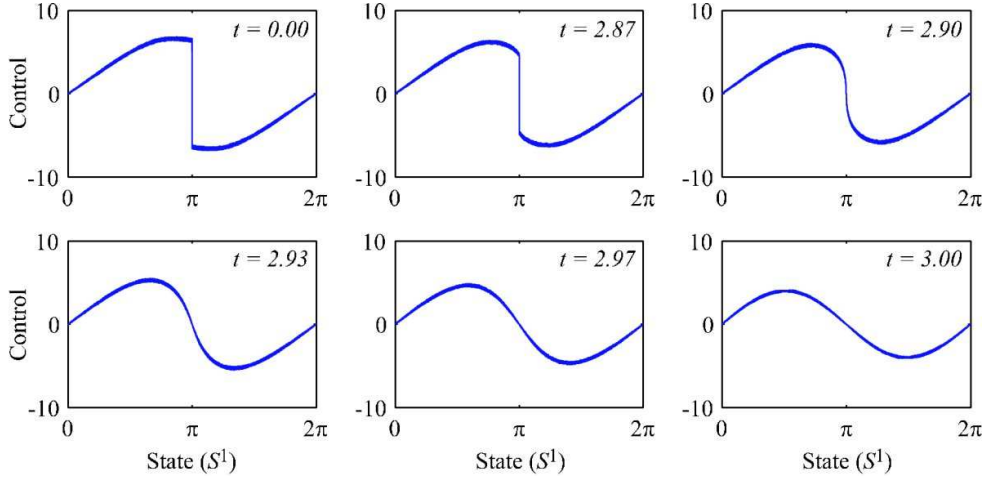}
\caption{\label{fig:dynprogorig} 
The optimal control function $g^*_l(\theta)$, for selected $l$, obtained 
by applying the dynamic programming algorithm to the cost function Eq.\ 
(\ref{eq:democost1}).  During the major portion of the control interval, 
the feedback function is essentially equal to the first plot; the control 
passes through zero at the state of maximal energy, and has a break at the 
point of minimal energy.  Close to the terminal time, the optimal control 
varies a little to accommodate the terminal cost.}
\end{figure}

As a first control goal, suppose we would like to maximize the expected 
energy, i.e.\ we would like to drive ${\rm Tr}[\rho_l^\mu\sigma_z]$ to 
$+1$.  To this end, we  use the cost (\ref{eq:cost}) with
\begin{eqnarray}\label{eq:democost1}
Q(u) &=& C u^2 + D (I-\sigma_z) ,
\nonumber \\	
K &=& 	I-\sigma_z
\end{eqnarray}
The first term in $Q(u)$ penalizes the use of large control strengths, which is 
necessary in any practical feedback loop.  The second term tries to 
minimize $1-{\rm Tr}[\rho_l^\mu\sigma_z]$ during the running time of the 
system, whereas the terminal cost tries to minimize the terminal value of 
$1-{\rm Tr}[\rho_l^\mu\sigma_z]$.  The constants $C,D>0$ determine the 
relative weights attributed to these control goals.  As an example, we 
have chosen $C=0.25$, $D=5$.  The corresponding optimal feedback function 
$g_l(\theta)$ is plotted for several times $l$ in Figure 
\ref{fig:dynprogorig}.  Note that the control passes through zero at the 
state of maximal energy $\theta=0$, whereas the break at $\theta=\pi$ 
drives the system away from the state of minimal energy.  This is not 
unexpected, as both $\theta=\{0,\pi\}$ are fixed points of the filter 
(see section \ref{sec:examplesfilter}); the break in the control 
destabilizes the undesired minimal energy state.

\begin{figure}
\centering
\includegraphics[width=\textwidth]{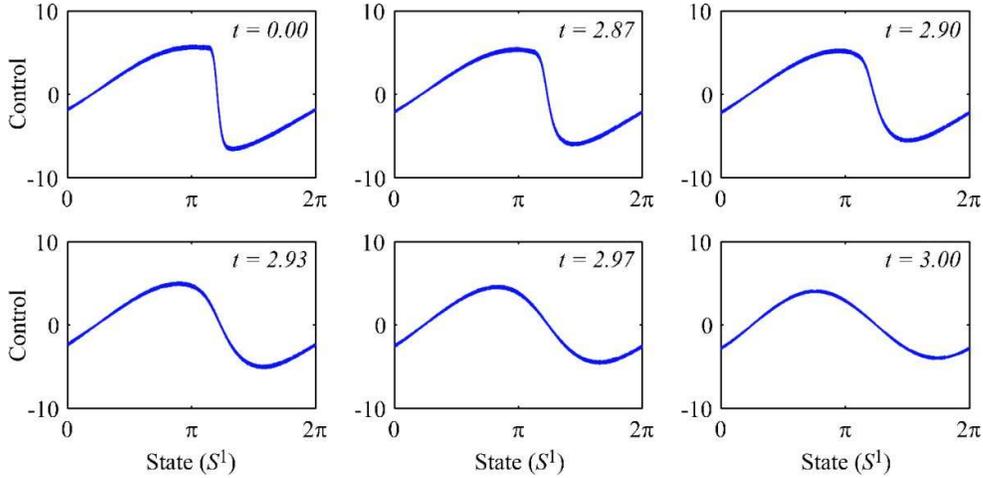}
\caption{\label{fig:dynprog45} 
The optimal control function $g^*_l(\theta)$ for the cost Eq.\ 
(\ref{eq:democost2}).  As before, the feedback function is essentially 
equal to the first plot during the major portion of the control interval, 
and the control passes through zero at the target state.  In this case, 
however, the control function is not discontinuous at the point of the 
circle opposite to the target state.}
\end{figure}

As a second example, suppose we would like to drive the system close to 
the point where ${\rm Tr}[\rho_l^\mu\sigma_z]={\rm 
Tr}[\rho_l^\mu\sigma_x]=2^{-1/2}$, i.e.\ $\theta=\pi/4$.  This point is 
uniquely characterized by ${\rm Tr}[\rho_l^\mu(\sigma_z+
\sigma_x)]=\sqrt{2}$, so we choose 
\begin{eqnarray}\label{eq:democost2}
Q(u) &=& C u^2 + D (I-X) ,
\nonumber \\	
K &=& 	I-X
\end{eqnarray}
where $X=2^{-1/2}(\sigma_x+\sigma_z)$ and, for example, $C=0.25$ and 
$D=5$.  The optimal feedback function for this case is plotted in Figure 
\ref{fig:dynprog45}.  Once again, the feedback function passes through 
zero at the target point $\theta=\pi/4$.  However, note that the function 
is no longer singular; as the opposite point on the circle is not a fixed 
point of the filter in this case, it is evidently more efficient to let 
the filter drive itself toward the target point without expending 
additional control effort.

\begin{figure}
\centering
\includegraphics[width=\textwidth]{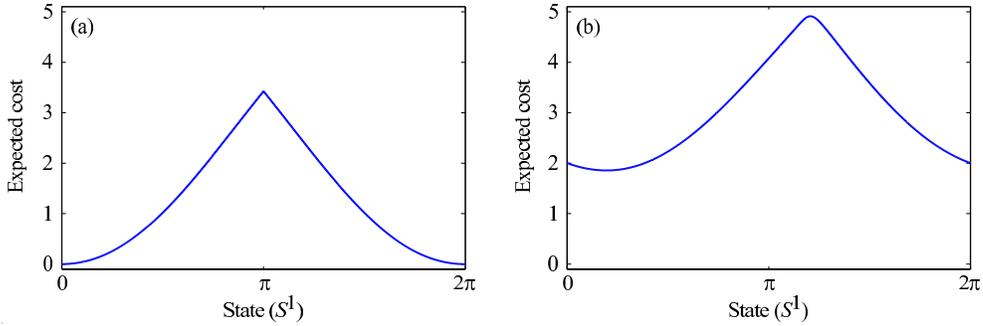}
\caption{\label{fig:dynprogvalue} 
The cost $J(\mu^*)$, as a function of $\rho_0\in S^1$, for the optimal 
control strategy $\mu^*$ minimizing the cost {\rm (a)} Eq.\ 
(\ref{eq:democost1}), and {\rm (b)} Eq.\ (\ref{eq:democost2}).  The cost 
is simply obtained as the final value function in the dynamic programming 
recursion, $J(\mu^*)(\rho_0)=V_0[\rho_0]$.  In case {\rm (a)} the cost is 
zero at $\theta=0$: after all, the target point is a fixed point of the 
filter, so no cost is accumulated if we start at the target point.  This 
is not the case in {\rm (b)} where the target point is not a fixed point.}
\end{figure}

Finally we have plotted the optimal cost $J(\mu^*)$, as a function of the 
initial state $\rho\in S^1$, in Figure \ref{fig:dynprogvalue}.  The cost 
is easily obtained as a byproduct of dynamic programming, as it is simply 
given by the value function at time zero.  Note that in the case of our 
first example, zero cost is possible: as the target $\theta=0$ is a fixed 
point of the filter, no cost is accumulated if the system is initially in 
its state of maximal energy.  In our second example, on the other hand, 
this is not the case: even if we start at the target state $\theta=\pi/4$, 
the filter will fluctuate around this point and a total cost (just under 
$J[\mu^*]=2$) is accumulated over the control interval.

\section{Lyapunov control}
\label{sec:lyapunov}

Many control design methods other than optimal control have been developed
and are applied in control engineering, for a variety of reasons. Indeed,
the applicability of dynamic programming is limited by computational
complexity for use in low dimensional problems, or in situations where
explicit solutions are available. Among the alternative methods are {\em
Lyapunov} design methods, and the purpose of this section is to
demonstrate their use in our quantum context.

 Consider for example the following scenario.  
Recall that the dispersively interacting atom, when untouched by the 
experimenter, has constant but unknown energy; our control goal is to 
drive the energy to a particular value of our choice, say 
$+\hbar\omega_0/2$.  In the absence of sharp time constraints this problem 
is not uniquely expressed as an optimal control problem: what cost should 
one choose?  All we want to achieve is that 
$\mathbb{P}(j_l(H))\to\hbar\omega_0/2$ for large $l$, a much more modest
control goal than the absolute minimization of a particular cost function.  
Hence dynamic programming is unnecessarily complicated for the problem at 
hand, and we can resort to a simpler method for control design.
In this section we will design a controller for the control problem
described above through a very simple Lyapunov function method, similar in
spirit to \cite{HSM05}.  The simplicity of the method is its chief
feature: not only will we find a feedback function that works, but the 
resulting feedback function is also of a very simple form and is easily 
implementable.

As before, we can easily express the control goal in terms of the filter:
$\mathbb{P}(j_l(H))=\mathbb{P}(\pi_l(H))\to\hbar\omega_0/2$ for large $l$.
We will design a separated control that achieves this goal.  The Markov 
property of the filter with separated controls is key to the procedure: it 
allows us to treat the problem entirely at the level of the filter, 
without reference to the original repeated interaction model.  

\subsection{Lyapunov stability}
\label{sec:lyapunov-stab}

The main 
tool is the following Lyapunov theorem for Markov processes, taken 
directly from \cite{Kus71}.

\begin{theorem}[\bf Lyapunov theorem]
Let $x_l$ be a (classical) Markov process.  Suppose there exists a 
nonnegative (Lyapunov) function $V(x)\ge 0$ that satisfies
$$
	\mathbf{E}(\Delta V(x_l)|\sigma\{x_{l-1}\})=
	\mathbf{E}(V(x_l)|\sigma\{x_{l-1}\}) 
		- V(x_{l-1}) = -k(x_{l-1})\le 0
	\qquad\forall\,l,
$$
where $k(x)\ge 0$ is another nonnegative function.  Then $k(x_l)\to 0$ as 
$l\to\infty$ a.s.
\end{theorem}

\begin{proof}
Fix some $n>0$.  As $V(x)$ is a nonnegative function, clearly
$$
	V(x_0)\ge V(x_0)-\mathbf{E}(V(x_n)).
$$
But by the condition of the Theorem, we obtain
$$
	\mathbf{E}(V(x_n))
	=\mathbf{E}(V(x_{n-1})-k(x_{n-1}))
	=\mathbf{E}(V(x_{n-2})-k(x_{n-2})-k(x_{n-1}))=\ldots
$$
Iterating this procedure, we obtain
$$
	V(x_0)\ge \mathbf{E}\left(\sum_{l=0}^{n-1} k(x_l)
	\right).
$$
As this holds for any $n$ and as $V(x_0)<\infty$, the right hand side
is finite for any $n$.  But
$$
	\sum_{l=0}^{\infty} \mathbf{P}(k(x_l)\ge\varepsilon)\le
	\frac{1}{\varepsilon}\,\mathbf{E}\left(\sum_{l=0}^{\infty} k(x_l)
        \right)<\infty\qquad
	\forall\,\varepsilon>0
$$
by Chebyshev's inequality, so the Borel-Cantelli Lemma gives $k(x_l)\to 
0$ a.s.
\qquad
\end{proof}

\begin{remark}
Previously we did everything on a fixed time horizon $l=0,\ldots,k$.  
Now, however, we are considering what happens as $l\to\infty$, so
technically a little more care is needed.  The reader should convince
himself that not much changes in this case.  In particular, there is no
need to deal directly with an infinite-dimensional algebra 
$\mathscr{W}_\infty=\mathscr{M}^{\ten\infty}$ corresponding to an 
infinite number of time slices of the electromagnetic field.  Rather, one 
can use the natural embedding $i:\mathscr{M}\otimes\mathscr{W}_k\to
\mathscr{M}\otimes\mathscr{W}_{k+1}$, $i(X)=X\otimes I$ and the sequence
of probability spaces $\Omega^k\subset\Omega^{k+1}$ obtained by applying
the spectral theorem to the corresponding observation algebras, to give 
meaning to infinite time limits without functional analytic complications 
(provided the time step $\lambda$ and repeated interaction model is 
fixed).  We leave the details as an exercise.
\end{remark}

Before we design the controller, let us verify that without control it is
indeed the case that $\pi_l(H)\to\pm\hbar\omega_0/2$ a.s. This already 
demonstrates the Lyapunov theorem.

\begin{lemma}
If $u_l^\mu=0$, then ${\rm Tr}[\rho_l^\mu\sigma_z]\to\pm 1$ with unit 
probability.
\end{lemma}

\begin{proof}
We will use the Lyapunov function $V(\rho_l^\mu)=1-{\rm 
Tr}[\rho_l^\mu\sigma_z]^2$, which is nonnegative and is zero precisely 
when ${\rm Tr}[\rho_l^\mu\sigma_z]=\pm 1$.  Using the filtering equation 
with $u_l^\mu=0$ and Eq.\ (\ref{eq:invsetcalc}), we find that
$$
	\Delta\, {\rm Tr}[\rho_l^\mu\sigma_z]=
	\frac{\sin(2\lambda)}{\lambda} \frac{V(\rho_{l-1}^\mu)}{
		\cos^2(2\lambda) + \sin^2(2\lambda) V(\rho_{l-1}^\mu)
	}\Delta\tilde y^\mu_l,
$$
where $\tilde y^\mu_l$ is the innovations process.  Using the quantum 
It\^o rule
$$
	\Delta V(\rho_{l}^\mu)=
	-\frac{\sin^2(2\lambda)}{\lambda^2} \frac{V(\rho_{l-1}^\mu)^2}{
		(\cos^2(2\lambda) + \sin^2(2\lambda) V(\rho_{l-1}^\mu))^2
	}(\Delta\tilde y^\mu_l)^2
	+(\cdots)\,\Delta\tilde y^\mu_l.
$$
But $\mathbf{E}^\mu(\Delta\tilde y^\mu_l|\mathscr{Y}^\mu_{l-1})=0$ by the 
martingale property, and furthermore
$$
	\mathbf{E}^\mu((\Delta\tilde y^\mu_l)^2|\mathscr{Y}^\mu_{l-1})=
	(1-\lambda^2\,{\rm Tr}[
		\mathcal{\overline{J}}(\rho^\mu_{l-1})]^2)\,\Delta t \\
	=
	(\cos^2(2\lambda) + \sin^2(2\lambda) V(\rho_{l-1}^\mu))\,\Delta t,
$$
where we have used $(\Delta\tilde Y(l))^2=(\cdots)\,\Delta\tilde Y(l)
+(\Delta Y(l))^2-(\Delta C(l))^2$ as in the proof of Lemma \ref{lem xi 
homo}.  Hence we find
$$
	\mathbf{E}^\mu(\Delta V(\rho_{l}^\mu)|\sigma\{\rho^\mu_{l-1}\})=
	-\frac{\sin^2(2\lambda)}{\lambda^2} \frac{V(\rho_{l-1}^\mu)^2}{
		\cos^2(2\lambda) + \sin^2(2\lambda) V(\rho_{l-1}^\mu)
	}\Delta t\le 0.
$$
The Lemma now follows from the Lyapunov theorem.
\qquad
\end{proof}

\subsection{Construction of a Lyapunov control}
\label{sec:lyapunov-control}

We now turn to the control design.  We wish to find a (time-invariant) feedback control 
$u_l^\mu=g(\rho_{l-1}^\mu)$ so that $\mathbb{P}(\pi_l^\mu(\sigma_z))\to+1$ 
as $l\to\infty$.  The way we approach this problem is to use a trial 
Lyapunov function $V(\rho)$ without fixing the control $g$.  By inspecting 
the expression for $\mathbf{E}^\mu(\Delta V(\rho_l^\mu)|
\sigma\{\rho_{l-1}^\mu\})$, we can subsequently try to choose $g$ so that 
this expression is nonpositive and is zero only at (or in a small 
neighborhood of) the point $\tilde\rho\in S^1$ where ${\rm 
Tr}[\tilde\rho\sigma_z]=+1$.  The desired result follows by dominated 
convergence.

Let us implement this procedure.  Choose the trial Lyapunov function 
$V(\tilde\rho)=1-{\rm Tr}[\tilde\rho\sigma_z]$.   Using the filtering 
equation with $u_l^\mu=g(\rho_{l-1}^\mu)$, Eq.\ (\ref{eq:invsetcalc}), and 
the martingale property of the innovations, we find that
\begin{multline*}
	\mathbf{E}^\mu(\Delta V(\rho_l^\mu)|\sigma\{\rho_{l-1}^\mu\})
	=
	-\frac{\sin(2\lambda w(f(\rho_{l-1}^\mu)))}
		{\lambda w(f(\rho_{l-1}^\mu))}
	\,f(\rho_{l-1}^\mu)\,{\rm Tr}[\rho_{l-1}^\mu(\sigma_++\sigma_-)]\,\Delta t
\\
	+ \frac{2f(\rho_{l-1}^\mu)^2\sin^2(\lambda w(f(\rho_{l-1}^\mu)))}
		{w(f(\rho_{l-1}^\mu))^2}\,
	{\rm Tr}[\rho_{l-1}^\mu\sigma_z]\,\Delta t.
\end{multline*}
As a first attempt, consider the following feedback function:
$$
	f(\tilde\rho)=\left\{\begin{array}{ll}
		-1 & \qquad \mbox{if } 
			{\rm Tr}[\tilde\rho(\sigma_++\sigma_-)] < 0, \\
		+1 & \qquad \mbox{if }
			{\rm Tr}[\tilde\rho(\sigma_++\sigma_-)] \ge 0.
	\end{array}\right.
$$
Then
$$
	\frac{\mathbf{E}^\mu(\Delta V(\rho_l^\mu)|\sigma\{\rho_{l-1}^\mu\})
	}{\Delta t}
	=
	-\frac{\sin(2\lambda w(1))}{\lambda w(1)}
	\,|{\rm Tr}[\rho_{l-1}^\mu(\sigma_++\sigma_-)]|
	+\frac{2\sin^2(\lambda w(1))}{w(1)^2}\,
	{\rm Tr}[\rho_{l-1}^\mu\sigma_z].
$$
Keeping in mind that $\rho_l^\mu\in S^1$ for all $l$, clearly there exists 
some $\delta>0$ such that this expression is strictly negative for ${\rm 
Tr}[\rho_{l-1}^\mu\sigma_z]<1-\delta$; in fact, $\delta$ is the solution 
of
$$
	-\frac{\sin(2\lambda w(1))}{\lambda w(1)}
	\,\sqrt{(2-\delta)\delta}
	+\frac{2\sin^2(\lambda w(1))}{w(1)^2}\,
	(1-\delta)=0.
$$
To keep $\mathbf{E}^\mu(\Delta V(\rho_l^\mu)|\sigma\{\rho_{l-1}^\mu\})$
nonpositive everywhere, we now modify the control function to turn off the 
feedback in the set ${\rm Tr}[\rho_{l-1}^\mu\sigma_z]\ge 1-\delta$:
\begin{equation}\label{eq:lyapfeedback}
	g(\tilde\rho)=\left\{\begin{array}{ll}
		0  & \qquad \mbox{if } 
			{\rm Tr}[\tilde\rho\sigma_z] \ge 1-\delta, \\
		-1 & \qquad \mbox{if } 
			{\rm Tr}[\tilde\rho(\sigma_++\sigma_-)] < 0
		\mbox{ and } {\rm Tr}[\tilde\rho\sigma_z] < 1-\delta, \\
		+1 & \qquad \mbox{if }
			{\rm Tr}[\tilde\rho(\sigma_++\sigma_-)] \ge 0
		\mbox{ and } {\rm Tr}[\tilde\rho\sigma_z] < 1-\delta.
	\end{array}\right.
\end{equation}
This ensures that $\mathbf{E}^\mu(\Delta V(\rho_l^\mu)|
\sigma\{\rho_{l-1}^\mu\})\le 0$, where the equality holds only when ${\rm 
Tr}[\rho_{l-1}^\mu\sigma_z]\ge 1-\delta$.  The Lyapunov theorem then 
guarantees that ${\rm Tr}[\rho_{l}^\mu\sigma_z]$ converges to the 
set $[1-\delta,1]$: i.e.\ $\liminf_{l\to\infty}{\rm 
Tr}[\rho_{l}^\mu\sigma_z]\ge 1-\delta$ with unit probability.  For 
$\lambda\ll 1$, the threshold $\delta$ will be very small and hence the 
Lyapunov theorem guarantees convergence to a tiny neighborhood of the 
desired control goal.  For example, with the choice $\lambda^{-2}=300$ 
which we have used in the simulations, $\delta\approx 5.6\times 10^{-6}$.

\begin{remark}
We can do better than proving convergence to a small neighborhood of the 
target point.  In addition to the convergence $k(x_t)\to 0$, the invariant 
set theorems \cite{Kus71} tell us that $x_t$ itself converges to the 
largest invariant set contained in $\{x:k(x)=0\}$ (the convergence is in 
probability, but this is easily strengthened to a.s.\ convergence).  As 
the only invariant set inside the set ${\rm Tr}[\tilde\rho\sigma_z] \ge 
1-\delta$ is the target point ${\rm Tr}[\tilde\rho\sigma_z]=1$, 
convergence is guaranteed.  A full discussion of the required theorems is 
beyond the scope of this article, and we refer to \cite{Kus71} for further 
details.
\end{remark}

\begin{figure}
\centering
\includegraphics[width=\textwidth]{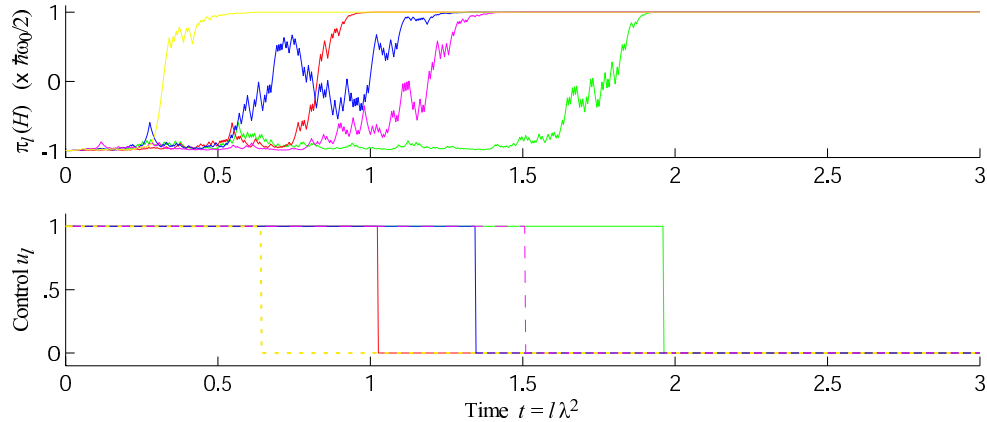}
\caption{\label{fig:stabil} 
Five typical sample paths of the conditional expectation $\pi_l(H)$ of the 
atomic energy (top plot) for the controlled quantum flow on the circle 
with separated feedback law Eq.\ (\ref{eq:lyapfeedback}).  Note that
$\pi_l(H)=+\hbar\omega_0/2$ is always attained.  The bottom plot 
shows the feedback signal $u_l^\mu=g(\rho_{l-1}^\mu)$.  The initial state 
was $\rho(x)=\langle\Phi,X\Phi\rangle$, $\rho\in S^1$, the time 
scale is $\lambda^{-2}=300$, and $\delta=6\times 10^{-6}$.} 
\end{figure}

To illustrate the effectiveness of Eq.\ (\ref{eq:lyapfeedback}), we have 
used the Monte Carlo method to simulate in Figure \ref{fig:stabil} several 
sample paths of the controlled filter.  It is immediately evident that the 
control goal was attained for the five sample paths plotted, and the 
Lyapunov theorem guarantees that this is indeed the case for any sample 
path.

The discussion above is mainly intended as an illustration of an
alternative to optimal control theory, and we have only used the Lyapunov 
theory in its simplest form.  Lyapunov methods can be extended to a large 
class of systems, see e.g.\ \cite{mazyar}, whose treatment using optimal 
control would be out of the question due to the high dimensionality of the 
state space.  Some control goals other than strict convergence can be 
treated using similar methods; for example, we could try to find a control 
law so that a quantity such as $\mathbf{E}(|{\rm Tr}[\rho_l^\mu\sigma_z]-
\alpha|)$, for some $\alpha\in(-1,1)$, becomes small (but usually nonzero)\
as $l\to\infty$, see e.g.\ \cite{naoki}.  Beside the restricted set of 
control goals that can be treated by Lyapunov methods, a drawback of 
such methods is that there is no general recipe which one can follow for 
the design of controls; the choice of a proper control and Lyapunov 
function has to be investigated on a case-by-case basis.  Ultimately the 
control goal of interest and the available resources determine which 
method of control design---be it optimal control, Lyapunov methods, or 
some other approach---is most suitable for the problem at hand.


\section{References for further reading}\label{continuous time}

We have come to the end of our exposition on the filtering and control of 
discrete quantum models, though in many ways we have only scratched the 
surface of the theory of quantum probability, filtering and stochastic 
control.  The goal of this final section is to provide some entry points 
into the literature for further reading.

\subsection*{Quantum probability}

We have used quantum probability theory in its simplest finite-dimensional
form, where only finite state random variables are available.  In practice
such a description is very restrictive, and one needs a theory that admits
continuous random variables.  The theory of operator algebras (e.g.\
Kadison and Ringrose \cite{KaR83,KaR86} or Bratteli and Robinson 
\cite{BrR87}), particularly Von Neumann algebras, provides the
proper generalization of the concept of a quantum probability spaces to
the infinite-dimensional context.  For an accessible introduction to 
quantum probability we refer to the excellent lecture notes by Maassen 
\cite{Maa03}.  See also \cite{Bia95,BHJ06} and the references therein.  
The required background on functional analysis can be found e.g.\ in the 
textbook by Reed and Simon \cite{ReS80}.

Readers with a background in physics might wonder why we did not use the 
projection postulate in this paper.  We do not need this postulate to 
develop the theory, and in fact the projection postulate can be obtained 
as a special case within the quantum probability framework used in this 
paper.  See \cite{BHJ06} for an example.

\subsection*{Quantum noise and stochastic calculus}

The three discrete noises---we have suggestively denoted them as $A(l)$, 
$A^*(l)$ and $\Lambda(l)$---are replaced in the continuous time theory by 
the standard noises $A_t$ (the annihilation process), $A_t^*$ (the 
creation process) and $\Lambda_t$ (the gauge process).  For each time $t$, 
$A_t,A^*_t,\Lambda_t$ are defined as operators acting on a Hilbert space 
which is known as the Fock space, and the latter plays a central role in 
the theory of quantum noise.  The processes $(A_t+A_t^*)_{t\ge 0}$ or 
$(\Lambda_t)_{t\ge 0}$ each generate a commutative algebra; the spectral 
theorem allows us to interpret the former as a diffusion process (a Wiener 
process in the vacuum state), and the latter as a counting process (a 
Poisson process in a so-called coherent state, and a.s.\ zero in the 
vacuum).  A brief introduction is given in e.g.\ \cite{BHJ06}.

Based on these processes, one can now proceed to define quantum stochastic
integrals and obtain the quantum It\^o rule: $dA_t\,dA_t^*=dt$,
$d\Lambda_t\,dA_t^*=dA_t^*$, $dA_t\,d\Lambda_t=dA_t$, and all other
combinations are zero.  Next one defines quantum stochastic differential
equations, which provide "quantum noisy" versions of the famous
Schr{\"o}dinger equation and form good models of actual physical systems,
particularly those used in quantum optics (see \cite{GaZ04} for a
physicist's perspective).  An accessible starting point in the literature
on quantum stochastic calculus are Hudson's lecture notes \cite{Hud03};
see also the original article by Hudson and Parthasarathy \cite{HuP84} and
Parthasarathy's book \cite{Par92}.  Several other approaches to the theory
have been developed since; some of these appear in the book by Meyer
\cite{Mey93} and the lecture notes of Biane \cite{Bia95}.

\subsection*{Physical models}

The physical theory that describes the interaction of light (the 
electromagnetic field) with matter (atoms, molecules, etc.)\ is called 
quantum electrodynamics; see the book by Cohen-Tannoudji \emph{et al.}\ 
\cite{CT89}.   Very often, and in particular for the purposes of filtering 
and control, quantum electrodynamics does not directly provide a usable 
model; these models are non-Markov.  In many systems, however, the Markov 
approximation is quite good, and applying this approximation to models 
from quantum electrodynamics gives rise precisely to quantum stochastic 
differential equations.  The Markov approximation can be pursued at 
different levels of rigor, ranging from the ad-hoc ``whitening'' of the 
noise which is common in the physics literature \cite{GaZ04} to rigorous 
Wong-Zakai type limits \cite{Gou05,AGLu95}.  As in classical probability 
theory it is also common to begin with a Markov model in the form of a 
quantum stochastic differential equation, which is justified either on 
phenomenological grounds or through a separate modelling effort.

The modelling of the detectors, e.g.\ photodetection or homodyne 
detection, is another important topic.  The operating principles of 
optical detectors are described in many textbooks on quantum optics, see 
e.g.\ \cite{WaM94,MaWo95}.  In the quantum stochastic context we refer to 
e.g.\ \cite{GaZ04}, or to the more mathematical perspective of Barchielli 
\cite{Bar03}.

\subsection*{Continuous limits of discrete models}

The discrete models which we have used throughout the paper are
conceptually very similar to the more realistic continuous models. There
is more than a conceptual resemblance, however: these models in fact
converge to corresponding continuous models as we let the time step
$\lambda^2\to 0$.  Though there is no particular physical intuition behind
this convergence except the one we have already given---that small time 
slices of the field approximately contain at most one photon---the 
convergence demonstrates that even results obtained using these highly 
simplified models have some relevance to the physical world.

A simple example of this convergence is easily verified.  Recall that
the classical stochastic process $x(l)=\iota(A(l)+A^*(l))$, with the 
measure induced by the vacuum state, defines a symmetric random walk with 
step size $\lambda$.  First, let us express this process in terms of time 
$t=l\lambda^2$ rather than the time step number $l$, i.e., $x_t=x(\lfloor 
t/\lambda^2\rfloor)$.  Taking the limit (in law) as $\lambda\to 0$ 
gives, by the usual functional central limit argument (e.g.\ \cite[pp.\ 
452--]{GiSk}), a Wiener process.  But $\iota(A_t+A_t^*)$ is a Wiener 
process, so we see that $A(l)+A^*(l)$ converges to $A_t+A_t^*$ at least in 
a weak sense.  

Similarly, most of the expressions given in this article converge as 
$\lambda\to 0$ to their continuous time counterparts.  In particular the 
difference equation Eq.\ (\ref{eq U}), the filtering equations, etc.\ 
converge to (quantum) stochastic differential equations.  Even the 
discrete It\^o table ``converges'' to the quantum It\^o rules 
$dA_t\,dA_t^*=dt$ etc., if we formally set $\Delta\Lambda\to d\Lambda$, 
$\Delta A\to dA$, $\Delta A^*\to dA^*$, $\Delta t\to dt$, and terms such 
as $\lambda^2\Delta A\to 0$ (this was noticed in \cite{Att03}).
Evidently the discrete theory mirrors its continuous counterpart quite 
faithfully, which is precisely the motivation for this article.

Convergence of the discrete models is investigated by Lindsay and 
Parthasarathy \cite{LP88}, Attal and Pautrat \cite{AtP06}, 
and Gough \cite{Gough04}.  These articles demonstrate convergence (of 
various types) of the solution of the difference equation Eq.\ (\ref{eq U}) 
to the solution of a corresponding quantum stochastic differential 
equation.  The convergence of discrete filters to their continuous 
counterparts is investigated in Gough and Sobolev \cite{GS04}.

\subsection*{Conditional expectations and quantum filtering}

The concept of conditional expectations in Von Neumann algebras is an old
one, see e.g.\ \cite{Tak71}.  However, a much more pragmatic notion of
conditional expectations on which the definition we have given is based
first appeared in the pioneering work of Belavkin, see e.g.\
\cite{Bel92b}.  The necessity of the nondemolition property, which is the
key to the development of both filtering and feedback control, appeared
already in \cite{Bel80}.  This opens the door for the development of
nonlinear filtering, which is done in detail in the difficult and
technical paper \cite{Bel92b} using martingale methods.  The reference
probability approach to quantum filtering, using the Bayes formula,
appears in \cite{BH05} (but see also \cite{Bel92a,GG94}).  An introduction 
to quantum filtering is given in \cite{BHJ06} using both methods.

\subsection*{Feedback control}

An early investigation on optimal feedback control in discrete time can be
found in Belavkin's 1983 paper \cite{Bel83}.  The Bellman equation for
continuous optimal control appears in \cite{Bel88}.  A recent exposition
on quantum optimal control, and in particular quantum LQG control, is
given in \cite{EdB05}, see also \cite{BEB05}, and see \cite{DJ99,DHJMT00} 
for a physicist's perspective.  The separation theorem for the continuous 
time case appears in \cite{BH05}.  A different type of optimal control 
problem, the risk-sensitive control problem with exponential cost, is 
investigated in \cite{Jam05}.  Finally, Lyapunov function methods 
for control design are investigated in \cite{HSM05,mazyar}.

Many applications of quantum feedback control have been considered in the
physics literature; we refer to the introduction for references.  Note 
that the ``quantum trajectory equations'' or ``stochastic master 
equations'' used in the physics literature are precisely filtering 
equations, written in terms of the innovations process as the driving 
noise rather than being driven by the observations.  The notion of 
filtering has only recently gained acceptance in the physics literature, 
though the relevant equations had already been developed by physicists.  
The reader should beware of the discrepancy in terminology and 
interpretation between the physical and mathematical literature.  However, 
as the Markov property of the filter allows us to pretend that the filter 
itself is the system to be controlled, most of the control schemes found 
in the physics literature make sense in both contexts.

\section*{Acknowledgment}
The authors would like to thank John Stockton for providing a figure.  
L.B.\ thanks Hans Maassen for introducing him to discrete systems.

\bibliographystyle{siam}
\bibliography{ref}           

\begin{thebibliography}{10}

\bibitem{AGLu95}
{\sc L.~Accardi, J.~Gough, and Y.~Lu}, {\em On the stochastic limit for quantum
  theory}, Rep. Math. Phys., 36 (1995), pp.~155--187.

\bibitem{ASL04}
{\sc A.~Andr{\'e}, A.~S. S{\o}rensen, and M.~D. Lukin}, {\em Stability of
  atomic clocks based on entangled atoms}, Phys. Rev. Lett., 92 (2004),
  p.~230801.

\bibitem{AASDM02}
{\sc M.~A. Armen, J.~K. Au, J.~K. Stockton, A.~C. Doherty, and H.~Mabuchi},
  {\em Adaptive homodyne measurement of optical phase}, Phys.~Rev.~Lett., 89
  (2002), p.~133602.

\bibitem{KA70}
{\sc K.~Astrom}, {\em Introduction to Stochastic Control Theory}, Academic
  Press, New York, 1970.

\bibitem{Att03}
{\sc S.~Attal}, {\em Approximating the {F}ock space with the toy {F}ock space},
  in S\'eminaire de Probabilit\'es, XXXVI, vol.~1801 of Lecture Notes in Math.,
  Springer, Berlin, 2003, pp.~477--491.

\bibitem{AtP06}
{\sc S.~Attal and Y.~Pautrat}, {\em From repeated to continuous quantum
  interactions}, Ann. Henri Poincar\'e, 7 (2006), pp.~59--104.

\bibitem{Bar03}
{\sc A.~Barchielli}, {\em Continual measurements in quantum mechanics and
  quantum stochastic calculus}, in Open Quantum Systems III: Recent
  Developments, S.~Attal, A.~Joye, and C.-A. Pillet, eds., Springer, 2006,
  pp.~207--292.

\bibitem{Bel80}
{\sc V.~P. Belavkin}, {\em Quantum filtering of {Markov} signals with white
  quantum noise}, Radiotechnika i Electronika, 25 (1980), pp.~1445--1453.

\bibitem{Bel83}
\leavevmode\vrule height 2pt depth -1.6pt width 23pt, {\em Theory of the
  control of observable quantum systems}, Autom. Rem. Control, 44 (1983),
  pp.~178--188.

\bibitem{Bel88}
\leavevmode\vrule height 2pt depth -1.6pt width 23pt, {\em Nondemolition
  stochastic calculus in {Fock} space and nonlinear filtering and control in
  quantum systems}, in Proceedings {XXIV} Karpacz winter school, R.~Guelerak
  and W.~Karwowski, eds., Stochastic methods in mathematics and physics, World
  Scientific, Singapore, 1988, pp.~310--324.

\bibitem{Bel92a}
\leavevmode\vrule height 2pt depth -1.6pt width 23pt, {\em Quantum continual
  measurements and a posteriori collapse on {CCR}}, Commun. Math. Phys., 146
  (1992), pp.~611--635.

\bibitem{Bel92b}
\leavevmode\vrule height 2pt depth -1.6pt width 23pt, {\em Quantum stochastic
  calculus and quantum nonlinear filtering}, J. Multivar. Anal., 42 (1992),
  pp.~171--201.

\bibitem{Bensoussan}
{\sc A.~Bensoussan}, {\em Stochastic Control of Partially Observable Systems},
  Cambridge University Press, 1992.

\bibitem{DB95}
{\sc D.~Bertsekas}, {\em Dynamic Programming and Optimal Control}, vol.~1,
  Athena Scientific, Boston, 1995.

\bibitem{BS78}
{\sc D.~Bertsekas and S.~Shreve}, {\em Stochastic Optimal Control: The
  Discrete-Time Case}, Academic Press, New York, 1978.

\bibitem{Bia95}
{\sc P.~Biane}, {\em Calcul stochastique noncommutatif}, in Lectures on
  Probability Theory. (Saint-Flour, 1993), P.~Bernard, ed., vol.~1608 of
  Lecture Notes in Mathematics, Springer, Berlin, 1995, pp.~1--96.

\bibitem{BEB05}
{\sc L.~M. Bouten, S.~C. Edwards, and V.~P. Belavkin}, {\em Bellman equations
  for optimal feedback control of qubit states}, J. Phys. B, At. Mol. Opt.
  Phys., 38 (2005), pp.~151--160.

\bibitem{BGM04}
{\sc L.~M. Bouten, M.~I. Gu\c{t}\u{a}, and H.~Maassen}, {\em Stochastic
  {Schr\"odinger} equations}, J. Phys. A: Math. Gen., 37 (2004),
  pp.~3189--3209.

\bibitem{BH05}
{\sc L.~M. Bouten and R.~{Van Handel}}, {\em Controlled quantum stochastic
  processes}, 2006.
\newblock In preparation; see ``On the separation principle of quantum
  control'', math-ph/0511021, and ``Quantum filtering: a reference probability
  approach'', math-ph/0508006.

\bibitem{BHJ06}
{\sc L.~M. Bouten, R.~{Van Handel}, and M.~R. James}, {\em An introduction to
  quantum filtering}, preprint, {\tt http://arxiv.org/math.OC/0601741},
  (2006).

\bibitem{BrR87}
{\sc O.~Bratteli and D.~Robinson}, {\em Operator algebras and quantum
  statistical mechanics 1}, Springer-Verlag, Berlin Heidelberg, second~ed.,
  1987.

\bibitem{Brun02}
{\sc T.~A. Brun}, {\em A simple model of quantum trajectories}, Am. J. Phys.,
  70 (2002), pp.~719--737.

\bibitem{ExpCool06}
{\sc P.~Bushev, D.~Rotter, A.~Wilson, F.~Dubin, C.~Becher, J.~Eschner,
  R.~Blatt, V.~Steixner, P.~Rabl, and P.~Zoller}, {\em Feedback cooling of a
  single trapped ion}, Phys. Rev. Lett., 96 (2006), p.~043003.

\bibitem{Car93}
{\sc H.~J. Carmichael}, {\em An Open Systems Approach to Quantum Optics},
  Springer-Verlag, Berlin Heidelberg New-York, 1993.

\bibitem{CT89}
{\sc C.~{C}ohen {T}annoudji, J.~{D}upont {R}oc, and G.~Grynberg}, {\em Photons
  and Atoms: Introduction to Quantum Electrodynamics}, Wiley, 1989.

\bibitem{CRR79}
{\sc J.~C. Cox, S.~A. Ross, and M.~Rubinstein}, {\em Option pricing: A
  simplified approach}, Journal of Financial Economics, 7 (1979), pp.~229--263.

\bibitem{Dav76}
{\sc E.~B. Davies}, {\em Quantum Theory of Open Systems}, Academic Press,
  London New-York San Francisco, 1976.

\bibitem{Dav77}
\leavevmode\vrule height 2pt depth -1.6pt width 23pt, {\em Quantum
  communication systems}, IEEE Trans. Inf. Th., IT-23 (1977), pp.~530--534.

\bibitem{DM81}
{\sc M.~H.~A. Davis and S.~I. Marcus}, {\em An introduction to nonlinear
  filtering}, in Stochastic Systems: The Mathematics of Filtering and
  Identification and Applications, M.~Hazewinkel and J.~C. Willems, eds., D.
  Reidel, 1981, pp.~53--75.

\bibitem{DHJMT00}
{\sc A.~C. Doherty, S.~Habib, K.~Jacobs, H.~Mabuchi, and S.~M. Tan}, {\em
  Quantum feedback and classical control theory}, Phys. Rev. A, 62 (2000),
  p.~012105.

\bibitem{DJ99}
{\sc A.~C. Doherty and K.~Jacobs}, {\em Feedback-control of quantum systems
  using continuous state-estimation}, Phys. Rev. A, 60 (1999), pp.~2700--2711.

\bibitem{TD68}
{\sc T.~Duncan}, {\em Evaluation of likelihood functions}, Information and
  Control,  (1968), pp.~62--74.

\bibitem{EdB05}
{\sc S.~C. Edwards and V.~P. Belavkin}, {\em Optimal quantum feedback control
  via quantum dynamic programming}, quant-ph/0506018, University of Nottingham,
  2005.

\bibitem{AEM95}
{\sc R.~J. Elliott, L.~Aggoun, and J.~B. Moore}, {\em Hidden {M}arkov Models:
  Estimation and Control}, Springer, New York, 1995.

\bibitem{FR75}
{\sc W.~Fleming and R.~Rishel}, {\em Deterministic and Stochastic Optimal
  Control}, Springer Verlag, New York, 1975.

\bibitem{GaZ04}
{\sc C.~Gardiner and P.~Zoller}, {\em Quantum Noise}, Springer, third~ed.,
  2004.

\bibitem{Ger04}
{\sc J.~M. Geremia}, {\em Distinguishing between optical coherent states with
  imperfect detection}, Phys. Rev. A, 70 (2004), p.~062303.

\bibitem{GSDM03}
{\sc J.~M. Geremia, J.~K. Stockton, A.~C. Doherty, and H.~Mabuchi}, {\em
  Quantum {Kalman} filtering and the {Heisenberg} limit in atomic
  magnetometry}, Phys. Rev. Lett., 91 (2003), p.~250801.

\bibitem{GSM05}
{\sc J.~M. Geremia, J.~K. Stockton, and H.~Mabuchi}, {\em Suppression of spin
  projection noise in broadband atomic magnetometry}, Phys. Rev. Lett., 94
  (2005), p.~203002.

\bibitem{GiSk}
{\sc I.~I. Gikhman and A.~V. Skorokhod}, {\em Introduction to the theory of
  random processes}, Dover, 1996.

\bibitem{GG94}
{\sc P.~Goetsch and R.~Graham}, {\em Linear stochastic wave equations for
  continuously measured quantum systems}, Phys. Rev. A, 50 (1994),
  pp.~5242--5255.

\bibitem{Gough04}
{\sc J.~Gough}, {\em Holevo-ordering and the continuous-time limit for open
  {F}loquet dynamics}, Lett. Math. Phys., 67 (2004), pp.~207--221.

\bibitem{Gou05}
\leavevmode\vrule height 2pt depth -1.6pt width 23pt, {\em Quantum flows as
  {Markovian} limit of emission, absorption and scattering interactions},
  Commun. Math. Phys., 254 (2005), pp.~489--512.

\bibitem{GS04}
{\sc J.~Gough and A.~Sobolev}, {\em Stochastic {S}chr\"odinger equations as
  limit of discrete filtering}, Open Syst. Inf. Dyn., 11 (2004), pp.~235--255.

\bibitem{PRH74}
{\sc P.~Halmos}, {\em Finite-Dimensional Vector Spaces}, Springer Verlag, New
  York, 1974.

\bibitem{Hol90}
{\sc A.~Holevo}, {\em Quantum stochastic calculus}, J. Soviet Math., 56 (1991),
  pp.~2609--2624.
\newblock Translation of Itogi Nauki i Tekhniki, ser. sovr. prob. mat. 36,
  3--28, 1990.

\bibitem{Asa03}
{\sc A.~Hopkins, K.~Jacobs, S.~Habib, and K.~Schwab}, {\em Feedback cooling of
  a nanomechanical resonator}, Phys. Rev. B, 68 (2003), p.~235328.

\bibitem{HJ85}
{\sc R.~A. Horn and C.~R. Johnson}, {\em Matrix analysis}, Cambridge University
  Press, 1985.

\bibitem{Hud03}
{\sc R.~L. Hudson}, {\em An introduction to quantum stochastic calculus and
  some of its applications}, in Quantum Probability Communications, S.~Attal
  and J.~Lindsay, eds., vol.~XI, World Scientific, Singapore, 2003,
  pp.~221--271.

\bibitem{HuP84}
{\sc R.~L. Hudson and K.~R. Parthasarathy}, {\em Quantum {It\^o's} formula and
  stochastic evolutions}, Commun. Math. Phys., 93 (1984), pp.~301--323.

\bibitem{JP05}
{\sc M.~James and I.~Petersen}, {\em Robustness properties of a class of
  optimal risk-sensitive controllers for quantum systems}, in Proc. 16th IFAC
  World Congress, 2005.

\bibitem{JP06a}
\leavevmode\vrule height 2pt depth -1.6pt width 23pt, {\em ${H}^\infty$ control
  of linear quantum systems}, in Proc. IEEE CDC, 2006.

\bibitem{Jam05}
{\sc M.~R. James}, {\em A quantum {Langevin} formulation of risk-sensitive
  optimal control}, J. Opt. B: Quantum Semiclass. Opt., 7 (2005),
  pp.~S198--S207.

\bibitem{KaR83}
{\sc R.~V. Kadison and J.~R. Ringrose}, {\em Fundamentals of the Theory of
  Operator Algebras}, vol.~I, Academic Press, San Diego, 1983.

\bibitem{KaR86}
\leavevmode\vrule height 2pt depth -1.6pt width 23pt, {\em Fundamentals of the
  Theory of Operator Algebras}, vol.~II, Academic Press, San Diego, 1986.

\bibitem{Kal80}
{\sc G.~Kallianpur}, {\em Stochastic Filtering Theory}, Springer, Berlin, 1980.

\bibitem{KalStr68}
{\sc G.~Kallianpur and C.~Striebel}, {\em Estimation of stochastic systems:
  {A}rbitrary system process with additive white noise observation errors},
  Ann. Math. Statist., 39 (1968), pp.~785--801.

\bibitem{KalStr69}
\leavevmode\vrule height 2pt depth -1.6pt width 23pt, {\em Stochastic
  differential equations occurring in the estimation of continuous parameter
  stochastic processes}, Teor. Verojatnost. i Primenen, 14 (1969),
  pp.~597--622.

\bibitem{navin}
{\sc N.~Khaneja, T.~Reiss, B.~Luy, and S.~J. Glaser}, {\em Optimal control of
  spin dynamics in the presence of relaxation}, J. Magnet. Res., 162 (2003),
  pp.~311--319.

\bibitem{Kr05}
{\sc V.~Krishnan}, {\em Nonlinear filtering and smoothing: An Introduction to
  Martingales, Stochastic Integrals and Estimation}, Dover, 2005.

\bibitem{KV86}
{\sc P.~Kumar and P.~Varaiya}, {\em Stochastic Systems: Estimation,
  Identification and Adaptive Control}, Prentice-Hall, Englewood Cliffs, NJ,
  1986.

\bibitem{Kum85}
{\sc B.~K\"ummerer}, {\em Markov dilations on {$W^*$-algebras}}, J. Funct.
  Anal., 63 (1985), pp.~139--177.

\bibitem{Kus71}
{\sc H.~Kushner}, {\em Introduction to stochastic control}, Holt, Rinehart and
  Winston, Inc, New York, 1971.

\bibitem{LP88}
{\sc J.~M. Lindsay and K.~R. Parthasarathy}, {\em The passage from random walk
  to diffusion in quantum probability. {II}}, Sankhy\=a Ser. A, 50 (1988),
  pp.~151--170.

\bibitem{LiS01}
{\sc R.~S. Liptser and A.~N. Shiryaev}, {\em Statistics of Random Processes I:
  General Theory}, Springer-Verlag, 2001.

\bibitem{Maa03}
{\sc H.~Maassen}, {\em Quantum probability applied to the damped harmonic
  oscillator}, in Quantum Probability Communications, S.~Attal and J.~Lindsay,
  eds., vol.~XII, World Scientific, Singapore, 2003, pp.~23--58.

\bibitem{MaWo95}
{\sc L.~Mandel and E.~Wolf}, {\em Optical coherence and quantum optics},
  Cambridge University Press, 1995.

\bibitem{EM98}
{\sc E.~Merzbacher}, {\em Quantum Mechanics}, Wiley, New York, third~ed., 1998.

\bibitem{Mey93}
{\sc P.-A. Meyer}, {\em Quantum Probability for Probabilists}, Springer,
  Berlin, 1993.

\bibitem{mazyar}
{\sc M.~Mirrahimi and R.~{Van Handel}}, {\em Stabilizing feedback controls for
  quantum systems}, 2005.
\newblock Submitted, http://arxiv.org/abs/math-ph/0510066.

\bibitem{RM66}
{\sc R.~Mortensen}, {\em Optimal Control of Continuous-Time Stochastic
  Systems}, PhD thesis, Univ. California, Berkeley, 1966.

\bibitem{Par92}
{\sc K.~R. Parthasarathy}, {\em An Introduction to Quantum Stochastic
  Calculus}, Birkh\"auser, Basel, 1992.

\bibitem{Pro04}
{\sc P.~E. Protter}, {\em Stochastic Integration and Differential Equations},
  Springer-Verlag, second~ed., 2004.

\bibitem{ReS80}
{\sc M.~Reed and B.~Simon}, {\em Functional Analysis}, vol.~1 of Methods of
  Modern Mathematical Physics, Elsevier, 1980.

\bibitem{Shr04}
{\sc S.~E. Shreve}, {\em Stochastic calculus for finance {I}. The binomial
  asset pricing model}, Springer-Verlag, New York, 2004.

\bibitem{Steck04}
{\sc D.~A. Steck, K.~Jacobs, H.~Mabuchi, T.~Bhattacharya, and S.~Habib}, {\em
  Quantum feedback control of atomic motion in an optical cavity}, Phys. Rev.
  Lett., 92 (2004), p.~223004.

\bibitem{SGDM04}
{\sc J.~K. Stockton, J.~M. Geremia, A.~C. Doherty, and H.~Mabuchi}, {\em Robust
  quantum parameter estimation: Coherent magnetometry with feedback}, Phys.
  Rev. A, 69 (2004), p.~032109.

\bibitem{Tak71}
{\sc M.~Takesaki}, {\em Conditional expectations in von {Neumann} algebras}, J.
  Funct. Anal., 9 (1971), pp.~306--321.

\bibitem{HSM05}
{\sc R.~{Van Handel}, J.~K. Stockton, and H.~Mabuchi}, {\em Feedback control of
  quantum state reduction}, IEEE Trans. Automat. Control, 50 (2005),
  pp.~768--780.

\bibitem{WaM94}
{\sc D.~Walls and G.~Milburn}, {\em Quantum Optics}, Springer Verlag, Berlin
  Heidelberg, 1994.

\bibitem{W81}
{\sc P.~Whittle}, {\em Risk-sensitive linear/quadratic/{G}aussian control},
  Advances in Applied Probability, 13 (1981), pp.~764--777.

\bibitem{Wil91}
{\sc D.~Williams}, {\em Probability with Martingales}, Cambridge University
  Press, Cambridge, 1991.

\bibitem{wisemanadaptive}
{\sc H.~M. Wiseman}, {\em Adaptive phase measurements of optical modes: Going
  beyond the marginal $q$ distribution}, Phys. Rev. Lett., 75 (1995),
  pp.~4587--4590.

\bibitem{WM94b}
{\sc H.~M. Wiseman and G.~J. Milburn}, {\em All-optical versus electro-optical
  quantum-limited feedback}, Phys. Rev. A, 49 (1994), pp.~4110--4125.

\bibitem{naoki}
{\sc N.~Yamamoto, K.~Tsumura, and S.~Hara}, {\em Feedback control of quantum
  entanglement in a two-spin system}, Proc. 44th IEEE CDC,  (2005),
  pp.~3182--3187.

\bibitem{YK03}
{\sc M.~Yanagisawa and H.~Kimura}, {\em Transfer function approach to quantum
  control---{P}art {II}: Control concepts and applications}, IEEE Trans.
  Automat. Control, 48 (2003), pp.~2121--2132.

\bibitem{Zak69}
{\sc M.~Zakai}, {\em On the optimal filtering of diffusion processes}, Z.
  Wahrsch. Verw. Geb., 11 (1969), pp.~230--243.

\end{thebibliography}
\end{document}